\documentclass[a4paper, a4wide, 10pt, reqno]{amsart}

\usepackage{a4wide}
\usepackage{amsfonts}
\usepackage{amsmath}
\usepackage{amssymb}
\usepackage{amsthm}
\usepackage{bm}
\usepackage{dsfont}
\usepackage{enumitem}
\usepackage{hyperref}
\usepackage[capitalise]{cleveref}
\usepackage[utf8]{inputenc}
\usepackage{mathtools}
\usepackage{nccmath}
\usepackage{relsize}
\usepackage{subcaption}
\theoremstyle{plain}
\newtheorem{theorem}{Theorem}[section]
\newtheorem{metatheorem}[theorem]{Meta-theorem}
\newtheorem{corollary}[theorem]{Corollary}

\newtheorem{proposition}[theorem]{Proposition}

\newtheorem{lemma}[theorem]{Lemma}

\newtheorem{definition}[theorem]{Definition}

\theoremstyle{definition}
\newtheorem{remark}[theorem]{Remark}
\numberwithin{equation}{section}

\newcommand{\R}{{\mathbb R}}

\newcommand{\N}{{\mathbb N}}
\newcommand{\E}{\mathbb E}
\newcommand{\Prob}{\mathbb{P}}
\newcommand\eps{\varepsilon}
\newcommand{\CA}{{\mathcal A}}
\newcommand{\CB}{{\mathcal B}}
\newcommand{\CC}{{\mathcal C}}

\newcommand{\CE}{{\mathcal E}}

\newcommand{\CG}{{\mathcal G}}

\newcommand{\CM}{{\mathcal M}}

\newcommand{\CV}{{\mathcal V}}

\newcommand{\CY}{{\mathcal Y}}

\newcommand{\re}{{\mathrm e}}

\newcommand{\ind}[1]{\mathds{1}_{\{#1\}}}
\newcommand{\rd}{\mathrm{d}}
\newcommand{\sss}[1]{{\scriptscriptstyle #1}}

\title[Large deviations of the giant in scale-free random graphs]{Large deviations of the giant component in\\ scale-free inhomogeneous random graphs}
\author[J.\ Jorritsma]{Joost Jorritsma$^{1,2}$}
\address{$^1$Centre for Mathematics \& Computer Science (CWI), The Netherlands}

\author[B. Zwart]{Bert Zwart$^{1,2}$}
\address{$^2$Department of Mathematics \& Computer Science, Eindhoven University of Technology}

\begin{document}
\begin{abstract}
    We study large deviations of the size of the largest connected component in a general class of inhomogeneous random graphs with iid weights, parametrized so that the degree distribution is regularly varying. We derive a large-deviation principle with logarithmic speed: the rare event that the largest component contains linearly more vertices than expected is caused by the presence of constantly many vertices with linear degree.  Conditionally on this rare event, we prove distributional limits of the  weight distribution and component-size distribution. 
\end{abstract}
\maketitle
\vspace*{-0.65cm}
{\footnotesize
\hspace{1em}Keywords: Large-deviation principle, giant component, inhomogeneous random graph, scale-free network.

\hspace{1em}MSC Class: 05C80, 60F10.
}
\vspace*{-0.2cm}
\section{Introduction}
\vspace*{-0.05cm}
This paper studies large deviations of the size of the largest connected component in a class of inhomogeneous random graphs (IRGs) in which the degree distribution is regularly varying. IRGs were introduced in a seminal paper by Bollob\'as, Janson, and Riordan~\cite{bollobas2007irg}. In short, each vertex $u\in[n]:=\{1,\ldots, n\}$ has a mark $m_u$ from a mark set $\CM$, such that the empirical mark-distribution converges weakly. Typically, marks are either independent and identically distributed (iid), or deterministic for each $n$. Given the marks, each pair of vertices $uv$ is independently connected by an edge  with probability $\min\big(\kappa(m_u, m_v)/n, 1\big)$, where $\kappa$ is the so-called kernel function that encodes the influence of marks: it is non-negative, and symmetric in its arguments. 

It is well-known that the size of the largest connected component $|\CC_n^\sss{(1)}|$ satisfies a weak law of large numbers as $n\to\infty$. The proportion of vertices in $\CC_n^\sss{(1)}$ converges  in probability to $\theta$: the probability that the branching process describing the local limit survives infinitely long, see~\cite{bollobas2007irg, Hofbook2} for details. When this branching process is supercritical, the graph contains a \emph{giant},  a linear-sized component, with probability tending to one as $n\to\infty$ (with high probability; whp).  As its asymptotic size is determined by the local limit, the size of the giant is ``almost local''\cite{hofstad2021giantlocal}.
Under regularity conditions, the giant (if it exists) is unique, and all other components are at most of logarithmic size whp~\cite{bollobas2007irg,ChungLu02.1, NorRei06}.

A large-deviation principle (LDP) of $|\CC_n^\sss{(1)}|/n$ (among other results) has been established recently by Andreis \emph{et al.}~\cite{andreis2023irg} for a subclass of IRGs: they consider \emph{bounded} kernels, and  \emph{deterministic} marks. This setting leads to a graph in which the degree distribution has an \emph{exponential} tail.  They identify a non-negative rate function $I_\mathrm{det}$ such that for each Borel set $B\subseteq\R$ (writing $\bar B$ for its closure and $B^\circ$ for its interior), 
    \begin{equation}\label{eq:ldp-intro-exp}
    -\!\! \inf_{\rho\in B^\circ}I_\mathrm{det}(\rho)\le \liminf_{n\to\infty}\frac{1}n\log \Prob\bigg(\frac{|\CC_{n}^\sss{(1)}|}n\in B\bigg)
    \le \limsup_{n\to\infty}\frac{1}{n}\log \Prob\bigg(\frac{|\CC_{n}^\sss{(1)}|}n \in B\bigg)\le -\!\!\inf_{\rho \in \bar B} I_\mathrm{det}(\rho).
    \end{equation}
Thus,~\cite{andreis2023irg} proves an LDP of the giant with \emph{linear} speed. Moreover, the rate-function $ I_\mathrm{det}$ is determined uniquely by the local limit of the graph. The present paper proves an LDP for $|\CC_n^\sss{(1)}|/n$ in a complementary setting, in which the marks are \emph{iid}, and the kernel is $\emph{unbounded}$. In our setting, the mark distribution and kernel lead to a degree distribution that has a  \emph{regularly varying} tail, similar to many real-world networks~\cite{New03,voitalov2019scale}. We summarize our main results informally. 
\begin{metatheorem}[Logarithmic speed and non-local rate function]
    Consider an inhomogeneous random graph with iid marks and a regularly-varying degree distribution. Then 
    \begin{enumerate}
        \item $|\CC_n^\sss{(1)}|/n$ satisfies an LDP with logarithmic speed and rate function $I_\mathrm{iid}(\cdot)$. 
        \item the rate function of this LDP may be different for two inhomogeneous random graphs with the same local limit.
        \item conditionally on the event $\{|\CC_n^\sss{(1)}|>\rho n\}$, for any $\rho\in(\theta, 1)$, the component-size distribution converges to a random measure as $n\to\infty$,  and $|\CC_n^\sss{(1)}|/n$ converges to a random variable.
    \end{enumerate}
\end{metatheorem}
\subsection*{Notation}
If two vertices $u$ and $v$ are connected by an edge, we write $u\sim v$. Similarly, for sets of vertices $A$ and $B$ we write $u\sim B$ (resp.\  $A\sim B$) if $u$ is connected to a vertex  $v\in B$ (resp.\ there exists $u\in A$ and $v\in B$ such that $u\sim v$). We use standard Bachmann--Landau notation: we say that $f=o(g)$ if $|f(x)/g(x)|\to 0$, $f=\omega(g)$ if $g=o(f)$, $f=O(g)$ if $\limsup_{x\to\infty}|f(x)/g(x)|<\infty$, $f=\Omega(g)$ if $g=O(f)$, and $f=\Theta(g)$ if both $f=O(g)$ and $f=\Omega(g)$. Similarly, we write $f\sim g$ if $f(x)/g(x)\to1$, $f\lesssim g$ if $\limsup f(x)/g(x)\le 1$, and $f\gtrsim g$ if $\liminf f(x)/g(x)\ge 1$. We say that a random variable $X$ stochastically dominates $Y$ if $\Prob(X\ge x)\ge \Prob(Y\ge x)$ for all $x\in\R$, in which case we write $X\succcurlyeq Y$ or $Y\preccurlyeq X$. We abbreviate $a\vee b:=\max(a,b)$, $a\wedge b:=\min(a,b)$, $[n]:=\{1,\ldots, n\}$, and write  $A^\ell$ for the $\ell$-fold Cartesian product of a set $A$.

\subsection{Model description} We formalize the model.  We call a function $L(x)$ slowly varying if $L(cx)/L(x)\to 1$ as $x\to\infty$ for any constant $c>0$, while we say that $F(x)$ is regularly varying with index $-\alpha$ if there exists a slowly varying function $L(x)$ such that $F(x)=L(x)x^{-\alpha}$. 
\begin{definition}[Inhomogeneous scale-free random graph]\label{def:irg}
Let $\alpha>1$, $\sigma<2\alpha-1$, and $q\in(0,1]$ be three constants.
    Consider the vertex set $\CV_n:=[n]:=\{1,\ldots,n\}$, and equip each vertex $u\in[n]$ with a weight $W_u$, which is an iid copy of the non-negative random variable $W$ whose distribution is given by 
    \begin{equation}
    1- F_W(w):=\Prob\big(W>w\big)=L(w)w^{-\alpha}, \qquad w\ge \underline w,\label{eq:weight-dist}
    \end{equation}
    for some slowly varying function $L(w)$ and constant $\underline w:=\inf_x\{x: F_W(x)>0\}>0$. Conditionally on all weights, two vertices $u$, $v$ are independently connected by an edge in $\CG_n=(\CV_n,\CE_n)$ with probability 
    \begin{equation}
    p_{uv}:=q\bigg(\frac{\kappa_\sigma(W_u, W_v)}{n}\wedge  1\bigg):=q\bigg(\frac{(W_u\vee W_v)(W_u\wedge W_v)^\sigma}{n}\wedge  1\bigg),\label{eq:conn-prob}
    \end{equation}
    We denote the largest connected component in $\CG_{n}$ by $\CC_{n}^\sss{(1)}$, and the component of vertex $u$ by $\CC_{n}(u)$. 
\end{definition}
Let us comment on the setup. Our assumption that $\underline w>0$ with strict inequality ensures that $\lim_{w\downarrow \underline w}\kappa_\sigma(w,z)<\infty$ for all $z$ when $\sigma<0$. Moreover, it simplifies the technicalities when $\sigma\ge 0$, and we expect that all results carry through when one allows $\underline w=0$ for non-negative $\sigma$. 
The kernel $\kappa_\sigma$ was recently introduced for related random graph models that are embedded in  Euclidean space with the additional restriction that $\sigma\ge 0$, see e.g.~\cite{gracar2019age, clusterI}. 
The more general setting here (including $\sigma<0$) includes several models of interest: when $\sigma=1$, the kernel $\kappa_\sigma(W_u, W_v)$ is simply a product, and we obtain a rank-one inhomogeneous random graph~\cite{ChungLu02.1, NorRei06}. When $\sigma=\alpha-1$, the graph is closely related to preferential attachment models~\cite{dereich2012typical,gracar2019age,jacobLinkerMortersFastEvolving}. The parameter $\sigma$ can be seen as an (dis)\emph{assortativity} parameter that tunes the correlation between the degrees of vertices incident to the same edges. Small or negative values of $\sigma$ decrease the probability that high-weight vertices are connected by an edge, without necessarily affecting the tail of the degree distribution.
Indeed, the parameters $\sigma, \alpha$ jointly determine the tail of the degree distribution of the resulting graph: it is regularly varying with index $-\alpha/\max(1, 1+\sigma-\alpha)$, see e.g.\ \cite{luchtrathThesis2022} on a related model. Moreover, when $\sigma$ is negative, the set of edges is no longer monotone in the weights under a natural coupling that encodes the presence of edges using uniform random variables. Hence, also the set of vertices in the giant is no longer monotone in the weights when $\sigma<0$. Throughout the paper we assume that $\sigma<2\alpha-1$ and $\alpha>1$, which ensures that the asymptotic degree distribution has finite first moment. At the expense of additional technicalities, the results can be extended to $\sigma\ge 2\alpha-1$. Lastly, the parameter $q$ makes the model closed under edge percolation.

The following multi-type branching process describes the local limit of the graph, i.e., it describes the graph structure around a uniform vertex in $\CG_n$ (see~\cite{Hofbook2} for an introduction and references on local convergence).
\begin{definition}[Associated multi-type branching process]\label{def:branching}Consider an inhomogeneous scale-free random graph $\CG_n$ with kernel $\kappa_\sigma$, weight-distribution $F_W$, and percolation parameter $q$ as in \cref{def:irg}. The associated multi-type branching process $\mathrm{BP}=\mathrm{BP}(\kappa_\sigma, F_W, q)$ is a branching process that starts with a single vertex $\varnothing$ with random type $w_\varnothing$ distributed according to $F_W$. In each generation, each particle $v$ of type $w_v$ gives independently birth to new particles according to a Poisson point process on $[\underline w,\infty)$ with intensity $q\kappa_\sigma(w_v, w)\rd F_W(w)$.  The atoms in the union of these Poisson point processes form the vertex types of the vertices in the next generation. We write $T_q$ for the set of vertices in the progeny of the branching process, and define $\theta_q:=\Prob\big(|T_q|=\infty\big)$. 
\end{definition}

By \cite[Theorems 3.1 and 9.1]{bollobas2007irg}, the size of the largest connected component in $\CG_n$ satisfies a Law of Large Numbers (LLN). By local convergence, also the number of vertices in components of size $\ell\in\N$ satisfy an LLN, see for instance~\cite{hofstad2021giantlocal, Hofbook2} and~\eqref{eq:local-l2} below. Combined, we have as $n\to\infty$
\begin{equation}\label{eq:lln-intro}
    |\CC_n^\sss{(1)}|\,\big/\, n\overset{\Prob}\longrightarrow \theta_q, \qquad 
    S_{n,\ell}\,\big /\, n:=\big|\big\{u: |\CC_n(u)|=\ell\big\}\big|\,\big/\,n\overset{\Prob}\longrightarrow \Prob\big(|T_q|=\ell\big).
\end{equation}
Thus, the limits are uniquely determined by the local branching process.
Our large-deviation principle that we formalize in the following section, implies a convergence rate for the first LLN. Contrary to the LLN, we will see shortly that the convergence rate is not uniquely determined by the local limit. In contrast with the exponential decay in~\eqref{eq:ldp-intro-exp} for IRGs with bounded and deterministic weights, the convergence rate of the LLN is polynomial when the weights are iid and regularly varying. 
We explain the reason for this polynomial decay: the most likely way to have a large $\CC_n^\sss{(1)}$ (resp.\ small $S_{n,\ell}$) is when vertices of weight $\Theta(n)$ ---called hubs--- are present in the graph. Small components in the induced graph on non-hub vertices connect with constant probability by an edge to the hubs, increasing the size of the giant (that connect to the hubs with sufficiently high probability), and decreasing the number of small components. Since the weights are regularly varying, the probability of having $h$ hubs is of the order $(n\Prob(W>n))^h=(L(n)n^{1-\alpha})^h$, i.e., the decay is much slower than exponential. To motivate notation below, we next present a back-of-the-envelope calculation for the required number of hubs that realize the event $\{|\CC_n^\sss{(1)}|>\rho n\}$.
\smallskip

\noindent
\emph{Back-of-the-envelope calculation for the number of hubs.}\phantomsection{\label{back-of-the-envelope}}
Let $\CV_n[a,b)$ denote the set of vertices with weight in the interval $[a,b)$, and let $\CC_n(v)[a,b)$ denote the connected component of vertex $v$ in the induced subgraph on $\CV_n[a,b)$; write $A\sim B$ if there exists an edge between two sets of vertices $A$ and $B$, and $A\nsim B$ otherwise.  Instead of analyzing $\CC_n^\sss{(1)}$, we analyze the size of its complement $n-|\CC_n^\sss{(1)}|$. 
Assuming that the hubs (for this computation, these are the vertices with weight at least $\underline w^{-\sigma} n$) are part of the giant,  small components of $\CG_n[\underline w,\underline w^{-\sigma} n)$ do not merge with the giant in $\CG_n$ if there is no edge between the component and the hubs.  So, 
\begin{equation}\label{eq:random-clusters}
n-|\CC_n^\sss{(1)}|\approx \hspace{-10pt}\sum_{v\in\CV_n[1,\underline w^{-\sigma} n)}\hspace{-15pt}\ind{\CC_n(v)[1,\underline w^{-\sigma} n)\nsim \CV_n[\underline w^{-\sigma} n, \infty)} 
=\sum_{\ell\ge 1}\sum_{\CC\in \CG_n[1,\underline w^{-\sigma} n):|\CC|=\ell}\hspace{-15pt}\ell\cdot\ind{\CC\nsim \CV_n[\underline w^{-\sigma} n, \infty)} .
\end{equation}
Each hub connects to any vertex with probability $q$ by~\eqref{eq:conn-prob}. Hence, there is no edge between a  component of size $\ell$ and any of the hubs with probability $(1-q)^{h\ell}$. 
Moreover, as stated in~\eqref{eq:lln-intro}, the number of vertices in components of size $\ell$ in $\CG_n[1,\underline w^{-\sigma} n)$ is roughly $n\cdot\Prob(|T_q|=\ell)$. Therefore, the total number of size-$\ell$ components is about $(n/\ell)\cdot\Prob(|T_q|=\ell)$. Hence, 
\begin{equation}\label{eq:concentrated-clusters}
n-|\CC_n^\sss{(1)}| \approx \sum_{\ell\ge 1} (n/\ell)\Prob(|T_q|=\ell)\cdot\ell\cdot (1-q)^{h\ell}=n\E\big[(1-q)^{h|T_q|}\big].
\end{equation}
For the event $\{|\CC_n^\sss{(1)}|>\rho n\}$ to occur,
the complement of the giant should  contain at most $n(1-\rho)$ vertices. Thus, if the number of hubs $h$ is the smallest integer satisfying
\begin{equation}\label{eq:hub-informal}
\E\big[(1-q)^{h|T_q|}\big] \le 1-\rho,
\end{equation}
then we expect that the largest component has size at least $\rho n$. 
In our proofs below, we formalize this reasoning and show that any other ``strategy'' is less effective in increasing the size of the giant. Our proof makes the above reasoning more precise: by controlling vertex weights in size-$\ell$ components, we estimate the impact of adding one hub with any weight $yn$ as a function of $y$. 
As a result, we can analyze the joint distribution of the weights of the hubs that lead to a giant of size at least $\rho n$, and show that the number of components of constant size decreases as described above. 
In the next section, we formalize  the large-deviation principle for the giant; the organization of the remainder of the paper is given there.

\section{Main results}\label{sec:ldp-giant}
Motivated by~\eqref{eq:hub-informal}, we define for $z\in[0,1]$ the probability generating function of the associated branching process restricted to be finite:
\begin{equation}
    H_{T_q}(z):=\E\big[z^{|T_q|}\ind{|T_q|<\infty}\big].
\end{equation}
This function is increasing, continuous, and has range $[0, 1-\theta_q]$ for $z\in[0,1]$. Hence, its inverse $y\mapsto H_{T_q}^{(-1)}(y)$ is well-defined for $y\in[0, 1-\theta_q]$.
We define for $\rho\in[0, 1)$
\begin{equation}
    \mathrm{hubs}(\rho, q):=\begin{dcases}
        \frac{\log H_{T_q}^{(-1)}(1-\rho)}{\log(1-q)},&\text{if }\rho>\theta_q, \text{and } q<1,\\
        1,&\text{if }\rho>\theta_q, \text{and } q=1,\\
        0,&\text{if }\rho\le\theta_q.
    \end{dcases}\label{eq:hubs}
\end{equation}
While the definition as the inverse of a generating function may seem abstract, its asymptotics are explicitly computable as shown in Lemma~\ref{lemma:properties-hubs} below.
By the back-of-the envelope calculation, $\lceil \mathrm{hubs}(\rho, q)\rceil$ hubs are able to increase the proportion of vertices in the giant component from $\theta_q$ to $\rho\in(\theta_q,1)$. 
The following set describes the  weights (rescaled by a factor $1/n$) of the hubs that are jointly able to increase the proportion of the vertices in the giant from $\theta_q$ to $\rho\in(\theta_q, 1)$ with sufficiently large probability. 
 Let, for $h\in\N$,
\begin{equation}
\CY_{\rho, q}(h) := \bigg\{(y_1,\ldots, y_{h})\in (0,\infty)^{h}: \E\bigg[\prod_{x\in T_q,\\ i\in[h]}\big(1-q\cdot(y_iW_x^\sigma\wedge 1)\big)\bigg]\le 1-\rho\bigg\}.
    \label{eq:y-rho-set}
\end{equation}
 The expectation in this definition is similar to the expectation in~\eqref{eq:hub-informal}, and represents the probability that the component of a uniform vertex (represented by $T_q$, containing vertices with weight $(W_x)_{x\in T_q}$) does not connect to any of the hubs with weights $\{y_1n,\ldots, y_h n\}$. 
In \cref{lemma:leading-term} below, we show that $\CY_{\rho, q}(\lceil\mathrm{hubs}(\rho, q)\rceil)$ is non-empty, and that the set does not contain any points in a small neighborhood around the origin if $\mathrm{hubs}(\rho, q)\notin\N$ or $q=1$. 

Let $\rho\in(\theta_q, 1)$, $q\in(0,1]$,   $h=\lceil\mathrm{hubs}(\rho, q)\rceil$, and $\alpha$ as in~\eqref{eq:weight-dist}.
We define the constant 
\begin{equation}
    C_{\rho, q}:=\frac{\alpha^h}{h!}\int_{y_1=0}^\infty\cdots\int_{y_h=0}^\infty\mathds{1}\{(y_1,\ldots, y_h)\in\CY_{\rho, q}(h)\}\cdot (y_1\cdot\ldots\cdot y_h)^{-(\alpha+1)}\rd y_1\cdot\ldots\cdot\rd y_h.\label{eq:constant-value}
\end{equation}
Our main result is the following theorem.
\begin{theorem}[Upper tail for the giant]\label{thm:upper-tail}
Consider an inhomogeneous scale-free random graph as in \cref{def:irg}. Fix a constant $\rho\in(\theta_q, 1)$. If $\mathrm{hubs}(\rho,q)\notin\N$ or $q=1$, then, as $n\to\infty$, 
\begin{equation}
\Prob\big(|\CC_n^\sss{(1)}|> \rho n\big)
\sim C_{\rho, q}\big(n\Prob(W>
 n)\big)^{\lceil \mathrm{hubs}(\rho, q)\rceil }.\label{eq:upper-tail-1}
\end{equation}
If $\mathrm{hubs}(\rho,q)\in \N$ and $q<1$, there exists a constant $c>0$ such that, as $n\to\infty$, 
\begin{equation}
c \big(n\Prob(W_1>
 n)\big)^{\mathrm{hubs}(\rho, q)+1}\le 
     \Prob\big(|\CC_n^\sss{(1)}|> \rho n\big)
\lesssim C_{\rho, q} \big(n\Prob(W_1>
 n)\big)^{\mathrm{hubs}(\rho, q)}.\label{eq:upper-tail-2}
\end{equation}
\end{theorem}
Thus, the lower bound and upper bound coincide up to smaller order terms when $\mathrm{hubs}(\rho, q)\notin\N$ or $q=1$. At the discontinuity points of $\lceil \mathrm{hubs}(\rho, q)\rceil$, i.e., when $\mathrm{hubs}(\rho, q)\in\N$ and $q<1$, the decay rate of the upper and lower bound differ by a regularly-varying factor.
The following theorem illustrates that the lower tail of large deviation decays exponentially fast. 
\begin{theorem}[Lower tail for the giant]\label{thm:lower-tail}
Consider an inhomogeneous scale-free random graph as in \cref{def:irg}. For all $\rho<\theta_q$, there exists a constant $c=c(\theta_q)>0$ such that for all $n\ge1$
\begin{equation}\label{eq:thm-lower-tail}
\Prob\big(|\CC_{n}^\sss{(1)}|< \rho n\big)\le \exp\big(-cn).
\end{equation}
\end{theorem}
We remark that more precise results than~\eqref{eq:thm-lower-tail} could be derived by adjusting the quenched LDP (i.e., with fixed weight sequence) from~\cite{andreis2023irg} to the annealed setting with iid weights. We focus on the upper tail, and leave this open. 
Theorems \ref{thm:upper-tail} and \cref{thm:lower-tail} imply a large-deviation principle with logarithmic speed, and rate function
\begin{equation}\label{eq:rate-function}
I_q(\rho):= \begin{dcases}
    (\alpha-1)\lceil\mathrm{hubs}(\rho, q)\rceil,&\text{if }\rho\in\big[\theta_q, 1\big),\\
    +\infty,&\text{otherwise.}
\end{dcases}
\end{equation}
Because $\lim_{\rho\uparrow 1}\mathrm{hubs}(\rho, q)=+\infty$, the function $I_q(\rho)$ is lower semi-continuous, 
and finite for $\rho\in[\theta_q,1)$.
Let $B^\circ$ and $\bar B$ denote the interior and closure of a Borel set $B\subseteq \R$.
\begin{corollary}[Large-deviation principle]\label{cor:ldp-giant}
    Consider an inhomogeneous scale-free random graph as in \cref{def:irg}. Then for any Borel set $B\subseteq \R$ 
        \[
    \begin{aligned}
    - \inf_{\rho\in B^\circ}I_q(\rho)\le \liminf_{n\to\infty}\frac{1}{\log n}\log \Prob\bigg(\frac{|\CC_{n}^\sss{(1)}}n\in B\bigg)
    \le \limsup_{n\to\infty}\frac{1}{\log n}\log \Prob\bigg(\frac{|\CC_n^\sss{(1)}|}n \in B\bigg)\le -\inf_{\rho \in \bar B} I_q(\rho).
    \end{aligned}
    \]
\end{corollary}
\begin{remark}[The rate function is global] By definition of $\mathrm{BP}(\kappa_\sigma, F_W,q)$ in \cref{def:branching}, there are multiple pairs $(F_W, q)$ that lead to the same distribution of $T_q$. 
A straightforward coupling yields an alternative construction for $\mathrm{BP}(\kappa_\sigma, F_W, q)$: first sample $\mathrm{BP}(\kappa_\sigma, F_W, 1)$, then remove each edge independently with probability $q$. The tree in this forest that contains $\varnothing$ is equal in distribution to  $\mathrm{BP}(\kappa_\sigma, F_W, q)$. Moreover, the parameter $q$ can be encapsulated by a reparametrization of $F_W$ to some $F_{W,q}$ without affecting the distribution of the associated branching process, i.e., $\mathrm{BP}(\kappa_\sigma, F_W, q)=\mathrm{BP}(\kappa_\sigma, F_{W, q}, 1)$ in distribution. Therefore, the two IRGs with the parametrizations $(\kappa_\sigma, F_W, q)$ and $(\kappa_\sigma, F_{W,q}, 1)$ have the same local limit.

However, the IRGs with parameters $(\kappa_\sigma, F_W, q)$ and $(\kappa_\sigma, F_{W,q}, 1)$ are different in distribution: for example, when $\sigma>0$ and $q<1$, vertices with weight $\omega(n^{1/(\sigma+1)})$ are not connected by an edge with probability $1-q>0$, while they are connected almost surely when $q=1$. Therefore, the function $\mathrm{hubs}(\rho, q)$ and the rate function $I$ do not agree for such parameter pairs, and we conclude that the rate function is not determined uniquely by the local limit, and is a global quantity instead. This contrasts the result that the \emph{typical} size of the giant is uniquely described by its local limit~\cite{bollobas2007irg,Hofbook2}.
\end{remark}

We end this section with two corollaries that describe the graph structure  \emph{conditionally} on the rare event $\{|\CC_n^\sss{(1)}|>\rho n\}$. They illustrate that the intuition given above is correct, and that the rare event is indeed caused by exactly $\lceil\mathrm{hubs}(\rho,q)\rceil$ hubs of which the rescaled weights are in the set $\CY_{\rho, q}$, which impacts the empirical component-size distribution depending on the exact weights. All other vertices have sublinear weight.
\begin{corollary}[Conditional weight distribution of the hubs]\label{cor:cond-dist-weight}
    Consider an inhomogeneous scale-free random graph as in \cref{def:irg}. Fix $\rho\in\big(\theta_q, 1\big)$ such that  $\mathrm{hubs}(\rho, q)\notin\N$ or $q=1$. Let $(Y_i)_{i\ge 1}$ be independent copies of $Y$ following distribution $\Prob(Y\ge y)=y^{-\alpha}$ for $y\ge 1$. For any sufficiently small constants   $\eps,\phi>0$, as $n\to\infty$,
    \[
    \big\{W_u/n: u\in \CV_n[n^{1-\eps},\infty)\big\} \,\, \Big|\,\,|\CC_n^\sss{(1)}|>\rho n \quad \overset d\longrightarrow \quad \big\{\phi Y_i\big\}_{i\le \lceil\mathrm{hubs}(\rho, q)\rceil}\, \Big|\, \big\{\phi Y_i\big\}_{i\le \lceil\mathrm{hubs}(\rho, q)\rceil}\in \CY_{\rho, q}.
    \]
\end{corollary}
The above corollary does not describe the empirical weight distribution of the vertices with weight at most $n^{1-\eps}$. Using methods similar to~\cite{kerriouMaitraCondensation} it can be proven that the empirical weight distribution of these vertices, conditionally on $\{|\CC_n^\sss{(1)}|>\rho n\}$, converges weakly to $F_W$. We leave the formal proof open, as it would lead to many more technicalities below. 

We proceed to the conditional component-size distribution, of which we show that it converges to a random sequence. 
Define for $h\in\N$ and $(y_1,\ldots, y_h)\in(0,\infty)^h$
\begin{equation}\label{eq:g-comp-size}
g_{\ell}\big((y_i)_{i\le h}\big):=\frac{1}{\ell}\E\Big[\ind{|T_q|=\ell}\prod_{x\in T_q, i\le h}\big(1-q(W_x^\sigma y_i\wedge 1)\big)\Big].
\end{equation}
We write $N_{n,\ell}$ for the number of components of size $\ell$ in $\CG_n$.  Let $R^\infty$ denote the space of all sequences ${\bf{x}}=(x_1,x_2,\ldots)$ of real numbers, metrized by $d_\infty({\bf{x}}, {\bf{y}})=\sum_i(|x_i-y_i|\wedge1 )2^{-i}$. 
\begin{corollary}[Conditional component-size distribution]\label{cor:cond-dist-comp}
    Consider an inhomogeneous scale-free random graph as in \cref{def:irg}. Fix $\rho\in\big(\theta_q, 1\big)$ such that  $\mathrm{hubs}(\rho, q)\notin\N$  or $q=1$. Let $(Y_i)_{i\ge 1}$ be independent copies of $Y$ following distribution $\Prob(Y\ge y)=y^{-\alpha}$ for $y\ge 1$. For any sufficiently small constant   $\phi>0$, as $n\to\infty$,
    \[
    \big(N_{n,\ell}/n, \ell\ge 1\big) \,\, \Big|\,\,|\CC_n^\sss{(1)}|>\rho n \quad \overset {\mathcal{D}}\longrightarrow \quad\big(g_{\ell}\big((\phi Y_i)_{i\le \lceil\mathrm{hubs}(\rho, q)\rceil}\big), \ell\ge 1\big)\, \big|\, \big(\phi Y_i\big)_{i\le \lceil\mathrm{hubs}(\rho, q)\rceil}\in \CY_{\rho, q}.
    \]       
    Here, ${\mathcal{D}}\rightarrow$ denotes weak convergence in $R^\infty$. 
    Moreover, conditionally on $|\CC_n^\sss{(1)}|>\rho n$, $|\CC_n^\sss{(1)}|/n$ converges in distribution to a random variable $Q$ supported on $\big[\rho, \inf\{\rho'>\rho: \mathrm{hubs}(\rho', q)\in \N\}\big]$  with distribution
    \begin{equation}\label{eq:cond-size-giant}
    \Prob\big(Q>s\big)=\frac{C_{s,q}}{C_{\rho, q}}, \qquad \text{}\quad 
    s\in\big[\rho, \inf\{\rho'>\rho: \mathrm{hubs}(\rho', q)\in \N\}\big).
    \end{equation}
\end{corollary}

We finish this section with a lemma that illustrates some properties of $\mathrm{hubs}$ defined in~\eqref{eq:hubs}.
 Only the first two items of the lemma will be required in the proofs of our main results. The proof is postponed to the appendix on page~\pageref{proof:properties-hubs}.
\begin{lemma}[Properties of $\mathrm{hubs}(\rho, q)$]\label{lemma:properties-hubs}
    Consider an inhomogeneous scale-free  random graph as in \cref{def:irg}. Let $W_\varnothing, W$ be two independent random variables following distribution $F_W$, and assume $\rho>\theta_q$. Then, 
    \begin{enumerate}[leftmargin=2.2em]    
        \item[(i)] the following identity holds:        
\begin{equation}
    \mathrm{hubs}(\rho, q)=\inf\Big\{h'>0: \E\big[(1-q)^{|T_q|h'}\big]\le 1-\rho\Big\}.\label{eq:hubs-expectation}
\end{equation}
        \item[(ii)] $\mathrm{hubs}(\rho, q)$ is continuous, positive, non-decreasing in $\rho$, and decreasing in $q$. 
        \item[(iii)] if either $\rho$ is fixed and $q\downarrow0$, or $q<1$ is fixed and $\rho\uparrow 1$, then $\mathrm{hubs}(\rho, q)$ tends to infinity, and 
        \[
        \mathrm{hubs}(\rho, q)\sim\frac{\log\big(1/(1-\rho)\big)}{\log\big(1/(1-q)\big)}.
        \]
        \item[(iv)] if $q<1$ is fixed, as $\rho\uparrow 1$,
        \[
        \mathrm{hubs}(\rho, q)=\frac{\log\big(1/(1-\rho)\big) - \log\big(\E\big[\exp\big(-q\E[\kappa_\sigma(W_\varnothing,W)\mid W_\varnothing]\big)\big]\big)}{\log\big(1/(1-q)\big)} + o(1),
        \]
        which simplifies for rank-one inhomogeneous scale-free random graphs, i.e., $\sigma=1$, to 
        \[
        \mathrm{hubs}(\rho, q)=\frac{\log\big(1/(1-\rho)\big) - \log\big(\E\big[\exp\big(-q\E[W]W_\varnothing\big)\big]\big)}{\log\big(1/(1-q)\big)} + o(1), \qquad\text{as }\rho\uparrow 1.
        \]
    \end{enumerate}
\end{lemma}

\subsection{Discussion and related work}
This paper provides the first large-deviation principle (LDP) for component sizes in random graph models with iid regularly varying weights, and enriches the emerging literature on LDPs for other graph properties. 
The rate function in Corollary~\ref{cor:ldp-giant} has the most interesting behavior when the edge-percolation parameter $q$ is strictly smaller than $1$; when $q=1$, a single hub can increase the size of the giant to $n$.

An applied example where varying $q$ is of interest can be found in epidemiology: the largest component in an edge-percolated graph corresponds to the final size of a Reed--Frost epidemic~\cite{abbey1952examination, barbour2004approximating, lefevre1990stochastic} on the largest component of the graph before percolation. The parameter $q$ represents the probability that a disease is transmitted along an edge. Thus, Theorem~\ref{thm:upper-tail} translates to the asymptotic probability that the epidemic has a significantly larger size than expected; Theorem~\ref{thm:lower-tail} examplifies that a significantly smaller final size than expected is less likely to occur in IRGs. 

The rate function in Corollary~\ref{cor:ldp-giant} behaves drastically different compared to the rate function of LDPs for the number of edges $E_n$~\cite{kerriouMaitraCondensation, kerriouFewest, stegehuis2023triangles, stegehuis2023edges}, and  the number of triangles $\triangle_n$~\cite{stegehuis2023triangles}. Similar to Theorem~\ref{thm:upper-tail}, the LDP for $E_n$ requires constantly many hubs of linear weight to find more than $a n$ edges than expected, and yields a discontinuous rate function in $a$. The number of hubs scales linearly in $a$ as $a\to\infty$ in~\cite{kerriouMaitraCondensation, stegehuis2023edges}, while in our setting $\mathrm{hubs}(\rho, q)=\Theta(\log(1/(1-\rho)))$ as $\rho\uparrow 1$. For $\triangle_n$ the situation is different~\cite{stegehuis2023triangles}: when the parameters are such that $\E[\triangle_n]=o(n)$, a single vertex of large (but sublinear) weight is required to find $a\E[\triangle_n]$ additional triangles, which occurs with a probability decaying polynomially in $n$; when $\E[\triangle_n]=\omega(n)$, polynomially many additional vertices of weight $\Omega(n^{1/(\sigma+1)})$ are required, which occurs with a probability that decays stretched exponentially. 

The rate function $I_q(\rho)$ in~\eqref{eq:rate-function} also appears in the related work on the size of the giant component in scale-free random graphs that are embedded in Euclidean space~\cite{clusterII}, where it is shown that $\Prob(|\CC_n^\sss{(1)}|>\rho n)=\Theta(n^{-I_q(\rho)})$ when $\mathrm{hubs}(\rho, q)\notin\N$. Parts of our proofs are related to the techniques from~\cite{clusterII}. However, the way that we obtain concentration bounds for the number of components of size $\ell$ in the graph without hubs is different:  there is no geometry that we can leverage to obtain concentration via (almost independent) disjoint graphs in subboxes as in~\cite{clusterII}. Instead, we use truncation/discretization arguments in combination with the LDP from~\cite{andreis2023irg}. Compared to~\cite{clusterII},
Theorem~\ref{thm:upper-tail} additionally provides the constant prefactor $C_{\rho, q}$, which requires substantial additional analysis. Deriving this constant allows us to prove the LDP in Corollary~\ref{thm:upper-tail}, and to give a detailed description of the graph conditionally on $\{|\CC_n^\sss{(1)}|>\rho n\}$. 

A commonality among Theorem~\ref{thm:upper-tail} and the above-mentioned works~\cite{kerriouMaitraCondensation,  stegehuis2023triangles, stegehuis2023edges}, is that the lower bound and upper bound are non-matching on  discontinuity points of $I_q(\rho)$. As argued in~\cite[Section 2.4]{kerriouMaitraCondensation} and~\cite[Section 3]{stegehuis2023edges}, a delicate analysis is required to find the correct scaling for such $\rho$, and the decay may heavily depend on the precise form of the connection probability in~\eqref{eq:conn-prob}.

The setup of Definition~\ref{def:irg} contains (a non-spatial version of) the age-dependent random connection model (ADCM) by setting $\sigma=\alpha-1$~\cite{GraHeyMonMor19}, which in several contexts mimics preferential attachment models (PAMs)~\cite{dereich2012typical,gracar2019age,jacobLinkerMortersFastEvolving}. However, we expect that our results here do not extend to PAMs, and that $\Prob(|\CC_n^\sss{(1)}|>\rho n)$ decays exponentially. Informally speaking, the main difference is that the connection probability in ADCM depends on the \emph{exact age} of vertices, while in PAMs only the \emph{order of the ages} affects connection probabilities. On the contrary, we believe that the edge-percolated configuration model with iid regularly varying degrees behaves similar to  IRGs, and that polynomial decay of $\Prob(|\CC_n^\sss{(1)}|>\rho n)$ can be proven using adaptations of our techniques. This would contrast the exponential decay from~\cite{bhamidiConfLargeDev} for the configuration model with given degrees in which no linear-sized hubs are allowed. 

Theorems~\ref{thm:upper-tail} and~\ref{thm:lower-tail} provide the asymptotic probability that the size of the giant deviates by a constant factor from its expectation if the parameters of the IRG are such that $\theta_q>0$, i.e., when the model is supercritical. However, the theorems also apply when $\theta_q=0$, in which case the model is subcritical or critical and the largest component grows sublinear in $n$. For such parameters, Theorem~\ref{thm:upper-tail} describes deviations of $|\CC_n^\sss{(1)}|$ at a much larger scale than $\E[|\CC_n^\sss{(1)}|]$, and Theorem~\ref{thm:lower-tail} is trivial. This poses a natural follow-up question: when $\theta_q=0$, what is the decay of $\Prob\big(|\CC_n^\sss{(1)}|<\beta_n\big)$ and $\Prob\big(|\CC_n^\sss{(1)}|>\gamma_n\big)$ at intermediate scales, i.e.,  when $1\ll\beta_n\ll\E[|\CC_n^\sss{(1)}|]\ll\gamma_n\ll n$? For critical rank-one IRGs with infinite third moment ($\sigma=1, \alpha\in(2,3)$) the component-size tails were studied in~\cite{kliemClusterTails}, but we are not aware of further results in this direction.

We describe the organization of the remainder of the paper via a proof sketch of Theorem~\ref{thm:upper-tail}.

\subsection{Outline of the proof}\label{sec:methodology}
In Section~\ref{sec:size-ell} we analyze the graph without hubs (vertices with weight at least $\phi n$ for some small constant $\phi>0$) in two stages. 

We first analyze the induced subgraph of $\CG_n$ on all vertices with weight at most a large constant $ R$, which we denote by $\CG_n[\underline w, R)$. Following the calculation on page~\pageref{back-of-the-envelope}, we have to  show for each constant $\ell\in\N$ the proportion of vertices in size-$\ell$ components in this graph concentrates around the probability that the total progeny of the associated branching process has size $\ell$. To estimate the precise effect of adding one hub with any weight at least $\phi n$, we additionally control the vertex weights in size-$\ell$ components. 

To do so, we discretize the interval $[\underline w, R)$ into small intervals and categorize each component of size $\ell$ in $\CG_n[\underline w, R)$ by the intervals that contain the weights of the $\ell$ vertices. For each possible category ---that we call component type--- we establish concentration of the proportion of vertices in such a component around the probability that the associated branching process has size $\ell$, and that the weights of the vertices in the total progeny fall exactly in the same intervals. We achieve this via a coupling with inhomogeneous random graphs with given weights, and then rely on a large-deviation principle from~\cite{andreis2023irg}. This discretization is at the core of the proof, and allows us to derive the explicit constant $C_{\rho,q}$ in \cref{thm:upper-tail}.
Afterwards, we add the vertices with weight in $[R, \phi n)$ to the graph and show that the edges incident to such vertices barely affects the number of components of any type. The size of the largest component in this graph is still about $\theta_q n$.

In Section~\ref{sec:reduction} we analyze the impact of vertices of weight at least $\phi n$ to the graph. When $\mathrm{hubs}(\rho, q)\notin\N$, we show that if there are exactly $\lceil\mathrm{hubs}(\rho, q)\rceil$ vertices with weight in the set $\CY_{\rho+\delta, q}$ for any $\delta>0$, then the size of the giant increases from $\theta_q n$ to at least $\rho n$ whp. Since $\CY_{\rho+\delta, q}\subset\CY_\rho$ by~\eqref{eq:y-rho-set}, we require slightly larger weights than explained before. On the contrary, if the number of hubs is not equal to $\lceil\mathrm{hubs}(\rho, q)\rceil$, or if their weights are not contained in $\CY_{\rho-\delta, q}\supset\CY_{\rho, q}$, then the size of the giant is very unlikely to increase above $\rho n$:  the probability of this event is of the same order as the event that there are strictly more than $\lceil\mathrm{hubs}(\rho, q)\rceil$ hubs. The presence of the small constant $\delta>0$ mitigates the effect of the truncation and discretization from Section~\ref{sec:size-ell}. At the end of Section~\ref{sec:reduction} we show that we can take the limit $\delta\to0$ when $\mathrm{hubs}(\rho, q)\notin\N$, showing that $\Prob\big(|\CC_n^\sss{(1)}|>\rho n\big)\sim \Prob\big(|\CV_n[\phi n,\infty)|=\lceil\mathrm{hubs}(\rho, q)\rceil, \CV_n[\phi n,\infty)\in \CY_{\rho, q}\big)\sim C_{\rho,q}(n\Prob(W>n))^{\lceil\mathrm{hubs}(\rho,q)\rceil}$; proving the first bound in \cref{thm:upper-tail}. 

Section~\ref{sec:main-proofs} formally verifies Theorems~\ref{thm:upper-tail} and~\ref{thm:lower-tail} and proves the LDP for the giant in Corollary~\ref{cor:ldp-giant}.
Lastly, in Section~\ref{sec:remainder} we analyze the graph conditionally on the rare event $\{|\CC_n^\sss{(1)}|>\rho n\}$ to prove Corollaries~\ref{cor:cond-dist-weight} and~\ref{cor:cond-dist-comp}. The appendix contains the proof of Lemma~\ref{lemma:properties-hubs}, and proofs of some technical lemmas from the following sections, and standard concentration bounds that we frequently use.

\section{The graph without hubs}\label{sec:size-ell}
In this section we analyze the graph in which all vertices with weight at least $\phi n$ are removed. Our main goal is to obtain concentration bounds for components of size $\ell$, together with their weight configuration. We introduce notation to categorize components. 
\begin{definition}[Type-$(\bf{\widetilde w},\varepsilon)$ component] \label{def:type-w-comp}
Fix $\ell\in\N$ and a small constant $\eps>0$, and consider a partitioning of $[\underline w,\infty)$ into intervals of length $\varepsilon$. Let $\mathrm{VT}(\eps):=\{x_i\}_{i\ge 1}:=\{\underline w+i\eps\}_{i\ge 1}$ denote the lower boundaries of these intervals.
Let ${\bf{\widetilde w}}^\sss{(\ell)}=(\widetilde w_1,\ldots, \widetilde w_\ell)\in \mathrm{VT}(\eps)^{\ell}$ be a vector. We say that a connected component $\CC$ of a vertex-weighted graph $G$ has type $({ \bf{\widetilde w}}^\sss{(\ell)},\varepsilon)$, if its number of vertices is $\ell$, and if there exists an ordering of the vertices $(v_1,\ldots, v_\ell)$ such that $\tilde w_i\le w_{v_i}\le \tilde w_i+\eps$ for all $i\in[\ell]$. Let $\overline w>\underline w$ be a weight threshold. We write $N_n({\bf{\widetilde w}}^\sss{(\ell)},\varepsilon, \overline w)$ for the number of components of type $({\bf{\widetilde w}}^\sss{(\ell)}, \varepsilon)$ in $\CG_n[\underline w,\overline w)$. For the total progeny $T_q$ of the associated branching process of $\CG_n$, we define 
\begin{equation}\label{eq:branching-type-prob}
\theta({\bf{\widetilde w}}^\sss{(\ell)},\varepsilon):=\Prob\big(T_q\mbox{ has type }({\bf{\widetilde w}}^\sss{(\ell)}, \varepsilon)\big).
\end{equation}
 We write $\mathrm{CT}_\ell(\eps)\subseteq \mathrm{VT}(\eps)^{\ell}$ for the set of component types ${\bf{\widetilde w}}^\sss{(\ell)}$ with weights in $\mathrm{VT}_\ell(\eps)$ of size $\ell$. Let $R$ be a large constant so that the interval $[\underline w, R)$ can be partitioned into $(R-\underline w)/\varepsilon$ intervals of length $\varepsilon$. We write $\mathrm{CT}_\ell(\eps, R)\subseteq \mathrm{VT}(\eps, R)^{\ell}$ for the finite set of component types ${\bf{\widetilde w}}^\sss{(\ell)}$ with weights in $\mathrm{VT}(\eps, R):=\{\underline w+i\eps\}_{i\in[(R-\underline w)/\eps]}$.
\end{definition}
The following lemma establishes concentration bounds for the proportion of vertices in  type-$({\bf{w}}^\sss{(\ell)},\varepsilon)$ components. 
\begin{lemma}[Concentration of size-$\ell$ components]\label{lemma:ell-trunc-n}
    Consider an inhomogeneous scale-free random graph as in \cref{def:irg}. 
    For any constants $\psi,  R,  \ell_\ast>0$, there exists a constant $c>0$ such that for all $\eps>0$ such that $(R-\underline w)/\eps\in\N$, and all $n\ge 1$,
    \begin{align}
        \Prob\Bigg(\sum_{\ell\le \ell_\ast}\sum_{  ({\bf{x}}^\sss{(\ell)},\eps)\in\mathrm{CT}_\ell(\eps, R)}  \Big|\frac{\ell N_n({\bf{w}}^\sss{(\ell)}, \varepsilon, R )}n- \theta({\bf{w}}^\sss{(\ell)}, \varepsilon)\Big|>\psi\Bigg)&\le \exp\big(-cn\big).\label{eq:lemma-ell-lower-R}
    \intertext{Moreover, for any constants $\psi,  R, C, \ell_\ast>0$, there exist constants $\phi, c>0$ such that for all $\eps>0$ such that $(R-\underline w)/\eps\in\N$, as $n\to\infty$,}
    \Prob\Bigg(\sum_{\ell\le \ell_\ast}\sum_{  ({\bf{x}}^\sss{(\ell)},\eps)\in\mathrm{CT}_\ell(\eps, R)}   \Big|\frac{\ell N_n({\bf{w}}^\sss{(\ell)}, \varepsilon, \phi n)}n - \theta({\bf{w}}^\sss{(\ell)}, \varepsilon)\Big|>\psi\Bigg)&=o\big(n^{-C}\big).\label{eq:lemma-ell-lower}
    \end{align}
\end{lemma}
We proceed with a lemma for the lower tail of the largest component in truncated graphs.
\begin{lemma}[Lower tail of the truncated giant]\label{lemma:giant-trunc-lower}
    Consider an inhomogeneous scale-free random graph as in \cref{def:irg}. 
    For any constant $\rho\in(0,\theta_q)$, there exist constants $R, c>0$ such that for all $n\ge 1$ and all $\overline w\ge R$,
    \begin{equation}
    \Prob\big(|\CC_n^\sss{(1)}[\underline w, \overline{w})|<\rho n\big)\le\exp(-cn).
    \end{equation}
\end{lemma}

In the next section we analyze the graph on constant weight vertices, and the obtained results already imply \cref{lemma:giant-trunc-lower} and~\eqref{eq:lemma-ell-lower-R} in \cref{lemma:ell-trunc-n}. Then in Section~\ref{sec:intermediate} we analyze the effect of adding vertices of intermediate-weight, i.e., in the interval $[R, \phi n)$, allowing to prove~\eqref{eq:lemma-ell-lower} in \cref{lemma:ell-trunc-n}.

\subsection{The graph on constant-weight vertices}
In this section we will construct an ``approximating'' IRG ($\delta$-IRG) with discretized weights of at most a constant, and use  that  the number of type--$({\bf{x}}^\sss{(\ell)}, \eps)$ components concentrates around its expectation in the $\delta$-IRG for all types simultaneously. For this, we rely on the large-deviation principle for component sizes in inhomogeneous random graphs (IRGs) with bounded kernel and deterministic weights from~\cite{andreis2023irg}. 

Then, on the event that the weights in $\CG_n[\underline w,R)$ satisfy a good event, we  construct a coupling between the $\delta$-IRG with the original graph $\CG_n[\underline w, R)$. We show that the number of components of any type changes compared to the $\delta$-IRG only by a small number with sufficiently high probability under this coupling, so that \cref{lemma:giant-trunc-lower} and~\eqref{eq:lemma-ell-lower-R} in \cref{lemma:ell-trunc-n} follow.

\begin{definition}[Approximating inhomogeneous random graph and branching process]\label{def:approx}Let $\delta,\eps>0$ be  small constants such that $\eps/\delta\in\N$, and let $R$ be a large constant. Let the weight-distribution $F_W$, kernel $\kappa_\sigma$, and percolation parameter $q$ be as in \cref{def:irg}, and recall that $\underline w=\inf\{w: F_W(w)>0\}>0$. Let $z_i=\underline w+i\delta$ for $i\in\N_0$. Set
\[
f_W^\sss{(\delta)}(z_i):= \Prob\big(z_i\le W< w_{i+1}\big), \qquad  \underline n_i:=\big\lceil (1-\delta) n f_W(z_i)\big\rceil.
\]
For $i\in[(R-\underline w)/\delta)]$, we fix $\underline n_i$ to be the number of vertices of weight exactly $z_i$ in $\CV_n^\sss{(\delta, R)}$, and write 
\begin{equation}\label{eq:approx-n-under}
\underline n=\underline n(\delta, R):=\sum_{i\in [(R-\underline w)/\delta)]}\underline n_i
\end{equation}
for the total number of vertices in $\CG_n^\sss{(\delta, R)}$. 
For a weight $w\in[\underline w, R)$, let $\underline w^\sss{(\delta)}:=\sup_i\{z_i: z_i\le w\}$, $\bar w^\sss{(\delta)}:=\inf_i\{z_i: z_i> w\}$.
Define for two weights $w_1, w_2\in[\underline w, \infty)$ the approximating kernel as
\begin{equation}\label{eq:kernel-adj}
\kappa_\sigma^\sss{(\delta, R)}(w_1, w_2):= \ind{w_1\le R, w_2\le R}\inf\big\{\kappa_\sigma(\tilde w_1, \tilde w_2): \tilde w_1\in[\underline w_1^\sss{(\delta)}, \overline w_1^\sss{(\delta)}],  \tilde w_2\in[\underline w_2^\sss{(\delta)}, \overline w_2^\sss{(\delta)}]\big\}.
\end{equation}
 Two vertices $u, v$ of weight $w_u, w_v$ are connected in $\CG_n^\sss{(\delta, R)}$ with probability 
 \[
 p^\sss{(\delta, R)}_{uv}:= q\bigg(\frac{\kappa_\sigma^\sss{(\delta, R)}(w_u, w_v)}{n}\wedge 1\bigg)
 \] independently of other vertex pairs. 
 The approximating associated branching process is denoted by $\mathrm{BP}^\sss{(\delta, R)}$. 
 In this branching process, the root $\varnothing$ has type $W_\varnothing$ where \[\Prob\big(W_\varnothing=z_i\big)=f_W^\sss{(\delta)}(z_i)\,\big/\,F_W(R).\]
 In each generation, each particle $v$ of weight $w_v$ gives, for each $i\in[(R-\underline w)/\delta]$, independently birth to $\mathrm{Poi}\big(q\cdot\kappa^\sss{(\delta, R)}_\sigma(w_v, w)\cdot f_W^\sss{(\delta)}(z_i)\big)$ many particles of weight $z_i$. We denote the set of types in the total progeny by $T_q^\sss{(\delta, R)}$
\end{definition}

 We write $\CC_n^\sss{(1), (\delta, R)}$ for the largest connected component in $\CG_{n}^\sss{(\delta, R)}$, 
$N_n({\bf{w}}^\sss{(\ell)},\eps)$ for the number of type-$({\bf{w}}^\sss{(\ell)},\eps)$ components in the $\delta$-IRG,
$\theta_q^\sss{(\delta, R)}$ for the survival probability of the associated branching process ($\delta$-BP) of the $\delta$-IRG, and $\theta_q^\sss{(\delta, R)}({\bf{x}}^\sss{(\ell)}, \eps)$ for the probability that the $\delta$-BP has size $\ell$ and type $({\bf{x}}^\sss{(\ell)}, \eps)$.

We state a lemma that follows from a large-deviation principle by Andreis \emph{et al.}~\cite{andreis2023irg}. 
\begin{lemma}[Size-$\ell$ components in the $\delta$-IRG]\label{lemma:andreis}
    Consider the approximated inhomogeneous random graph and its associated branching process for some $R>0$, $\delta>0$. Fix $\ell_\ast\in\N$ and $\eps>0$. 
    For each $\psi>0$, there exist constants $c_{\psi, 1},\delta_0, R_0>0$ such that for all $\delta\in(0,\delta_0)$ such that $\eps/\delta\in\N$, and $n\ge 1$
    \begin{equation}\label{eq:andreis-ell}
    \Prob\Bigg(\sum_{\ell\le \ell_\ast}\sum_{  ({\bf{x}}^\sss{(\ell)},\eps)\in\mathrm{CT}_\ell(\eps, R)}  \Big|\frac{\ell N_{n}^\sss{(\delta, R)}({\bf{x}}^\sss{(\ell)}, \eps)}{n} - \theta^\sss{(\delta, R)}({\bf{x}}^\sss{(\ell)}, \eps)\Big|>\psi\Bigg)\le \exp\big( -c_{\psi, 1} n).
    \end{equation}
    Moreover, there exists $c_{\psi, 2}>0$ such that for all $n\ge 1$
    \begin{equation}\label{eq:andreis-largest}
    \Prob\big(\big||\CC_{n}^\sss{(1), (\delta, R)}|/ n-\theta_q^\sss{(\delta, R)} \big|>\psi\big)\le\exp\big(-c_{\psi, 2} n\big).
    \end{equation}
\end{lemma}
\begin{proof}
By \cref{def:approx}, the graph $\CG_n^\sss{(\delta, R)}$ contains $\underline n\le n$ vertices. Then, 
\[
\begin{aligned}
    \bigg\{ 
    \sum_{\ell\le \ell_\ast}&\sum_{  ({\bf{w}}^\sss{(\ell)},\eps)\in\mathrm{CT}_\ell(\eps, R)}  
    \Big|\frac{\ell N_{n}^\sss{(\delta, R)}({\bf{w}}^\sss{(\ell)}, \eps)}{n} - \theta^\sss{(\delta, R)}({\bf{w}}^\sss{(\ell)}, \eps)\Big|>\psi\bigg\}\\&=
    \bigg\{ 
    \sum_{\ell, ({\bf{w}}^\sss{(\ell)},\eps)}    \Big|\frac{\ell N_{n}^\sss{(\delta, R)}({\bf{w}}^\sss{(\ell)}, \eps)}{\underline n} - \theta^\sss{(\delta, R)}({\bf{w}}^\sss{(\ell)}, \eps)\frac{n}{\underline n}\Big|>\psi\frac{n}{\underline n}\bigg\}\\
    &\subseteq \bigg\{ 
    \sum_{\ell, ({\bf{w}}^\sss{(\ell)},\eps)}    \Big(\Big|\frac{\ell N_{n}^\sss{(\delta, R)}({\bf{w}}^\sss{(\ell)}, \eps)}{\underline n} - \theta^\sss{(\delta, R)}({\bf{w}}^\sss{(\ell)}, \eps)\Big|+\theta^\sss{(\delta, R)}({\bf{w}}^\sss{(\ell)}, \eps)|1-n/\underline n|\Big)>\psi\frac{n}{\underline n} \bigg\},
    \end{aligned}
    \]
    where in the second and third line we take the sums over the same sets as in the first line.
By definition of $\underline n$ in~\eqref{eq:approx-n-under}, $\underline n\ge (1-\delta)n\Prob(W\le R)$. Thus, we can take $\delta$ sufficiently small and $R$ sufficiently large such that $(n/\underline n - 1)\le \psi/4$. Since the sum over the probabilities $\theta^\sss{(\delta, R)}({\bf w}^\sss{(\ell)}, \eps)$ is at most one, we obtain  for these values of $\delta$ and $R$ that 
\begin{equation}\label{eq:andreis-n-bar-ell}
\begin{aligned}
    \Prob\bigg(&\sum_{\ell\le \ell_\ast}\sum_{  ({\bf{w}}^\sss{(\ell)},\eps)\in\mathrm{CT}_\ell(\eps, R)}  \Big|\frac{\ell N_{n}^\sss{(\delta, R)}({\bf{w}}^\sss{(\ell)}, \eps)}{n} - \theta^\sss{(\delta, R)}({\bf{w}}^\sss{(\ell)}, \eps)\Big|>\psi\bigg)\\&\le 
    \sum_{\ell\le \ell_\ast}\sum_{  ({\bf{w}}^\sss{(\ell)},\eps)\in\mathrm{CT}_\ell(\eps, R)}\Prob\bigg(  \Big|\frac{\ell N_{n}^\sss{(\delta, R)}({\bf{w}}^\sss{(\ell)}, \eps)}{\underline n} - \theta^\sss{(\delta, R)}({\bf{w}}^\sss{(\ell)}, \eps)\Big|>\frac{\psi}{2\sum_{\ell\le \ell_\ast}|\mathrm{CT}_\ell(\eps, R)|}\bigg).
    \end{aligned}
    \end{equation}
By the same argumentation, 
\begin{equation}\label{eq:andreis-n-bar-largest}
\Prob\big(\big||\CC_{n}^\sss{(1), (\delta, R)}|/ n-\theta^\sss{(\delta, R)}_q \big|>\psi\big)
\le 
\Prob\big(\big||\CC_{n}^\sss{(1), (\delta, R)}|/ \underline n-\theta^\sss{(\delta, R)}_q \big|>\psi/2\big).
\end{equation}
We now use the results from~\cite{andreis2023irg}, that (rephrased to our notation) derives among others a large-deviations principle (LDP) with speed $\underline n$ for the vector $\big(N_n^\sss{(\delta, R)}({\bf w}^\sss{(\ell)}, \eps)/n\big)_{\ell\ge 1, {\bf w}^\sss{(\ell}\in \mathrm{CT}_\ell(\eps, R)}$ and $|\CC_n^\sss{(1), (\delta, R)}|/n$. We present a corollary of this LDP and omit the full description as it would require significantly more notation. The LDP implies laws of large numbers with exponential convergence rate: there exist constants $(\eta({\bf w}^\sss{(\ell)}, \eps))_{\ell\ge 1, {\bf w}^\sss{(\ell}\in \mathrm{CT}_\ell(\eps, R)}$ such that for each ${\bf w}^\sss{(\ell)}\in \mathrm{CT}_\ell(\eps, R)$ and any $\psi'>0$
\[
\limsup_{\underline n\to\infty}\frac{1}{\underline n}\log \Prob\Big(\Big|\frac{N_n^\sss{(\delta, R)}({\bf w}^\sss{(\ell)}, \eps)}{\underline n}- \eta({\bf w}^\sss{(\ell)}, \eps)\Big|>\psi'\Big)<0.
\]
We choose $\psi' := \psi/(\ell\sum_{\ell\le \ell_\ast}|\mathrm{CT}_\ell(\eps, R))|$, and thus obtain that there exists $c'>0$ such that for all $n$ sufficiently large
\[
\begin{aligned}
    \sum_{\ell\le \ell_\ast}\sum_{  ({\bf{w}}^\sss{(\ell)},\eps)\in\mathrm{CT}_\ell(\eps, R)}\Prob\bigg(  \Big|\frac{\ell N_{n}^\sss{(\delta, R)}({\bf{w}}^\sss{(\ell)}, \eps)}{\underline n} - \ell\eta({\bf w}^\sss{(\ell)}, \eps)\Big|>\frac{\psi}{2\sum_{\ell\le \ell_\ast}|\mathrm{CT}_\ell(\eps, R)|}\bigg)\le \exp(-c'n).
    \end{aligned}
\]
Recalling~\eqref{eq:andreis-n-bar-ell},~\eqref{eq:andreis-ell} follows if we show that 
\begin{equation}\label{eq:local-l2}
\frac{\ell N_{n}^\sss{(\delta, R)}({\bf{w}}^\sss{(\ell)}, \eps)}{\underline n} \overset{\Prob}{\longrightarrow} \theta^\sss{(\delta, R)}({\bf w}^\sss{(\ell)}, \eps),\qquad\text{ as }\underline n\to\infty,
\end{equation}
as this implies that $\ell\eta({\bf w}^\sss{(\ell)}, \eps)=\theta^\sss{(\delta, R)}({\bf w}^\sss{(\ell)}, \eps)$. Showing~\eqref{eq:local-l2} is the goal of the remainder of the proof. We employ local convergence in probability for rooted vertex-marked graphs that we introduce briefly for finite mark sets (corresponding to the setting of \cref{def:approx}). We refer to~\cite{Hofbook2} for references and more elaborate descriptions.
A rooted vertex-marked graph is a couple $(G, o)$ of a graph $G$ and some, possibly random, distinguished vertex $o$ of $G$, which we call the root of $G$. We assume that the vertices $v\in V$ are given marks $m_v$ from a finite mark set $\CM$.  Let $\mathfrak{G}$ be the space of all vertex-marked rooted locally finite graphs.  We call two vertex-marked graphs $(G_1, o_1)$ and $(G_2, o_2)\in \mathfrak{G}$ isomorphic, i.e., $(G_1, o_1)\simeq(G_2, o_2)$, if there exists a bijection $\phi: V((G_1, o_1))\mapsto V((G_2, o_2))$ such that $\phi(o_1)=o_2$, $\{u,v\}$ is an edge in $(G_1, o_1)$ if and only if $\{\phi(u),\phi(v)\}$ is an edge in $(G_2, o_2)$, and $m_{u}=m_{\phi(u)}$ for all vertices $u$ in $G_1$. We write $B_G(v, r)$ for the induced subgraph of $G$ on all vertices that are at graph distance at most $r$ from a vertex $v$. Define
\[
\begin{aligned}
R\big((G_1, o_1), (G_2, o_2)\big)&:=\max\big\{r\in\N: B_{G_1}(o_1, r)\simeq B_{G_2}(o_2, r)\big\}, 
\\
d_{\mathfrak G}\big((G_1, o_1), (G_2, o_2)\big)&:=1/\big(R\big((G_1, o_1), (G_2, o_2)\big)\big).
\end{aligned}
\]   
Then, $(\mathfrak G, d_{\mathfrak G})$ constitutes a Polish space. 
We call a finite rooted graph $(G, o)$ uniformly rooted, if $o$ is chosen uniformly at random among the vertices of $G$.
We say that a sequence of uniformly rooted graphs $(G_n, o_n)_{n \ge 1}$ converges locally in probability towards $(G_{\infty}, o)$ having law $\mu$, if for every bounded and continuous function $h:\mathfrak G\mapsto\R$,
    \begin{equation}\label{eq:local-l1}
    \E\big[h(G_n, o_n) \mid G_n\big]\overset\Prob\longrightarrow\E_\mu[h(G_\infty, o)],\qquad \text{as }n\to\infty,
    \end{equation}   
    where the expectation on the left-hand side is only with respect to the uniform root $o_n$. Inhomogeneous random graphs as in \cref{def:approx} converge locally in probability to their associated branching process: it was essentially proven in~\cite{bollobas2007irg}, and a formal proof is given in \cite{Hofbook2}. We now use~\eqref{eq:local-l1} for the function specific function $h(\CG_n^\sss{(\delta, R)}, o_n)=\mathds{1}\big\{\CC_n^\sss{(\delta, R)}(o_n)\text{ has type }\theta^\sss{(\delta, R)}({\bf w}^\sss{(\ell)}, \eps)\big\}$. As the function only depends on the induced subgraph up to graph distance $\ell+1$ from the root, it follows that $h(G_1, o_1)=h(G_2, o_2)$ for any two rooted graphs that are within distance $1/(\ell+2)$ from each other. So $h(G_1, o_1)$ is continuous in the Polish space $(\mathfrak G, d_{\mathfrak G})$. It is clearly bounded, and therefore 
    \[
    \begin{aligned}
    \frac{|\{v\in [n]: \CC^\sss{(\delta, R)}(v)\text{ has type }({\bf w}^\sss{(\ell)}, \eps)\}|}{\underline n} &= \frac{1}{|\CV_n^\sss{(\delta, R)}|}\sum_{v\in \CV_n^\sss{(\delta, R)}}h(G_n, v) \\&= \E[h(\CG_n^\sss{(\delta, R)}, o_n)\mid \CG_n^\sss{(\delta, R)}]\\ 
    &\overset{\Prob}{\longrightarrow}\Prob\big(T_q^\sss{(\delta, R)}\text{ has type }({\bf w}^\sss{(\ell)}, \eps) \big)=\theta^\sss{(\delta, R)}({\bf w}^\sss{(\ell)}, \eps).
    \end{aligned}
    \]
    The number of components of type $({\bf w}^\sss{(\ell)}, \eps)$ is exactly a factor $\ell$ smaller than the left-hand side. This proves~\eqref{eq:local-l2} and therefore~\eqref{eq:andreis-ell} follows. The bound~\eqref{eq:andreis-largest} follows from~\eqref{eq:andreis-n-bar-largest} and the LDP in~\cite[Theorem 3.1]{andreis2023irg}, noting that the weak law of large numbers  $|\CC_n^\sss{(1), (\delta, R)}|/\underline n\to \theta_q^\sss{(\delta, R
)}$ was already proven in~\cite{bollobas2007irg}.
\end{proof}
The next lemma controls the relation between the IRG and the $\delta$-IRG.
Recall the definition of a type-$({\bf{w}}^\sss{(\ell)},\eps)$-component from \cref{def:type-w-comp}. By our assumption that $\eps/\delta\in\N$, the sequence $(w_i)_{i\ge 0}$ from \cref{def:type-w-comp} is a subsequence of the sequence $(z_i)_{i\ge 0}$ from \cref{def:approx}.
Define for fixed $\delta>0$
\begin{equation}
    \CA_\mathrm{reg}:=\big\{\forall i\in[(R-\underline w)/\delta]: |\CV_n[z_i, z_{i+1}) | \in\big((1-\delta) n f_W^\sss{(\delta)}(z_i), (1+\delta) n f_W^\sss{(\delta)}(z_i)\big)\big\}.\label{eq:event-reg}
\end{equation}
\begin{lemma}[Coupling of the graphs]\label{lemma:graph-coupling}
    Consider the inhomogeneous random graph from \cref{def:irg}, and its approximation from \cref{def:approx} given some $R>\underline w$ and $\delta>0$. There exists a constant $c_1>0$ such that for all $n\ge 1$
    \begin{equation}\label{eq:coupling-exp-decay}
    \Prob\big(\CA_\mathrm{reg}\big)\ge 1-\exp(-c_1n).
    \end{equation}
    On the event $\CA_\mathrm{reg}$, there exists a coupling between $\CG_n[\underline w,R)$ and $\CG_n^\sss{(\delta, R)}$ such that $\CG_n[\underline w,R)\supseteq \CG_n^\sss{(\delta, R)}$ and $|\CV_n[\underline w,R)\setminus \CV_n^\sss{(\delta, R)}|\le 2\delta n$. If under the coupling a set of vertices is a component of both $\CG_n[\underline w, R)$ and $\CG_n^\sss{(\delta, R)}$, then the type of the component, cf.\ \cref{def:type-w-comp}, coincides in the two graphs.    
    Moreover, for all $\psi, R>0$ there exists $\delta_0, c_2>0$ such that for all $\delta\in(0,\delta_0)$ and $n\ge1$, under this coupling, 
    \begin{equation}
        \Prob\big(|\CE_n[\underline w,R)\setminus \CE_n^\sss{(\delta, R)}|\ge \psi n\mid \CA_\mathrm{reg}\big)\le \exp\big(-c_2n\big).\label{eq:coupling-edges}
    \end{equation}
\end{lemma}
\begin{proof}
The number of vertices with weight in $\big((1-\delta) n f_W^\sss{(\delta)}(z_i), (1+\delta) n f_W^\sss{(\delta)}(z_i)\big)$ is distributed as $\mathrm{Bin}(n,f_W^\sss{(\delta)}(z_i))$ by $F_W$ in \cref{def:irg} and $f_W$ in \cref{def:approx}. Therefore, the exponential decay of $\Prob\big(\neg \CA_\mathrm{reg}\big)$ in~\eqref{eq:coupling-exp-decay} follows by a union bound and afterwards applying a Chernoff bound for all $i\in[(R-\underline w)/\delta]$. 

We will now construct the coupling between the two graphs. To do so, we work conditionally on the realization of the vertex set $\CV_n[\underline w, R)$ that satisfies $\CA_\mathrm{reg}$. Fix a subset $\CV_n^\mathrm{sub}[\underline w, R)\subseteq \CV_n[\underline w, R)$ such that $|\CV_n^\mathrm{sub}[\underline w+i\delta, \underline w+(i+1)\delta)|=\underline n_i$ for all $i$, which is possible by definition of $\underline n_i$ in \cref{def:approx} and $\CA_\mathrm{reg}$ in~\eqref{eq:event-reg}. Then, the vertex set $\CV_n^\mathrm{sub}[\underline w, R)$ has almost the same weight distribution as $\CV_n^\sss{(\delta, R)}[\underline w, R)$. Indeed, with $(z_i)_{i\ge 0}$ as in \cref{def:approx}, define $w^\sss{(\delta)}:=\max\{z_i: z_i\le w\}$ for $w\ge \underline w$. Then, 
    \begin{equation}\label{eq:weights-rounded-set}
    \{w_v^\sss{(\delta)}:v\in\CV_n^\mathrm{sub}[\underline w, R)\}=\{w_v: v\in \CV_n^\sss{(\delta, R)}\}.
    \end{equation}
    Thus, we set the weights of the vertices $[n]$ in $\CV_n^\sss{(\delta, R)}$ to be the set on the left-hand side for the coupling. 
    We first show that the set of the remaining vertices from $\CV_n^\mathrm{sprinkle}[\underline w,R):=\CV_n[\underline w, R)\setminus\CV_n^\mathrm{sub}[\underline w, R)=\CV_n[\underline w, R)\setminus\CV_n^\sss{(\delta, R)}[\underline w, R)$ is small. Recall $\underline n_i$ from \cref{def:approx}. On the event $\CA_\mathrm{reg}$,
\begin{equation*}
|\CV_n^\mathrm{sprinkle}[\underline w,R)|\le (1+\delta)\Prob\big(W\le R\big) n-(1-\delta)n\Prob\big(W\le R\big)<2\delta n.
\end{equation*}
Now we simultaneously construct the edges in the induced subgraph $\CG_n^{\mathrm{sub}}[\underline w, R)$ and $\CG_n^\sss{(\delta, R)}$. For each pair of vertices $u,v$ in $\CV_n^\mathrm{sub}[\underline w, R)$, we include the edge $\{u,v\}$ in both graphs with probability $p_{uv}$ from \cref{def:irg}, but then independently remove the edge from $\CG_n^\sss{(\delta, R)}$ with probability \[
1-\frac{\kappa_\sigma^\sss{(\delta, R)}( w^\sss{(\delta)}_u,  w^\sss{(\delta)}_v)}{\kappa_\sigma(w_u, w_v)}=1- \frac{\kappa_\sigma^\sss{(\delta, R)}(w_u,  w_v)}{\kappa_\sigma(w_u, w_v)}.
\] 
 Here, the equality follows from the definition of $\kappa_\sigma^\sss{(\delta, R)}$ in~\eqref{eq:kernel-adj}. The constructed edge set $\CE_n^\sss{(\delta, R)}$ has the same distribution as desired by \cref{def:approx}. Moreover, if a set of vertices is a component in both $\CG_n^\sss{(\delta, R)}$ and $\CG_n[\underline w, R)$, then the type of this set is the same by the choice of $\CV_n^\sss{(\delta, R)}$ below \eqref{eq:weights-rounded-set}, and since we assume $\eps/\delta\in\N$. It remains to prove~\eqref{eq:coupling-edges}.

The probability that an edge is present in $\CG_n^\mathrm{sub}[\underline w, R)$, but absent in $\CG_n^\sss{(\delta, R)}$, is at most
\[
\sup_{w_1, w_2\in[\underline w, R]}q\big(\kappa_{\sigma}(w_1,  w_2)- \kappa_{\sigma}^\sss{(\delta, R)}( w_1,  w_2)\big)\,\big/\,n.
\]
The kernel $\kappa_\sigma^\sss{(\delta, R)}(w_1, w_2)$ converges uniformly to $\kappa_\sigma$ for $w_1,w_2\in[\underline w,R)$ by its definition in~\eqref{eq:kernel-adj} as $\delta\downarrow 0$. Therefore, the right-hand side is at most $\psi/(4n)$ for any sufficiently small $\delta$. 
Let 
\[
\CA_\mathrm{edge\textnormal{-}1}:=\Big\{|\CE_n^{\mathrm{sub}}[\underline w, R)\setminus \CE_n^\sss{(\delta, R)}[\underline w, R)|\le (\psi/2)n\Big\}.
\]
There are at most $n^2$ edges that can be in the symmetric difference. By independence of the edge-removals, we obtain conditionally on the vertex set $\CV_n[\underline w, R)$ satisfying  $\CA_\mathrm{reg}$, 
\begin{equation}\label{eq:coupling-edge-diff-1}
\begin{aligned}
\Prob\big(\neg\CA_\mathrm{edge\textnormal{-}1} \mid \CV_n[\underline w, R)\big)\le \Prob\big(\mathrm{Bin}\big(n^2, \psi/(4 n)\big)\ge (\psi/2)n\big)=\exp\big(-\Omega(n)\big).
\end{aligned}
\end{equation}
Next, we bound the total number of edges incident to one of the vertices not in $\CV_n^{\mathrm{sub}}[\underline w, R)$, i.e., incident to at least one vertex in $\CV_n^\mathrm{sprinkle}[\underline w, R)=\CV_n[\underline w, R)\setminus \CV_n^\mathrm{sub}[\underline w,R)$. Let 
\begin{equation*}
    \CA_\mathrm{edge\textnormal{-}2}:=\Big\{\big|\{u,v\in \CV_n^\mathrm{sprinkle}[\underline w,R)\times \CV_n[\underline w,R): u\sim v\}\big| \ge 4\delta n\sup_{x,y\in[\underline w,R]}\kappa_\sigma(x,y)\Big\}.
\end{equation*}
On the event $\CA_\mathrm{reg}$, there are at most $2\delta n^2$ potential edges, each occurring with probability at most $\sup_{x,y\in[\underline w,R]}\kappa_\sigma(x,y)/n$. By another application of the Chernoff bound, 
\begin{equation*}
    \Prob\big(\CA_\mathrm{edge\textnormal{-}2}\big)\le \Prob\Big(\mathrm{Bin}\Big(2\delta n^2, \sup_{x,y\in[\underline w,R]}\kappa_\sigma(x,y)/n\Big)\ge 4\delta n\sup_{x,y\in[\underline w,R]}\kappa_\sigma(x,y)\Big) \le \exp\big(-\Omega(n)\big).
\end{equation*}
Assume that $\delta=\delta(R,\psi)$ is so small that the constant factor $4\delta \sup_{x,y\in[\underline w,R]}\kappa_\sigma(x,y)$ is at most $\psi/2$, then~\eqref{eq:coupling-edges} follows when this bound is combined with~\eqref{eq:coupling-edge-diff-1}.
\end{proof}

The next lemma compares the approximating branching process from \cref{def:approx} to the branching process from \cref{def:branching}.
\begin{lemma}[Coupling of the branching processes]\label{lemma:branching-approx}
    Consider the associated branching process of a scale-free inhomogeneous random graph from \cref{def:branching}, and its approximation from \cref{def:approx}. Fix $\ell_\ast\in\N$. For all $\psi>0$ there exists $\delta_0>0$   such that for each $\eps>0$, $\delta\in(0,\delta_0)$ such that $(\eps/\delta)\in\N$, and $R>\underline w$, 
    \begin{equation}\label{eq:branching-approx-1}
    \sum_{\ell\le \ell_\ast}\sum_{({\bf{x}}^\sss{(\ell)},\eps)\in\mathrm{CT}_\ell(\eps)} \big|\theta({\bf{x}}^\sss{(\ell)}, \eps) - \theta^\sss{(\delta, R)}({\bf{x}}^\sss{(\ell)}, \eps)\big|<\psi.
    \end{equation}
    Moreover, there exists $R_0$ such that for any $\delta\in(0,\delta_0)$ and $R>R_0$
    \begin{equation}\label{eq:branching-approx-2}
    |\theta_q - \theta^\sss{(\delta, R)}_q|<\psi.
    \end{equation}
\end{lemma}
\begin{proof}
We will couple the branching process 
$\mathrm{BP}$ and $\mathrm{BP}^\sss{(\delta, R)}$ until the branching processes die out or have size larger than $\ell_\ast$. On a good event that we construct (and which holds with probability arbitrary close to 1 by choosing $\delta$ small), the two branching processes have the same type. A branching process has exactly one type in $\mathrm{CT}_\ell(\eps)$, so the sum on the left-hand side in~\eqref{eq:branching-approx-1} is at most the probability that the branching processes have a different type.  \smallskip 

\emph{Approximating infinite-type branching process. } We start with an observation. By definition of $f^\sss{(\delta)}_W$ and $\kappa_\sigma^\sss{(\delta, R)}$ in \cref{def:approx}, the approximating \emph{finite}-type branching process from \cref{def:approx} can be coupled with the following infinite-type branching process that we denote by $\widetilde{\mathrm{BP}}^\sss{(\delta, R)}$. The root of the branching process has type $\widetilde W_\varnothing$ with distribution satisfying $\Prob(\widetilde W_\varnothing > w)=\Prob\big(W>w \mid W<R\big)$. Here, $W$ has distribution $F_W$ defined in \cref{def:irg}. In each generation, each particle $v$ of type $w_v$ gives independently birth to new particles according to a Poisson point process (PPP) on $[\underline w,\infty)$ with intensity $q\kappa_\sigma(w_v, w)\rd F_W(w)$. The atoms in the union of these PPPs form the vertex types of the vertices in the next generation.  

The coupling with the branching process $\mathrm{BP}^\sss{(\delta, R)}$ from \cref{def:approx} follows from rounding down each type $W$ to the largest $z_i\le W$ with $z_i$ from \cref{def:approx}, and by coupling the PPP determining the offspring (which has constant intensity in sets $[z_i, z_i+\delta)\times[z_j, z_j+\delta)$) with the Poisson random variables determining the offspring in \cref{def:approx}. This rounding operation does not affect the type of the total progeny of the branching process, since we assumed that $\eps/\delta\in\N$, see Definitions~\ref{def:type-w-comp} and~\ref{def:approx}.

Now, we sample the original branching process $\mathrm{BP}$ from \cref{def:irg}. If it contains a vertex with weight at least $R$ in one of the first $\ell$ generations, we say that the coupling has failed. If $R=R(\psi)$ is sufficiently large, this occurs with probability at most $\psi/2$. 
The branching process $\widetilde{\mathrm{BP}}^\sss{(\delta, R)}$  can be created from BP as follows. Consider the  particles of the PPP determining the offspring of the root. 
We remove each particle with probability $1-\kappa_\sigma^\sss{(\delta, R)}(W_j, W_\varnothing)/\kappa_\sigma( W_j, W_\varnothing)$ independently of the rest. The obtained offspring has the same distribution as in $\widetilde{\mathrm{BP}}^\sss{(\delta, R)}$.
If none of the at most $\ell_\ast-1$ particles is removed, the coupling step of generation zero is successful. By the definition of $\kappa^\sss{(\delta, R)}$ in~\eqref{eq:kernel-adj}, we obtain for $\delta=\delta(\psi,\ell_\ast)$ sufficiently small that
\[
1-\kappa_\sigma^\sss{(\delta, R)}(W_j, W_\varnothing)/\kappa_\sigma( W_j, W_\varnothing) < \psi/(2\ell_\ast^2).
\] Hence, with probability at most $\psi/(2\ell_\ast)$, one of the particles is removed. Iterating this procedure at most $\ell_\ast$ times yields that no particle is removed with probability at least $1-\psi/2$. Thus, the coupling is successful with probability at least $1-\psi$, and~\eqref{eq:branching-approx-1} follows.

We turn to~\eqref{eq:branching-approx-2}. Let $\theta_q^\sss{(0,R)}$ denote the survival probability of a slightly modified associated branching process compared to \cref{def:branching}: the new particles are formed by a PPP on $[\underline w,\infty)$ with intensity $\ind{w\le R}q\kappa_\sigma(w_v, w)\rd F_W(w)$. 
By the triangle inequality 
\begin{equation}\label{eq:triangle}
|\theta_q-\theta_q^\sss{(\delta, R)}|\le |\theta_q-\theta_q^\sss{(0, R)}| + |\theta_q^\sss{(0,R)}-\theta_q^\sss{(\delta, R)}|
\end{equation}
Since $\ind{w\le R}\kappa_\sigma(w_v, w)\uparrow \kappa_\sigma(w_v, w)$ as $R\to\infty$, it follows by \cite[Theorem 6.3]{bollobas2007irg} that $\theta_q^\sss{(0,R)}\uparrow \theta_q$. Fix $R$ sufficiently large that~\eqref{eq:branching-approx-1} holds, and moreover, $|\theta_q-\theta_q^\sss{(0, R)}|<\psi/2$. The same argument works for the second term on the right-hand side in~\eqref{eq:triangle}: by definition of $\kappa_\sigma^\sss{(\delta, R)}$ in~\eqref{eq:kernel-adj}, $\kappa_\sigma^\sss{(\delta, R)}(w_u, w_v)\uparrow \ind{w_u, w_v\le R}\kappa_\sigma(w_u, w_v)$ as $\delta$ tends to 0, so  $|\theta_q^\sss{(\delta, R)}- \theta_q^\sss{(0,R)}|\to 0$ as $\delta$ tends to 0. Hence, we  can choose $\delta$ sufficiently small so that~\eqref{eq:branching-approx-1} holds, and so that $|\theta_q^\sss{(0,R)}-\theta_q^\sss{(\delta, R)}|\le \psi/2$. 
\end{proof}
 
Using Lemmas~\ref{lemma:andreis}---\ref{lemma:branching-approx} one can prove the bound~\eqref{eq:lemma-ell-lower-R} in \cref{lemma:ell-trunc-n}. We postpone the proof to the end of the following subsection, after we analyzed the effect of the vertices with weight in $[R, \phi n)$ so that we can immediately prove the other bound in \cref{lemma:ell-trunc-n}.

\subsection{Vertices with intermediate weight} \label{sec:intermediate}
The next lemma shows that the amount of edges from vertices with weight in $[R,\phi n)$ is small if $R$ is large and $\phi$ is small. As a consequence, such vertices do not significantly change the number of type-$({\bf{x}}^\sss{(\ell)},\eps)$ components in $\CG_n[\underline w,\phi n)$ when compared to $\CG_n[\underline w, R)$. We write $\mathrm{deg}_v[a,b)$ for the number of neighbors of $v$ with weight in the interval $[a,b)$.
\begin{lemma}[Excess edges]\label{lemma:total-high-degree}
    Consider an inhomogeneous scale-free random graph as in \cref{def:irg}. For all constants $\psi, C>0$ there exist constants $R_0, \phi_0>0$ such that for all $R\ge R_0$ and $\phi\in(0,\phi_0)$, as $n\to\infty$,
    \begin{equation}\label{eq:lemma-total-high}
    \Prob\bigg(\sum_{v\in\CV_n[R, \phi n)}\mathrm{deg}_v[\underline w, R)> \psi n\bigg)=o(n^{-C}).
    \end{equation}
\end{lemma}
\begin{proof}
Let $B(u,v), u,v \in [n]$ be a collection of Bernoulli random variables with success probability $p_{uv}$ defined in~\eqref{eq:conn-prob}.
Observe that
\begin{equation}\nonumber
    \sum_{v\in\CV_n[R, \phi n)}\mathrm{deg}_v[\underline w, R) \stackrel{d}{=} 
    \sum_{u\in\CV_n[R, \phi n)} \sum_{v\in\CV_n[\underline w, R]} B(u,v) =: E_n(\phi, R). 
\end{equation}
The mean of the right-hand side conditional upon the weights $W_u, u\in [n]$,
equals 
\begin{equation}\label{eq:mn}K_n(R,\phi) := \frac{q}{n} \sum_{u\in\CV_n[R, \phi n)}W_u \sum_{v\in\CV_n[\underline w, R]}   W_v^\sigma\end{equation}
(making $\phi$ smaller if needed so that $\phi n R \leq n$ and we can drop the minimum operator in the definition of $p_{uv}$). Next, take  $\varepsilon>0$ and bound using the triangle inequality 
\begin{align}
    \Prob(E_n(R, \phi) >  \psi n)  & =   \Prob\big(\{E_n(R, \phi) >  \psi n
    \}\cap\{ |E_n(R, \phi)  - K_n(R, \phi)| \leq \varepsilon K_n( R, \phi)\}\big) \nonumber\\
    & \quad +  \Prob\big(\{E_n(R, \phi) >  \psi n\} \cap \{ |E_n(R, \phi) - M_n(R, \phi)| >\varepsilon K_n(R,\phi)\}\big) \nonumber\\
    &\le \Prob\big( K_n(R, \phi) >  \tfrac{\psi}{1+\eps} n\big) + \Prob\big(|E_n(R, \phi) - K_n(R, \phi)| >\varepsilon K_n(R,\phi)\big).\label{twoterms}
\end{align}
As a first step, we show that the first term on the right-hand side is of the order $o(n^{-C})$ if $R$ is sufficiently large and $\phi$ sufficiently small. 
Define $G(R):=2\E[W^\sigma\ind{W\le R}]$, and let 
\begin{equation}\label{eq:GR}
   g(R):= \frac1{2\E[W^\sigma \ind{W\leq R}]\E[W \ind{W \geq R}] }.
\end{equation}
Consequently, $\Prob\big( K_n(R, \phi) >  \psi n/(1+\eps)\big)$ is bounded from above by
\begin{equation}\label{eq:two-terms-gr}
     \Prob\bigg(\sum_{v\in [n]} W_v^\sigma \ind{W_v\leq R} > G(R)n\bigg) +\Prob\bigg(\sum_{v\in [n]} W_v \ind{W_v\in [R, \phi n)} > n \E[W \ind{W \geq R}] \frac{\psi  g(R)}{1+\eps}\bigg).
\end{equation}
The first term decreases exponentially in $n$, regardless of the choice of $R$, due to Cram\'ers bound. 
For the second term, we first show that $g(R)\to\infty$. Indeed, using Potter's bound, the product of the truncated expectations in~\eqref{eq:GR} is at most of order $O(R^{\max(\sigma-\alpha, 0)+\varepsilon}R^{\varepsilon-(\alpha-1)})$ for any $\varepsilon>0$ as $R\to\infty$, which is of order $o(1)$ when $\varepsilon\in(0, \min(2\alpha-1-\sigma, \alpha-1))$ (note that this interval is non-empty by the assumptions $\alpha>1$ and $\sigma<2\alpha-1$ in Definition~\ref{def:irg}). Thus, $g(R)\to\infty$, and hence we may choose $R_0$ such that $\psi g(R)/(1+\eps) \geq 2$
if $R\geq R_0$. Next, we apply a bound for sums of truncated heavy-tailed random variables (see Lemma \ref{lem-truncated}) to conclude that, given $C>0$, we can pick $\phi_0$ suitably small so that
 the second term in~\eqref{eq:two-terms-gr} is of order $o(n^{-C})$. Hence, the first term in~\eqref{twoterms} is of order $o(n^{-C})$. 

 We turn to the second term in~\eqref{twoterms}. Conditional on $(W_u)_{u\in [n]}$, the $B(u,v), u<v$, are independent. 
Therefore,
we can apply Lemma \ref{lem-bernoulli} to obtain that
\begin{equation}\nonumber
      \Prob\big( |E_n(R, \phi) - K_n(R, \phi) | > \eps K_n (R, \phi) \mid W_1,\ldots,W_n\big) \leq 2\re^{-K_n J(\eps)}
\end{equation}
almost surely, for some $J(\eps) >0$. We distinguish whether $K_n(R, \phi) \leq \zeta n$ for some $\zeta>0$: 
\begin{align}
\label{bound-EnEn}
 \Prob\big(|E_n&(R, \phi) - K_n(R, \phi)| >\varepsilon K_n(R, \phi)\big)\nonumber \\ &\leq
    \Prob\big(\{ |E_n(R, \phi) - K_n(R, \phi)| >\varepsilon K_n(R, \phi)\} \cap\{K_n(R, \phi) > \zeta n\}\big)
     + \Prob(K_n(R, \phi) \leq \zeta n) \nonumber\\
    & \leq 2\re^{-\zeta n J(\eps)} + \Prob(K_n(R, \phi) \leq \zeta n).
\end{align}
Recall the definition of $M_n$ from~\eqref{eq:mn}. By the union bound,
\begin{equation}\nonumber
     \Prob(K_n \le \zeta n) \leq \Prob\bigg( q \sum_{u\in [n]} 
     W_u \ind{W_u \in [R,\phi n]} \le \sqrt{\zeta} n \bigg)+\Prob\bigg(  \sum_{v\in[n]} W_v^\sigma \ind{W_v \leq R} \le \sqrt{\zeta} n\bigg).
 \end{equation}
 Applying Cram\'er's theorem to each of the terms on the right-hand side yields that the right-hand side tends to zero exponentially fast for $\zeta$ sufficiently small depending on $R$ and $\phi$. Therefore both terms in~\eqref{bound-EnEn} decay exponentially fast, and also the second term in~\eqref{twoterms} is of order $o(n^{-C})$. 
\end{proof}

\subsection{Size-$\ell$ components in the graph without hubs}\label{sec:size-ell-without-hubs}
We are ready to prove Lemmas~\ref{lemma:ell-trunc-n}--\ref{lemma:giant-trunc-lower}. 
\begin{proof}[Proof of \cref{lemma:ell-trunc-n}]
We start with the proof of~\eqref{eq:lemma-ell-lower-R}.
We use the coupling of the IRG restricted to weights in $[\underline w, R)$, and its $\delta$-approximation from \cref{lemma:graph-coupling}. Fix $\psi_{\ref{lemma:ell-trunc-n}}$ where $\delta$ is sufficiently small that \cref{lemma:graph-coupling} holds with $\psi_{\ref{lemma:ell-trunc-n}}/4$, and~\eqref{eq:branching-approx-1} in \cref{lemma:branching-approx} hold with $\psi=\psi_{\ref{lemma:ell-trunc-n}}/(16\ell_\ast)$, and $2\delta\le \psi_{\ref{lemma:ell-trunc-n}}/8$.  
We bound
\begin{equation}\nonumber \begin{aligned}
    \sum_{\ell, ({\bf{x}}^\sss{(\ell)},\eps)}|\ell N_n({\bf{x}}^\sss{(\ell)}, \eps, R) / n - \theta({\bf{x}}^\sss{(\ell)}, \eps)| &\le \sum_{\ell, ({\bf{x}}^\sss{(\ell)},\eps)} \ell\cdot |N_n({\bf{x}}^\sss{(\ell)}, \eps, R)  - N_n^\sss{(\delta, R)}({\bf{x}}^\sss{(\ell)}, \eps, R)|/n \\&\hspace{15pt}+ \sum_{\ell, ({\bf{x}}^\sss{(\ell)},\eps)}
    |\ell N_n^\sss{(\delta, R)}({\bf{x}}^\sss{(\ell)}, \eps, R)/n - \theta^\sss{(\delta, R)}({\bf{x}}^\sss{(\ell)}, \eps)|  \\&\hspace{15pt}
    +\sum_{\ell, ({\bf{x}}^\sss{(\ell)},\eps)} 
    |\theta^\sss{(\delta, R)}({\bf{x}}^\sss{(\ell)}, \eps) - \theta({\bf{x}}^\sss{(\ell)}, \eps)  |.
\end{aligned}\end{equation}
Here and in the remainder of the proof, all sums are for $\ell\le \ell_\ast$ and $({\bf{x}}^\sss{(\ell)}, \eps)\in \mathrm{CT}_\ell(\eps, R)$. For readability we omit this in the subscripts of the sums. The third term on the right-hand side is at most $\psi_{\ref{lemma:ell-trunc-n}}/4$ by \cref{lemma:branching-approx}. With probability tending to one exponentially fast, the second term is also at most  $\psi_{\ref{lemma:ell-trunc-n}}/4$ by \cref{lemma:andreis}. The coupling is successful with probability tending to one exponentially fast by~\eqref{eq:coupling-exp-decay}. Thus, it remains to show on the event that the coupling is successful that the first sum is at most $\psi_{\ref{lemma:ell-trunc-n}}/2$ with probability tending to one exponentially fast, i.e., for some $c>0$
\begin{equation}\nonumber 
    \Prob\Bigg(\sum_{\ell, ({\bf{x}}^\sss{(\ell)},\eps)} \ell |N_n({\bf{x}}^\sss{(\ell)}, \eps, R)  - N_n^\sss{(\delta, R)}({\bf{x}}^\sss{(\ell)}, \eps, R)|> (\psi_{\ref{lemma:ell-trunc-n}}/4)n \, \Big|\, \CA_\mathrm{reg}\Bigg)\le \exp(-cn).
\end{equation}
On the event $\CA_\mathrm{reg}$, we add the vertices from $\CV_n[1, R)\setminus\CV_n^\sss{(\delta, R)}$ and the edges from $\CE_n[1, R)\setminus\CE_n^\sss{(\delta, R)}$, and evaluate how the sum changes. After adding the vertices, the sum of the differences changes by at most $|\CV_n[1, R)\setminus\CV_n^\sss{(\delta, R)}|$. We add the edges iteratively: after each added edge, the summands corresponding to the two incident vertices changes its value by the size of the respective components. So the added edges change the sum in total by at most $2\ell_\ast|\CE_n[1, R)\setminus\CE_n^\sss{(\delta, R)}|$. We conclude that 
\begin{equation}\nonumber
\begin{aligned}
    \Prob\Bigg(\sum_{\ell, ({\bf{x}}^\sss{(\ell)},\eps)} |N_n({\bf{x}}^\sss{(\ell)}, \eps, R)  &- N_n^\sss{(\delta, R)}({\bf{x}}^\sss{(\ell)}, \eps, R)|> (\psi_{\ref{lemma:ell-trunc-n}}/4)n \, \Big|\, \CA_\mathrm{reg}\Bigg)\\
    &\le 
    \Prob\Big(|\CV_n[1, R)\setminus\CV_n^\sss{(\delta, R)}|+ 2\ell_\ast|\CE_n[1, R)\setminus\CE_n^\sss{(\delta, R)}|> (\psi_{\ref{lemma:ell-trunc-n}}/4)n \, \big|\, \CA_\mathrm{reg}\Big).
    \end{aligned}
\end{equation}
Since on $\CA_\mathrm{reg}$, $|\CV_n[1, R)\setminus\CV_n^\sss{(\delta, R)}|\le 2\delta n\le (\psi_{\ref{lemma:ell-trunc-n}}/8)n$ by \cref{lemma:graph-coupling} and the choice of $\delta$ at the beginning, the right-hand side is at most $\Prob\big(|\CE_n[1, R)\setminus\CE_n^\sss{(\delta, R)}|> (\psi_{\ref{lemma:ell-trunc-n}}/16\ell_\ast)n \, \big|\, \CA_\mathrm{reg}\big)$. This probability decays exponentially fast by \cref{lemma:graph-coupling} by the choice of $\psi$ at the beginning of the proof. This finishes the proof of \eqref{eq:lemma-ell-lower-R}.

We next turn to the bound in~\eqref{eq:lemma-ell-lower} for which we take a similar approach, but start from~\eqref{eq:lemma-ell-lower-R} instead, assuming that this bound holds with $\psi$ replaced by $\psi/2$. The probability on the left-hand side in~\eqref{eq:lemma-ell-lower} is increasing in $R$. Therefore, we may assume without loss of generality that $R$ is at least the large constant $R_0$ from  Lemma~\ref{lemma:total-high-degree} for the same value $C$ and $\psi$ replaced by $\psi/(2\ell_\ast)$. Let $\phi>0$ also be  from Lemma~\ref{lemma:total-high-degree}.

  We consider the number of type-$({\bf{x}}^\sss{(\ell)},\eps)$ for all component types $({\bf{x}}^\sss{(\ell)},\eps)\in \mathrm{CT}_\ell(\eps, R)\in \mathrm{CT}_\ell(\eps, R)$ in $\CG[\underline w, R)$ and iteratively add the edges incident to the vertices with weight in $[R, \phi n)$ to this graph.  None of the vertices in $\CV_n[R, \phi n)$ can be in a component of type $({\bf{x}}^\sss{(\ell)},\eps)\in \mathrm{CT}_\ell(\eps, R)$, since by definition all vertices in a type-$({\bf{x}}^\sss{(\ell)},\eps)$ component have weight at most $R$ by \cref{def:type-w-comp}. Therefore each edge in $\CE_n[\underline w, \phi n)\setminus\CE_n[\underline w, R)$ never increases the number of type-$({\bf{x}}^\sss{(\ell)},\eps)$ components, and decreases the number of at most one type $({\bf{x}}^\sss{(\ell)},\eps)\in \mathrm{CT}_\ell(\eps, R)$ by at most one. Therefore, 
  \[
  \begin{aligned}
  \sum_{\ell\le \ell_\ast}\sum_{({\bf x}^\sss{(\ell)}, \eps)\in \mathrm{CT}_\ell(\eps, R)}\hspace{-15pt}&|\ell N_n({\bf x}^\sss{(\ell)}, \eps, \phi n)-n\theta({\bf x}^\sss{(\ell)}, \eps)| \\
  &\hspace{-25pt}\le 
  \sum_{\ell\le \ell_\ast}\sum_{({\bf x}^\sss{(\ell)}, \eps)\in \mathrm{CT}_\ell(\eps, R)}\hspace{-15pt}\Big(\ell\cdot|N_n({\bf x}^\sss{(\ell)}, \eps, \phi n)-N_n({\bf x}^\sss{(\ell)}, \eps, R)| 
  +|\ell N_n({\bf x}^\sss{(\ell)}, \eps, R)- n\theta({\bf x}^\sss{(\ell)}, \eps)|\Big)\\
  &\hspace{-25pt}\le \ell_\ast\cdot\hspace{-10pt}\sum_{v\in \CV_n[R, \phi n)}\mathrm{deg}_v[\underline w, R)+ \sum_{\ell\le \ell_\ast}\sum_{({\bf x}^\sss{(\ell)}, \eps)\in \mathrm{CT}_\ell(\eps, R)}\hspace{-15pt}|\ell N_n({\bf x}^\sss{(\ell)}, \eps, \phi n)- n\theta({\bf x}^\sss{(\ell)}, \eps)|.
  \end{aligned}
  \]
  We bound
\begin{align*}
\Prob&\Bigg(\sum_{\ell\le \ell_\ast, ({\bf{x}}^\sss{(\ell)},\eps)\in \mathrm{CT}_\ell(\eps, R)}\hspace{-15pt}|\ell N_n({\bf{x}}^\sss{(\ell)}, \eps, \phi n) / n - \theta({\bf{x}}^\sss{(\ell)}, \eps)|>\psi\bigg) \\
&\hspace{-5pt}\le 
\Prob\Bigg(\sum_{\ell\le \ell_\ast, ({\bf{x}}^\sss{(\ell)},\eps)\in \mathrm{CT}_\ell(\eps, R)}\hspace{-29pt}|\ell N_n({\bf{x}}^\sss{(\ell)}, \eps, R) / n - \theta({\bf{x}}^\sss{(\ell)}, \eps)|>\psi/2\Bigg) + \Prob\Bigg(\sum_{v\in \CV_n[R, \phi n)}\hspace{-12pt}\mathrm{deg}_v[\underline w, R)>\psi/(2\ell_\ast)\Bigg).
\end{align*}
  The first term decays exponentially by~\eqref{eq:lemma-ell-lower-R}, while the second term is $o(n^{-C})$ by \cref{lemma:total-high-degree} and our assumptions on $R$ and $\phi$.
\end{proof}

We end this section by proving \cref{lemma:giant-trunc-lower}.
\begin{proof}[Proof of \cref{lemma:giant-trunc-lower}]
    Let $\CA_\mathrm{reg}$ be as in~\eqref{eq:event-reg}. Assume that $\delta$ and $R$ are such that Lemmas~\ref{lemma:andreis}--\ref{lemma:branching-approx} apply with $\psi=(\theta_q-\rho)/2$.  By \cref{lemma:graph-coupling} we can couple the graphs $\CG_n^\sss{(\delta, R)}$ and $\CG_n[\underline w, R)$ on the event $\CA_\mathrm{reg}$ such that $\CG_n[\underline w, R)\supseteq \CG_n^\sss{(\delta, R)}$. Therefore, the largest component of $\CG_n[\underline w, R)$ has at least the size of the largest component of $\CG_n^\sss{(\delta, R)}$ provided that the coupling is successful, which happens with probability increasing to $1$ exponentially fast.  Therefore, 
    \[
    \Prob\big(|\CC_{n}^\sss{(1)}[\underline w, R)|<\rho n\big) \le \Prob\big(\big\{|\CC_{n}^\sss{(\delta, R), (1)}|<(\theta_q-\psi)n\big\}\cap \CA_\mathrm{reg}\big) + \exp(-\Theta(n)).
    \]
    By \cref{lemma:andreis}, $\Prob\big(|\CC_{n}^\sss{(\delta, R), (1)}|<(\theta_q^\sss{(\delta, R)}-\psi/2)n\big)$ decays exponentially fast. Moreover, by \cref{lemma:branching-approx}, $|\theta_q^\sss{(\delta, R)}-\theta_q|\le \psi/2$. Thus, \cref{lemma:giant-trunc-lower} follows.
\end{proof}

\section{The graph with hubs}\label{sec:reduction}
In this section we prove our main result, the upper tail $\Prob(|\CC_n^\sss{(1)}|>\rho n)$ stated in \cref{thm:upper-tail}. After a technical lemma, we formalize that  the probability that the giant increases above $\rho n$ is of smaller order than $(n\Prob\big(W>n))^{\lceil\mathrm{hubs}(\rho, q)\rceil}$ if the number of hubs is not equal to $\lceil\mathrm{hubs}(\rho, q)\rceil$ or if the weights are not contained in $\CY_{\rho-\delta, q}$, leading to a lower bound for a suitably small $\delta$. Afterwards, we prove that the presence of exactly $\lceil\mathrm{hubs}(\rho, q)\rceil$ hubs with weights in $\CY_{\rho+\delta, q}$ increases the size of the giant to above $\rho$ to establish an upper bound. At the end of the section we show that we can take the limit of $\delta\to0$ to finish the proof of \cref{thm:upper-tail}.

In the next lemma that we state, we approximate the expectation in~\eqref{eq:y-rho-set} with a version where the weights are discretized and truncated, so that the number of types is finite. 
Define 
\begin{equation}
    \bar P ({\bf{w}}^\sss{(\ell)}, {\bf{y}}^{(h)}):= \prod_{j_1=1}^\ell \prod_{j_2=1}^h \Big(1-q\big(x_{j_1}^\sigma y_{j_2}\wedge 1\big)\Big).\label{eq:p-bar}
\end{equation}
$\bar P$  corresponds to the probability that a component whose vertices $\{v_1, \ldots, v_\ell\}$ have weight exactly $\{w_1,\ldots, w_\ell\}$, do not connect to one of the $h$ hubs that have weight exactly $\{y_{1}n, \ldots, y_hn\}$, see the connection probability $p_{uv}$ in~\eqref{eq:conn-prob}. Let  $\theta({\bf{ w}}^\sss{(\ell)},\varepsilon)$ be the probability that the progeny of the associated branching process has type $({\bf{ w}}^\sss{(\ell)},\varepsilon)$, as defined in Definition~\ref{def:type-w-comp}. In the next lemma we  implicitly show  that
\begin{equation*}
    \E\bigg[\bar P\big((W_x)_{x\in T_q}, {\bf{y}}\big)\bigg]=\E\bigg[\prod_{\substack{x\in T_q,\\ i\in[h] }}\big(1-q\cdot(y_iW_x^\sigma\wedge 1)\big)\bigg]=\lim_{\varepsilon\downarrow 0}\sum_{\ell=1}^\infty \sum_{{\bf{w}}^\sss{(\ell)}\in \mathrm{CT}_\ell(\eps)}\hspace{-10pt}\bar P ({\bf{w}}^\sss{(\ell)}, {\bf{y}}^{(h)})\cdot \theta({\bf{ w}}^\sss{(\ell)},\varepsilon).
\end{equation*}
We will control the rate of convergence in the above limit uniformly in the vector ${\bf{y}}^{(h)}$. We recall from \cref{def:branching} that $\theta_q$ denotes the probability that the associated branching process survives. 
\begin{lemma}[Effect of truncation and discretization]\label{lemma:approx-non-connections}
Consider an inhomogeneous scale-free random graph as in \cref{def:irg}. For all constants $\psi, h>0$, there exist constants $\eps, \ell_\ast, \bar R>0$ such that
    \begin{equation}\label{eq:non-conn-exp-lower}
        \sum_{\ell=1}^{\ell_\ast}\sum_{{\bf{w}}^\sss{(\ell)}\in\mathrm{CT}_\ell(\eps, \bar R)}\theta({\bf{w}}^\sss{(\ell)}, \eps)
        >1-\theta_q-\psi,
    \end{equation}
    and
    \begin{equation}\label{eq:non-conn-upper}
    \sup_{{\bf{y}}^{(h)}>\phi{\bf{1}}^{(h)}}\bigg|\E\Big[\bar P\big((W_x)_{x\in T_q}, {\bf{y}}^{(h)}\big)\Big]\,-\, \sum_{\ell=1}^{\ell_\ast}\sum_{{\bf{w}}^\sss{(\ell)}\in\mathrm{CT}_\ell(\eps, \bar R)}\bar P({\bf{w}}^\sss{(\ell)}, {\bf{y}}^{(h)})\theta({\bf{w}}^\sss{(\ell)}, \eps) \bigg|<\psi,
    \end{equation}
    and 
        \begin{align}
        \hspace{-9pt}\sup_{\substack{{\bf{y}}^{(h)}>\phi{\bf{1}}^{(h)}, \\{\bf{w}^\sss{(\ell)}}\in\mathrm{CT}_\ell(\eps,\bar R)}} \left|
        \frac{\Prob\Big(\CC\!\nsim\! \CV_n[\phi n,\infty) \,\big|\, \CC\,\, \text{\normalfont has type }({\bf{w}}^\sss{(\ell)}, \eps), \CV_n[\phi n,\infty)=\{y_1 n, \ldots, y_hn\}\Big)}{\bar P ({\bf{w}}^\sss{(\ell)}, {\bf{y}}^{(h)})} - 1
        \right|
    <\psi. \label{eq:non-conn-component} 
    \end{align}
\end{lemma}
We postpone the proof to the appendix on page~\pageref{proof:truncation-effect}. We proceed to a proposition that {we will use to prove a lower bound on $\Prob(|\CC_n^\sss{(1)}|>\rho n)$.} Let $\CY_{\rho, q}(h)\subseteq [0,\infty)^h$ as in~\eqref{eq:y-rho-set}. We will slightly abuse notation, and identify the set of vertices $\CV_n[\phi n,\infty)$ with the weights of the vertices in the interval $[\phi n,\infty)$. Thus, if this set has size $h$, we
we write $\CV[\phi n,\infty)\in n\cdot \CY_{\rho,q}(h)$ if $(w_v/n)_{v\in\CV[\phi n,\infty)}\in\CY_{\rho, q}(h)$. 
\begin{proposition}[Absence of hubs implies no large giant]\label{prop:lower-bound}
Consider an inhomogeneous scale-free random graph as in \cref{def:irg}. Fix a constant $\rho\in(\theta_q, 1)$. There exists a constant $\phi_0>0$ such that for any $\phi\in(0,\phi_0)$ and $r\in(\theta_q,\rho)$, with $h=\lceil\mathrm{hubs}(r, q)\rceil$, as $n\to\infty$,
\begin{equation}\label{eq:prop-upper}
\begin{aligned}
\Prob\Big(|\CC_n^\sss{(1)}|> \rho n \, \big|\, \neg\big\{|\CV_n[\phi n,\infty)|=h, \CV_n[\phi n,\infty)\in n\cdot\CY_{r, q}(h)\big\}\Big)= O\big(\big(n\Prob(W_1>n)\big)^{h+1}\big).
\end{aligned}
\end{equation}
\end{proposition}
\begin{proof}
We describe the idea of the proof. We will reveal the graph $\CG_n$ in two stages: first we consider the graph $\CG_n[\underline w,\phi n)$, in which the number of size-$\ell$ components of finitely many types concentrates using \cref{lemma:ell-trunc-n}. In this graph, the complement of the giant has size approximately $(1-\theta_q)n$. Then we `add' the vertices of weight at least $\phi n$ to the graph, which we call hubs. If the number of hubs is smaller than $\lceil\mathrm{hubs}(r, q)\rceil$, or the set of their weights is not contained in $n\cdot\CY_{r, q}$ while there are exactly $\lceil\mathrm{hubs}(r, q)\rceil$ many hubs, then enough components do \emph{not} connect to one of the hubs with sufficiently high probability. Consequently, the complement of the giant may become much smaller than $(1-\theta_q)n$, but remains larger than $(1-\rho)n$ as long as $r<\rho$. 
The bound on the right-hand side in~\eqref{eq:prop-upper} corresponds to the event that there are more than $\lceil\mathrm{hubs}(r, q)\rceil$ hubs.

We now make this formal. 
Fix $r\in(\theta_q, \rho)$, let
$\psi=\psi(r), \varepsilon=\varepsilon(\psi), \phi=\phi(\psi)$  be sufficiently small constants, and let $\bar R=\bar R(\psi), \ell_\ast=\ell_\ast(\psi)$ be two large constants.
We first define subevents for the event on which we conditioned in~\eqref{eq:prop-upper}:
\begin{align}
        \CA_\mathrm{hubs, 1}^\le&:=\big\{ |\CV_n[\phi n, \infty)|=h, \CV_n[\phi n,\infty)\notin n\cdot\CY_{r, q}(h)\big\}\label{eq:hubs1-ub},\\
    \CA_\mathrm{hubs, 2}^\le&:=\big\{|\CV_n[\phi n, \infty)|< h\big\}, \label{eq:hubs2-ub}\\
    \CA_\mathrm{hubs, 3}^\le&:=\big\{|\CV_n[\phi n, \infty)|>\lceil\mathrm{hubs}(r, q)\rceil\big\},\label{eq:hubs3-ub}\\    
    \CA_\mathrm{hubs, 1\textnormal{-}2}^\le &:=\CA_\mathrm{hubs, 1}^\le\cup \CA_\mathrm{hubs, 2}^\le, \\
    \CA_\mathrm{hubs, 1\textnormal{-}3}^\le&:=\CA_\mathrm{hubs, 1}^\le\cup \CA_\mathrm{hubs, 2}^\le\cup\CA_\mathrm{hubs, 3}^\le.\label{eq:hubs-ub}
\end{align}
Observe that $\CA^\le_\mathrm{hubs, 1\textnormal{-}3}$ corresponds to the event on which we conditioned in~\eqref{eq:prop-upper}. 
We first bound
\begin{equation}\nonumber
    \Prob\big(|\CC_n^\sss{(1)}|>\rho n \mid \CA_\mathrm{hubs, 1\textnormal{-}3}^\le\big) \le \Prob\big(\{|\CC_n^\sss{(1)}|>\rho n \}\cap (\neg \CA_\mathrm{hubs, 3}^\le)  \mid \CA_\mathrm{hubs, 1\textnormal{-}3}^\le\big) + \Prob\big(\CA_\mathrm{hubs, 3}^\le  \mid \CA_\mathrm{hubs, 1\textnormal{-}3}^\le\big).
\end{equation}
Writing out the conditional probabilities and applying elementary operations, we obtain that
\begin{equation}\label{eq:ub-below-hubs}
    \Prob\big(|\CC_n^\sss{(1)}|>\rho n \mid \CA_\mathrm{hubs, 1\textnormal{-}3}^\le\big) \le \Prob\big(|\CC_n^\sss{(1)}|>\rho n  \mid  \CA_\mathrm{hubs, 1\textnormal{-}2}^\le\big) + \frac{\Prob\big(\CA_\mathrm{hubs, 3}^\le\big)}{\Prob\big(\CA_\mathrm{hubs, 2}^\le\big)}.
\end{equation}
The denominator of the second term tends to 1 as $n\to\infty$, while the numerator is of order $O\big(\big(n\Prob(W_1>n)\big)^{h+1}\big)$ by~\eqref{eq:weight-dist}. Thus, it suffices to show in the remainder that the first term on right-hand side is of smaller order.

We introduce more notation.
For each $\ell\le \ell_\ast$,
    we write $M_{n}({\bf{w}}^\sss{(\ell)},\eps)$ for the number of components of type ${\bf{w}}^\sss{(\ell)}$ in the induced subgraph $\CG_n[\underline w,\phi n)$ that are \emph{not} connected by an edge to the hubs in $\CG_n$, $h:=|\CV_n[\phi n,\infty)|$ for the number of hubs,  and ${\bf{Y}}^{(h)}=\{y_1,\ldots, y_h\}\ge\phi{\bf{1}}^{(h)}$ for the rescaled weights of the hubs, i.e.,  $\CV_n[\phi n,\infty)=:n\cdot{\bf{Y}}^{(h)}$. We define two events:
    \begin{align}
    \CA_\mathrm{comp}^\le&:=\big\{\forall \ell\in[\ell_\ast], {\bf{w}}^\sss{(\ell)}\in\mathrm{CT}_\ell(\varepsilon,\bar R): N_{n}({\bf{w}}^\sss{(\ell)}, \eps, \phi n)\ge (1-\psi)\theta({\bf{w}}^\sss{(\ell)}, \eps)n/\ell\big\}\label{eq:comp-ub}\\
    &\hspace{15pt}\cap \{|\CC_n^\sss{(1)}[\underline w,\phi n)|\ge (\theta_q-\psi)n\},\nonumber\\
    \CA_\mathrm{conn}^\le&:=
    \Bigg\{\begin{aligned}\forall \ell\!\in\![\ell_\ast], {\bf{w}}^\sss{(\ell)}\!\in\!\mathrm{CT}_\ell(\varepsilon,\bar R): M_{n}({\bf{w}}^\sss{(\ell)},\eps)&\ge (1-\psi)^3 \cdot(n/\ell)\cdot  \theta({\bf{w}}^\sss{(\ell)}, \eps)\\
    &\hspace{15pt}\cdot \bar P({\bf{w}}^\sss{(\ell)}, {\bf{Y}}^{(h)}) \ind{\bar P({\bf{w}}^\sss{(\ell)}, {\bf{Y}}^{(h)})\ge \psi}
    \end{aligned}\Bigg\}.\label{eq:conn-ub}
    \end{align}
    
    We first show that $\{|\CC_n^\sss{(1)}|\le\rho n\}\subseteq\CA^\le_\mathrm{comp}\cap\CA^\le_\mathrm{hubs, 1\textnormal{-}2}\cap\CA^\le_\mathrm{conn}$. To do so, we bound the size of the complement of the giant on this intersection. Since $\CG_n$ contains a component of size at least $(\theta_q-\psi)n$ by $\CA^\le_\mathrm{comp}$, the complement of the giant contains all components of size at most $\ell_\ast$. Given $\psi>0$, we assume that $\ell_\ast$ and $\bar R$ are sufficiently large, and $\eps$ is sufficiently small, so that \eqref{eq:non-conn-upper} in \cref{lemma:approx-non-connections} applies. On $\CA_\mathrm{comp}^\le\cap\CA_\mathrm{hubs, 1\textnormal{-}2}^\le\cap\CA_\mathrm{conn}^\le$, we obtain
    \[
    \begin{aligned}
    1-\frac{|\CC_n^\sss{(1)}|}n &\ge   \frac{1}{n}\sum_{\ell=1}^{\ell_\ast} \sum_{{\bf{w}}^\sss{(\ell)}\in\mathrm{CT}_\ell(\eps, \bar R)}\hspace{-20pt}\ell M_n({\bf{w}}^\sss{(\ell)}, \bar R, \varepsilon) \\
    &\ge 
    (1-\psi)^3\sum_{\ell=1}^{\ell_\ast}\sum_{{\bf{w}}^\sss{(\ell)}\in\mathrm{CT}_\ell(\eps, \bar R)}\hspace{-20pt}\bar P({\bf{w}}^\sss{(\ell)}, {\bf{Y}}^{(h)})\theta({\bf{w}}^\sss{(\ell)}, \eps)\ind{\bar P({\bf{w}}^\sss{(\ell)}, {\bf{Y}}^{(h)})\ge \psi}\\
    &\ge(1-\psi)^3\Big(\E\big[\bar P\big((W_x)_{x\in T_q}, {\bf{Y}}^{(h)})\big)\big]-\psi\\
    &\hspace{70pt}-\sum_{\ell=1}^{\ell_\ast}\sum_{{\bf{w}}^\sss{(\ell)}\in\mathrm{CT}_\ell(\eps, \bar R)}\hspace{-20pt}\bar P({\bf{w}}^\sss{(\ell)}, {\bf{Y}}^{(h)})\theta({\bf{w}}^\sss{(\ell)}, \eps)\ind{\bar P({\bf{w}}^\sss{(\ell)}, {\bf{Y}}^{(h)})< \psi}\Big) .
    \end{aligned}
    \]
Because of the indicator, we may bound $\bar P({\bf{w}}^\sss{(\ell)}, {\bf{Y}}^{(h)})$ in the sum from above by $\psi$.  We also use that the probabilities $\theta({\bf{w}}^\sss{(\ell)},\eps)$, see~\eqref{eq:branching-type-prob}, sum up to at most one, so that the sum is at most $\psi$. Therefore, on  $\CA_\mathrm{comp}^\le\cap\CA_\mathrm{hubs, 1\textnormal{-}2}^\le\cap\CA_\mathrm{conn}^\le$,
    \[
    1-\frac{|\CC_n^\sss{(1)}|}n \ge (1-\psi)^3\Big(\E\Big[\bar P\big((W_x)_{x\in T_q}, {\bf{Y}}^{(h)})\big)\Big]-2\psi\Big).
    \]
    When $h<\lceil\mathrm{hubs}(r, q)\rceil$, i.e., on $\CA_\mathrm{hubs,2}^\le$, the expectation is larger than $1-r$ by \eqref{eq:hubs-expectation}. When $h=\lceil\mathrm{hubs}(r, q)\rceil$ and ${\bf{Y}}^{(h)}\notin\CY_{r, q}(h)$, i.e., on $\CA_\mathrm{hubs,1}^\le$, the definition of $\CY_{r, q}(h)$ in~\eqref{eq:y-rho-set} implies that the expectation on the right-hand side is larger than $1-r$.  Thus, on 
    $\CA_\mathrm{comp}^\le\cap\CA_\mathrm{hubs, 1\textnormal{-}2}^\le\cap\CA_\mathrm{conn}^\le$
    \[
    1-|\CC_n^\sss{(1)}|/n 
    >(1-\psi)^3(1-r-2\psi) .
    \]
    As we fixed $r\in(\theta_q, \rho)$, the right-hand side is strictly larger than $1-\rho$ if $\psi$ is sufficiently small depending on $r$ and $\rho$. Rearranging yields that $\CA^\le_\mathrm{comp}\cap\CA^\le_\mathrm{hubs, 1\textnormal{-}2}\cap\CA^\le_\mathrm{conn}\subseteq\{|\CC_n^\sss{(1)}|\le\rho n\}$ for such values of $\psi$. This implies that $\{|\CC_n^\sss{(1)}|>\rho n\}\subseteq (\neg \CA^\le_\mathrm{comp})\cup (\neg \CA^\le_\mathrm{conn})$ when working conditionally on $\CA^\le_\mathrm{hubs, 1\textnormal{-}2}$.
    We use that $\Prob(A\mid B\cup C)=\Prob(A\cap (B\cup C))/\Prob(B\cup C)\le \Prob(A)/\Prob(B)$, to obtain 
    \begin{align*}
\Prob\big(|\CC_n^\sss{(1)}|>\rho n  \mid  \CA_\mathrm{hubs, 1\textnormal{-}2}^\le\big) &\le \Prob\big(\neg \CA^\le_\mathrm{comp}  \mid  \CA_\mathrm{hubs, 1\textnormal{-}2}^\le\big)+ \Prob\big((\neg \CA^\le_\mathrm{conn})\cap \CA^\le_\mathrm{comp} \mid  \CA_\mathrm{hubs, 1\textnormal{-}2}^\le\big) \\
     &\le \Prob\big(\neg \CA^\le_\mathrm{comp}\big)/\Prob\big(\CA_\mathrm{hubs, 2}^\le\big)+\Prob\big((\neg \CA^\le_\mathrm{conn})\cap\CA^\le_\mathrm{comp}\mid  \CA_\mathrm{hubs, 1\textnormal{-}2}^\le\big).\end{align*}
   Recall the definition of $\CA^\le_\mathrm{comp}$ from~\eqref{eq:comp-ub}. By Lemmas~\ref{lemma:ell-trunc-n} and~\ref{lemma:giant-trunc-lower} we obtain for $\phi=\phi(\psi)$ sufficiently small that
    \[
    \Prob\big(\neg\CA^\le_\mathrm{comp}\big)=O\big(\big(n\Prob(W_1>n)\big)^{h+1}\big).
    \]
    Since $\Prob\big(\CA_\mathrm{hubs, 2}^\le\big)\to1$, see~\eqref{eq:hubs2-ub} and~\eqref{eq:weight-dist}, it follows that 
    \begin{equation}\label{eq:upper-pr1-last-line}
\Prob\big(|\CC_n^\sss{(1)}|>\rho n  \mid  \CA_\mathrm{hubs, 1\textnormal{-}2}^\le\big)
     \le \Prob\big(\neg \CA^\le_\mathrm{conn}\mid  \CA_\mathrm{hubs, 1\textnormal{-}2}^\le\cap \CA^\le_\mathrm{comp}\big)+O\big(\big(n\Prob(W_1>n)\big)^{h+1}\big). \end{equation}   
    We now condition on the graph $\CG_n[\underline w,\phi n)$ satisfying $\CA^\le_\mathrm{comp}$ and the realization of $\CV_n[\phi n,\infty)$ being equal to $n{\bf{Y}}^{(h)}$ and satisfying $\CA^\le_\mathrm{hubs, 1\textnormal{-}2}$. We abbreviate 
    \begin{equation}\nonumber
        \Prob_{{\bf{Y}}, \CG}\big(\,\cdot\,\big):=\Prob\big(\, \cdot \mid \CG_n[\underline w,\phi n), \CV_n[\phi n,\infty)=n{\bf{Y}}^{(h)}, \CA^\le_\mathrm{comp}\cap \CA^\le_\mathrm{hubs, 1\textnormal{-}2}\big).
    \end{equation}    
    Consequently, 
    \begin{equation}\label{eq:upper-upper-pr2} 
    \Prob\big(\neg \CA^\le_\mathrm{conn}\mid  \CA_\mathrm{hubs, 1\textnormal{-}2}^\le\big) = 
    \E\Big[\ind{\CA^\le_\mathrm{comp}}\Prob_{{\bf{Y}}, \CG}\big(\neg \CA^\le_\mathrm{conn}\big)\Big].
    \end{equation}    
    We recall $\CA_\mathrm{conn}^\le$ from~\eqref{eq:conn-ub}. By a union bound over all component types, it follows that 
    \[
    \begin{aligned}
    \Prob_{{\bf{Y}}, \CG}&\big(\neg \CA^\le_\mathrm{conn}\big) \\
    &\hspace{-20pt}\le 
    \sum_{\ell=1}^{\ell_\ast}\sum_{{\bf{w}}^\sss{(\ell)}\in\mathrm{CT}_\ell(\eps, \bar R)}\hspace{-10pt}
    \Prob_{{\bf{Y}}, \CG}\Big(M_n({\bf{w}}^\sss{(\ell)}, \eps) < (1-\psi)^3(n/\ell) \theta({\bf{w}}^\sss{(\ell)}, \eps) \bar P({\bf{w}}^\sss{(\ell)}, {\bf{Y}}^{(h)}) \ind{\bar P({\bf{w}}^\sss{(\ell)}, {\bf{Y}}^{(h)})\ge \psi}\Big).
    \end{aligned}
    \]
    Since $M_n({\bf{w}}^\sss{(\ell)}, \eps)$ counts the number of components that does not connect by an edge to the hubs, $M_n({\bf{w}}^\sss{(\ell)}, \eps)$ is nonnegative. Therefore, when the indicator inside the probability equals 0, the probability also equals 0. Therefore, we only need to consider the cases in which the indicator equals one, i.e., 
    \[
        \begin{aligned}
    \Prob_{{\bf{Y}}, \CG}&\big(\neg \CA^\le_\mathrm{conn}\big) \\
    &\hspace{-20pt}\le 
    \sum_{\ell=1}^{\ell_\ast}\sum_{{\bf{w}}^\sss{(\ell)}\in\mathrm{CT}_\ell(\eps, \bar R)}\hspace{-10pt}\ind{\bar P({\bf{w}}^\sss{(\ell)}, {\bf{Y}}^{(h)})\ge \psi}
    \Prob_{{\bf{Y}}, \CG}\Big(M_n({\bf{w}}^\sss{(\ell)}, \eps) < (1-\psi)^3(n/\ell) \theta({\bf{w}}^\sss{(\ell)}, \eps) \bar P({\bf{w}}^\sss{(\ell)}, {\bf{Y}}^{(h)})\Big).
    \end{aligned}
    \]
    Components in $\CG_n[\underline w,\phi n)$ of type $({\bf{w}}^\sss{(\ell)}, \eps)$ connect independently by an edge to the hubs. The probability that a component does not connect to the hubs, is at least $(1-\psi)\bar P({\bf{w}}^\sss{(\ell)}, {\bf{Y}}^{(h)})$ by~\eqref{eq:non-conn-component}. Thus, conditionally on $\CG_n[\underline w,\phi n)$ satisfying $\CA^\le_\mathrm{comp}$ defined in~\eqref{eq:comp-ub}, and $\CV_n[\phi n,\infty)=n{\bf{Y}}^{(h)}$,
    \[
    \begin{aligned}
     M_{n}({\bf{w}}^\sss{(\ell)}, \bar R, \eps)&\succcurlyeq \mathrm{Bin}\big(N_{n}({\bf{w}}^\sss{(\ell)}, \bar R, \eps), (1-\psi)\bar P({\bf{w}}^\sss{(\ell)}, {\bf{Y}}^{(h)})\big)\\&\succcurlyeq\mathrm{Bin}\big((1-\psi)\theta({\bf{w}}^\sss{(\ell)}, \eps) n/\ell, (1-\psi)\bar P({\bf{w}}^\sss{(\ell)}, {\bf{Y}}^{(h)})\big).
     \end{aligned}
    \]
    We apply  a Chernoff bound for each type ${\bf{w}}^\sss{(\ell)}$: there exist constants $c>0$  depending on $\psi$ such that 
     \[
        \begin{aligned}
    \Prob_{{\bf{Y}}, \CG}\big(\neg \CA^\le_\mathrm{conn}\big) &\le 
    \sum_{\ell=1}^{\ell_\ast}\sum_{{\bf{w}}^\sss{(\ell)}\in\mathrm{CT}_\ell(\eps, \bar R)}\hspace{-10pt}\ind{\bar P({\bf{w}}^\sss{(\ell)}, {\bf{Y}}^{(h)})\ge \psi}
    \exp\Big(-c\theta({\bf{w}}^\sss{(\ell)}, \eps)(n/\ell)\bar P({\bf{w}}^\sss{(\ell)}, {\bf{Y}}^{(h)})\Big) \\
    &\le 
    \sum_{\ell=1}^{\ell_\ast}\sum_{{\bf{w}}^\sss{(\ell)}\in\mathrm{CT}_\ell(\eps, \bar R)}\hspace{-10pt}
    \exp\Big(-c\psi\theta({\bf{w}}^\sss{(\ell)}, \eps)(n/\ell)\Big).    
    \end{aligned}
    \]
    Since both $\ell_\ast$ and the number of considered component types for each $\ell$ are finite, the right-hand side decays exponentially fast in $n$. Thus, also~\eqref{eq:upper-upper-pr2}  decays exponentially fast, and~\eqref{eq:upper-pr1-last-line} is of order $O\big(\big(n\Prob(W_1>n)\big)^{h+1}\big)$. We substitute that bound into~\eqref{eq:ub-below-hubs}, which finishes the proof since the second term in~\eqref{eq:ub-below-hubs} is of the same order.
\end{proof}

Next, we state and prove a proposition that shows that the presence of $\lceil\mathrm{hubs}(r, q)\rceil$ hubs with  weights in $\CY_{r,q}(\lceil\mathrm{hubs}(\rho,q)\rceil)$ lead to a large giant for any $r>\rho$.

\begin{proposition}[Presence of hubs implies a large giant]\label{prop:upper-bound}
Consider an inhomogeneous scale-free random graph as in \cref{def:irg}. Fix a constant $\rho\in(\theta_q, 1)$. For any $r\in(\rho, 1)$ there exists a constant $c>0$ such that, with $h=\lceil\mathrm{hubs}(r,q)\rceil$, for any ${\bf{y}}^{(h)}\in\CY_{r,q}(h)$ and $n\ge 1$,
\begin{equation}\label{eq:prop-lower}
\Prob\Big(|\CC_n^\sss{(1)}|> \rho n \, \big|\, \CV_n[\phi n, \infty)=n{\bf{y}}^{(h)}\Big) \ge 1-\exp(-cn).
\end{equation}
\end{proposition}
\begin{proof}[Proof]
The proof uses a similar construction as  \cref{prop:lower-bound}.  We reveal the graph $\CG_n$ in two stages:  we first consider the graph $\CG_n[\underline w,\bar R)$ for a large constant $\bar R>\underline w$ (this is opposed to the proof of~\cref{prop:upper-bound}, where we first revealed the graph up to $\phi n$). \cref{lemma:ell-trunc-n} yields concentration for the number of size-$\ell$ components of finitely many types in this graph \cref{lemma:ell-trunc-n}. The giant has size approximately $\theta_q n$ by \cref{lemma:giant-trunc-lower}. Then we `add' the vertices of weight at least $\phi n$ to the graph, which we call hubs. Each hub connects with probability tending to one exponentially fast to one of the vertices in the giant of $\CG[\underline w, \bar R)$. Since we assume that there are exactly $\lceil\mathrm{hubs}(r,q)\rceil$ many hubs, whose set of their weights is contained in $n\cdot\CY_{r, q}$, sufficiently many components connect to one of the hubs. Consequently, the size of the giant becomes larger than $\rho n$ with sufficiently high probability. We will not consider the effect of vertices with weight in $[\bar R, \phi n)$: edges incident to those vertices can only increase the size of the largest component, but this effect is negligible.

We now formalize this. 
Fix $r\in(\theta_q, \rho)$, let
$\psi=\psi(r), \varepsilon=\varepsilon(\psi)$  be sufficiently small constants, and let $\bar R=\bar R(\psi), \ell_\ast=\ell_\ast(\psi)$ be two large constants.
For each $\ell\le \ell_\ast$,
    we write $M_{n}({\bf{w}}^\sss{(\ell)},\eps)$ for the number of components of type ${\bf{w}}^\sss{(\ell)}$ in the induced subgraph $\CG_n[\underline w,\bar R)$ that are \emph{not} connected by an edge to the hubs in $\CG_n$. Let $h:=|\CV_n[\phi n,\infty)|$ denote the number of hubs, which is equal to $\lceil\mathrm{hubs}(r, q)\rceil$ by assumption. By the conditioning in~\eqref{eq:prop-lower} all rescaled weights in ${\bf{y}}^{(h)}=\{y_1,\ldots, y_h\}$ are at least $\phi$.  Define 
    \begin{align}
    \CA^\ge_\mathrm{comp}&:=\big\{\forall \ell\in[\ell_\ast], {\bf{w}}^\sss{(\ell)}\in\mathrm{CT}_\ell(\varepsilon,\bar R): \big|\ell N_{n}({\bf{w}}^\sss{(\ell)}, \eps, \bar R)\big/(n\theta({\bf{w}}^\sss{(\ell)}, \eps)) - 1\big|\ge \psi\big\},\label{eq:lower-comp}\\
    &\hspace{15pt}\cap\{|\CC_n^\sss{(1)}[\underline w,\phi n)|\ge (\theta_q-\psi) n\}\nonumber\\
    \CA^\ge_\mathrm{hubs}&:=\big\{\forall v\in\CV_n[\phi n,\infty): v\sim \CC_n^\sss{(1)}[\underline w,\bar R)\big\},\label{eq:lower-hubs2} \\
    \CA^\ge_\mathrm{conn}&:=\Bigg\{\begin{aligned}\forall \ell\!\in\![\ell_\ast], {\bf{w}}^\sss{(\ell)}\!\in\!\mathrm{CT}_\ell(\varepsilon,\bar R):  M_{n}({\bf{w}}^\sss{(\ell)},\eps)&\le (1+\psi)^2\cdot\big(\psi+ \bar P({\bf w}^\sss{(\ell)}, {\bf{y}}^{(h)})\big)\\&\hspace{15pt}\cdot N_{n}({\bf{w}}^\sss{(\ell)}, \eps, \bar R)\end{aligned}\Bigg\}.\label{eq:lower-conn}
    \end{align}
    We bound the size of the giant on the intersection of the three events from below. Since all vertices in $\CV_n[\phi n,\infty)$ connect by an edge to the largest component in $\CG_n[\underline w,\bar R)$ on $\CA^\ge_\mathrm{hubs}$, the size of the largest component $|\CC_n^\sss{(1)}|$ increases from $|\CC_n^\sss{(1)}[1,\bar R)|$ by at least the  total number of vertices in a component of size at most $\ell_\ast$ in $\CG_n[\underline w,\bar R)$ that is connected by an edge to one of the hubs in $\CG_n$. Thus, using first the definition of $ \CA^\ge_\mathrm{conn}$, and then the definition of  $\CA^\ge_\mathrm{comp}$,
    \[
    \begin{aligned}
    |\CC_n^\sss{(1)}|/n&\ge \theta_q-\psi  + \frac{1}{n}\sum_{\ell=1}^{\ell_\ast}\sum_{{\bf{w}}^\sss{(\ell)}\in\mathrm{CT}_\ell(\eps, \bar R)}\ell\big(N_{n}({\bf{w}}^\sss{(\ell)}, \eps, \bar R) -  M_{n}({\bf{w}}^\sss{(\ell)},\eps)\big)\\
    &\ge 
    \theta_q-\psi + 
    \frac{1}{n}\sum_{\ell=1}^{\ell_\ast}\sum_{{\bf{w}}^\sss{(\ell)}\in\mathrm{CT}_\ell(\eps, \bar R)}\ell N_{n}({\bf{w}}^\sss{(\ell)}, \eps, \bar R)\Big(1-(1+\psi)^2\big(\psi+ \bar P({\bf{w}}^\sss{(\ell)}, {\bf{y}}^{(h)})\big)\Big) \\
    &\ge 
    \theta_q-\psi+(1-\psi)\sum_{\ell=1}^{\ell_\ast}\sum_{{\bf{w}}^\sss{(\ell)}\in\mathrm{CT}_\ell(\eps, \bar R)}\theta({\bf{w}}^\sss{(\ell)}, \eps)\\
    &\hspace{15pt}-
    (1+\psi)^3
    \psi\sum_{\ell=1}^{\ell_\ast}\sum_{{\bf{w}}^\sss{(\ell)}\in\mathrm{CT}_\ell(\eps, \bar R)}\hspace{-10pt}
    \theta({\bf{w}}^\sss{(\ell)}, \eps)-
    (1+\psi)^3
    \sum_{\ell=1}^{\ell_\ast}\sum_{{\bf{w}}^\sss{(\ell)}\in\mathrm{CT}_\ell(\eps, \bar R)}\hspace{-10pt}
    \theta({\bf{w}}^\sss{(\ell)}, \eps)\bar P({\bf{w}}^\sss{(\ell)}, {\bf{y}}^{(h)}).
    \end{aligned}
    \]
    We assume that $\eps$, $\bar R$, and $\ell_\ast$ satisfy the conditions of \cref{lemma:approx-non-connections}.
    We bound the first double sum via~\eqref{eq:non-conn-exp-lower}.  The second double sum over the probabilities $\theta({\bf{w}}^\sss{(\ell)},\eps)$ is at most 1 by~\eqref{eq:branching-type-prob}. We invoke~\eqref{eq:non-conn-upper} for the third double sum. 
    Therefore, on $\CA^\ge_\mathrm{comp}\cap\CA^\ge_\mathrm{hubs}\cap\CA_\mathrm{conn}^\ge$, 
    \[
        \begin{aligned}
    |\CC_n^\sss{(1)}|/n     
    &\ge \theta_q-\psi+(1-\psi)(1-\theta_q-\psi) -\psi(1+\psi)^3 - (1+\psi)^3\Big(\E\Big[\bar P\big((W_x)_{x\in T_q}, {\bf{y}}^{(h)}\big)\Big] + \psi\Big).
    \end{aligned}
    \]
     By assumption, we have ${\bf{y}}^{(h)}\in\CY_{r,q}$.
    Thus, the expectation is at most $1-r$ by definition of $\CY_{r, q}$ in~\eqref{eq:y-rho-set}. Hence,  for some constant $C>0$
    \[
    \begin{aligned}
    |\CC_n^\sss{(1)}|/n     
    &\ge 
    \theta_q-\psi+(1-\psi)(1-\theta_q-\psi)-\psi(1+\psi)^3  - (1+\psi)^3(1-r+\psi)\ge r-C\psi.
    \end{aligned}
    \]
    Since $r>\rho$, we can make the  right-hand side is strictly larger than $\rho$ when $\psi$ is fixed  sufficiently small. 
    So, for $\psi$ sufficiently small, and $\eps$, $\bar R$, and $\ell_\ast$ satisfying the conditions of \cref{lemma:approx-non-connections},
    \[
    \Prob\big(|\CC_n^\sss{(1)}|>\rho n \, \big|\, \CV_n[\phi n, \infty)=n{\bf{y}}^\sss{(\ell)}\big) \ge 
    \Prob\big(\CA^\ge_\mathrm{comp}\cap\CA^\ge_\mathrm{hubs}\cap\CA^\ge_\mathrm{conn}\, \big|\, \CV_n[\phi n, \infty)=n{\bf{y}}^\sss{(\ell)} \big).
    \]
    In the remainder, we bound the probability on the right-hand side from below. We condition on the graph $\CG_n[\underline w,\bar R)$ satisfying $\CA^\ge_\mathrm{comp}$ and the realization of $\CV_n[\phi n,\infty)=n{\bf{y}}^{(h)}$ satisfying $\CA^\ge_\mathrm{comp}$. We abbreviate 
    \begin{equation}\label{eq:cond-meas-lb}
        \Prob_{{\bf{y}}, \CG}\big(\,\cdot\,\big):=\Prob\big(\, \cdot \mid \CG_n[\underline w,\bar R), \CV_n[\phi n,\infty)=n{\bf{y}}^{(h)}, \CA^\le_\mathrm{comp}\big).
    \end{equation}    
    Consequently, 
    \begin{equation}\label{eq:lower-pr-1}
    \Prob\big(|\CC_n^\sss{(1)}|>\rho n\, \big|\, \CV_n[\phi n, \infty)=n{\bf{y}}^\sss{(\ell)}\big) 
    \ge 
    \E\Big[\ind{\CA^\ge_\mathrm{comp}}\Prob_{{\bf{y}}, \CG}\big(\CA^\ge_\mathrm{hubs}\cap\CA^\ge_\mathrm{conn}\big)\,\big|\,\CV_n[\phi n, \infty)=n{\bf{y}}^\sss{(\ell)}\Big].
    \end{equation}
    We first show that the conditional probability tends to $1$ for any realization of $\CG_n[\underline w,\bar R)$ satisfying $\CA^\ge_\mathrm{comp}$. 
    Since we condition on the weights of all vertices with weight in $[\underline w, \bar R)\cup[\phi n,\infty)$ in~\eqref{eq:cond-meas-lb}, edges between the hubs and components in $\CG[\underline w, \bar R)$ are present independently. Thus, by definition of $\CA^\ge_\mathrm{hubs}$ and $\CA^\ge_\mathrm{conn}$ in~\eqref{eq:lower-hubs2} and~\eqref{eq:lower-conn},
    \begin{equation}\label{eq:lower-pr-hub2}
    \Prob_{{\bf{y}}, \CG}\big(\CA^\ge_\mathrm{hubs}\cap\CA^\ge_\mathrm{conn}\big)
    =
    \Prob_{{\bf{y}}, \CG}\big(\CA^\ge_\mathrm{conn}\big)
    \cdot\Prob_{{\bf{y}}, \CG}\big( \CA^\ge_\mathrm{hubs}\big).
    \end{equation}
    We next show that the first probability on the right-hand side tends to one.
    By  definition of $\CA^\ge_\mathrm{conn}$ in~\eqref{eq:lower-conn}, we have to bound from above the number of components of type ${\bf{w}}^\sss{(\ell)}$ that do \emph{not} connect by an edge to one of the hubs.  We apply a union bound over all component types: 
    \begin{equation}\label{eq:conn-lower-pr-1}
    \begin{aligned}
    \Prob_{{\bf{y}}, \CG}&\big(\neg\CA^\ge_\mathrm{conn}\big)\\&\hspace{-20pt}\le 
    \sum_{\ell=1}^{\ell_\ast}\sum_{{\bf{w}}^\sss{(\ell)}\in\mathrm{CT}_\ell(\eps, \bar R)}\hspace{-15pt}
    \Prob_{{\bf{y}}, \CG}\Big(M_{n}({\bf{w}}^\sss{(\ell)}, \bar R, \eps)>(1+\psi)^2\big(\psi+ \bar P({\bf w}^\sss{(\ell)}, {\bf{y}}^{(h)})\big)N_n({\bf w}^\sss{(\ell)}, \eps, \bar R)\Big).
    \end{aligned}
    \end{equation}
    Components in $\CG_n[\underline w,\bar R)$ of type $({\bf{w}}^\sss{(\ell)}, \eps)$ connect independently by an edge to the hubs, and the probability that a component does not connect to the hubs, is at least $(1-\psi)\bar P({\bf{w}}^\sss{(\ell)}, {\bf{y}}^{(h)})$ by~\eqref{eq:non-conn-component}. Thus, under the conditional probability measure $\Prob_{{\bf{y}}, \CG}$,
    \[
    \begin{aligned}
     M_{n}({\bf{w}}^\sss{(\ell)}, \bar R, \eps)&\preccurlyeq \mathrm{Bin}\big(N_{n}({\bf{w}}^\sss{(\ell)}, \bar R, \eps), (1+\psi)\bar P({\bf{w}}^\sss{(\ell)}, {\bf{y}}^{(h)})\big)\\
     &\preccurlyeq  \mathrm{Bin}\Big(N_{n}({\bf{w}}^\sss{(\ell)}, \bar R, \eps), (1+\psi)\big(\psi+ \bar P({\bf{w}}^\sss{(\ell)}, {\bf{y}}^{(h)})\big)\Big).
     \end{aligned}
    \]
    We apply the stochastic domination to all terms in~\eqref{eq:conn-lower-pr-1} and apply a Chernoff bound for each component type. This yields for some $c=c(\psi)>0$
    \[
    \begin{aligned}
    \Prob_{{\bf{y}}, \CG}\big(\neg\CA^\ge_\mathrm{conn}\big)&\le \sum_{\ell=1}^{\ell_\ast}\sum_{{\bf{w}}^\sss{(\ell)}\in\mathrm{CT}_\ell(\eps, \bar R)}
    \exp\Big(-c(\bar P({\bf{w}}^\sss{(\ell)}, {\bf{y}}^{(h)}) +\psi) N_{n}({\bf{w}}^\sss{(\ell)}, \bar R, \eps)\Big)
    \\
    &\le \sum_{\ell=1}^{\ell_\ast}\sum_{{\bf{w}}^\sss{(\ell)}\in\mathrm{CT}_\ell(\eps, \bar R)}\exp\Big(-c\psi N_{n}({\bf{w}}^\sss{(\ell)}, \bar R, \eps)\Big) .
    \end{aligned}
    \]
    The probability measure $\Prob_{{\bf{y}}, \CG}$ defined in~\eqref{eq:cond-meas-lb} is conditional on the event $\CA^\le_\mathrm{comp}$ defined in~\eqref{eq:lower-comp}. On this event, each $N_{n}({\bf{w}}^\sss{(\ell)}, \bar R, \eps)$ increases linearly in $n$. Since the number of component types is finite, $\Prob_{{\bf{y}}, \CG}\big(\neg\CA^\ge_\mathrm{conn}\big)$ decays exponentially in $n$. We substitute this bound into~\eqref{eq:lower-pr-hub2} and obtain for some other $c>0$,
    \begin{equation}\label{eq:lower-pr-hub3}
    \Prob_{{\bf{y}}, \CG}\big(\CA^\ge_\mathrm{hubs}\cap\CA^\ge_\mathrm{conn}\big)
    \ge
    \big(1-\exp(-cn)\big)
    \cdot\Prob_{{\bf{y}}, \CG}\big( \CA^\ge_\mathrm{hubs}\big).
    \end{equation}
    We bound the second factor on the right-hand side. Recall $\CA^\ge_\mathrm{hubs}$ from~\eqref{eq:lower-hubs2}.    
    On the event $\CA^\ge_\mathrm{comp}$ on which we conditioned, $|\CC_n[\underline w, \bar R)|\ge (\theta_q-\psi)n$, and each vertex of weight at least $\phi n$ connects with probability at least $q(\phi \underline w\wedge 1)=\Theta(1)$ to each vertex in  $\CC_n^\sss{(1)}[\underline w,\bar R)$. Thus, the probability that a single vertex of weight at least $\phi n$ does not connect by an edge to the giant, decays exponentially in $n$. By assumption, there are exactly $\lceil\mathrm{hubs}(r, q)\rceil$ many vertices in $\CV_n[\phi n,\infty)$. By a union bound over these constantly many hubs, also the second factor in~\eqref{eq:lower-pr-hub3} tends to one on the event $\CA^\ge_\mathrm{comp}$ exponentially fast. Substituting this limit into~\eqref{eq:lower-pr-1} yields for some other $c>0$,
    \begin{equation}\label{eq:lower-pr-hub4}
    \begin{aligned}
    \Prob\big(|\CC_n^\sss{(1)}|>\rho n\, \big|\, \CV_n[\phi n, \infty)=n{\bf{y}}^\sss{(\ell)}\big) 
    &\ge \big(1-\exp(-cn)\big)\Prob\big(\CA^\ge_\mathrm{comp}\big).
    \end{aligned}
    \end{equation}
    By definition of $\CA^\ge_\mathrm{comp}$ in~\eqref{eq:lower-comp}, 
    \[
    \begin{aligned}
    \Prob(\neg\CA^\ge_\mathrm{comp}) &\le \Prob\big(\exists \ell\in[\ell_\ast], {\bf{w}}^\sss{(\ell)}\in\mathrm{CT}_\ell(\varepsilon,\bar R): N_{n}({\bf{w}}^\sss{(\ell)}, \eps, \bar R)< (1-\psi)\theta({\bf{w}}^\sss{(\ell)}, \eps)n/\ell\big)\\
    &\hspace{15pt}+\Prob\big(|\CC_n^\sss{(1)}[\underline w,\bar R)|< (1-\psi)\theta_q n\big).
    \end{aligned}
    \]     
    By \cref{lemma:ell-trunc-n} and \cref{lemma:giant-trunc-lower}, the two terms on the right-hand side decay exponentially in $n$. Combined with~\eqref{eq:lower-pr-hub4}, this proves~\eqref{eq:prop-lower}.
\end{proof}
Via similar proof strategies as Propositions~\ref{prop:lower-bound} and~\ref{prop:upper-bound}, we prove the following lemma in the appendix on page~\pageref{proof:concentration-components}. 
\begin{lemma}\label{lemma:size-ell-hubs}
Consider an inhomogeneous scale-free random graph as in \cref{def:irg}.
    Moreover, for any constants $C, \psi>0$ and $\ell, h\in\N$, there exists $\phi_1>0$ such that for any $\phi\in(0,\phi_1)$, as $n\to\infty$,
\[
\sup_{{\bf{y}}^{(h)}\ge\phi{\bf{1}}^{(h)}}\Prob\Big(\big|N_{n,\ell}/n -\tfrac1\ell\E\big[\ind{|T_q|=\ell}\bar P\big((W_x)_{x\in T_q}, {\bf{y}}^{(h)}\big)\big]\big|>\psi\, \big|\, \CV_n[\phi n, \infty)=n{\bf{y}}^{(h)}\Big)\le o(n^{-C}).
\]
\end{lemma}

\subsection{Asymptotics for the probability of having hubs}
In this section, we analyze the probabilities of the events on which we conditioned in Propositions~\ref{prop:lower-bound} and~\ref{prop:upper-bound}. Let $C_{\rho, q}(h)$ be the constant from~\eqref{eq:constant-value}.
\begin{lemma}[Leading constant]\label{lemma:leading-term}
Consider an inhomogeneous scale-free random graph as in \cref{def:irg}  with $\sigma<2\alpha - 1$. Let $\rho\in(\theta_q,1)$ and set $h=\lceil\mathrm{hubs}(\rho,q)\rceil$. There exists $\phi_0>0$ such that for any $\phi\in(0,\phi_0)$, as $n\to\infty$,
    \begin{equation}\label{eq:lemma-leading-term-1}
    \frac{\Prob\big(|\CV_n[\phi n, \infty)|=h, \CV_n[\phi n,\infty)\in n\cdot\CY_{\rho, q}(h)\big)}{\big(n\Prob(W_1> n)\big)^{h}}\longrightarrow C_{\rho, q}(h)\in(0,\infty).
    \end{equation}
\end{lemma}

\begin{proof}
    We start with analyzing the numerator of~\eqref{eq:lemma-leading-term-1}. Abbreviate $h=\lceil \mathrm{hubs}(\rho, q)\rceil$. By \cref{def:irg}, all weights are iid and follow distribution $F_W(w)=1-L(w)w^{-\alpha}$ for a parameter $\alpha>1$ and slowly varying function $L(w)$. Thus, 
    \[
    \begin{aligned}
    \Prob\big(|\CV_n[\phi n, \infty)|=h\big)&=\binom{n}{h}\big(L(\phi n)(\phi n)^{-\alpha}\big)^h\big(1-L(\phi n)(\phi n)^{-\alpha}\big)^{n-h}\\&\sim\frac{1}{h!}\big(nL(\phi n)(\phi n)^{-\alpha}\big)^h.
    \end{aligned}
    \]
    Since $L$ is slowly varying, $L(\phi n)n^{-\alpha}\sim L(n)n^{-\alpha}\sim\Prob\big(W_1>n\big)$. So, 
    \begin{equation}
    \begin{aligned}\label{eq:proof-limiting-constant-1}
    &\frac{\Prob\big(|\CV_n[\phi n, \infty)|=h, \CV_n[\phi n,\infty)\in n\cdot\CY_{\rho, q}(h)\big)}{\big(n\Prob(W_1> n)\big)^{h}}\nonumber\\&\hspace{80pt} \sim\frac{1}{\phi^{\alpha h}h!}\Prob\Big(\frac{1}{n}\CV_n[\phi n,\infty)\in \CY_{\rho, q}(h)\,\big|\, |\CV_n[\phi n, \infty)|=h\Big).
    \end{aligned}
    \end{equation}
    Let $(W_1/n,\ldots, W_h/n)$ denote the weights of the vertices with weight at least $\phi n$. Since the weights are iid and regularly varying with index $\alpha$, conditionally on $W_i\ge\phi n$ for all $i\in[h]$, 
    \[
    (W_1/n,\ldots, W_h/n)\overset{d}\longrightarrow ( Y_1,\ldots,  Y_h),
    \]
    where $(Y_1,\ldots, Y_h)$ are independent copies of $Y$ following distribution $\Prob\big(Y\ge y\big)=(y/\phi)^{-\alpha}$ for $y\ge \phi$. Writing out the probability in~\eqref{eq:proof-limiting-constant-1} as an integral, the factor $\phi^{-\alpha h}$ cancels, i.e., 
    \begin{equation}\label{eq:proof-limiting-constant-middle}
     \begin{aligned}
    &\frac{\Prob\big(|\CV_n[\phi n, \infty)|=h, \CV_n[\phi n,\infty)\in n\cdot\CY_{\rho, q}(h)\big)}{\big(n\Prob(W_1> n)\big)^{h}}\\&\hspace{40pt} \sim\frac{1}{h!}\int_{y_1=\phi}^\infty\cdots\int_{y_h=\phi}^\infty\frac{\alpha^h}{(y_1\cdot\ldots\cdot y_h)^{\alpha+1}}\mathds{1}\{( y_1,\ldots, y_h)\in\CY_{\rho, q}\}\rd y_1\cdot\ldots\cdot\rd y_h.
    \end{aligned}
    \end{equation}
    The only difference between the integrals here compared to the ones in the definition of $C_{\rho, q}$ in~\eqref{eq:constant-value}, is that the integrals in~\eqref{eq:constant-value} start at the value 0, rather than $\phi$. In the remainder of the proof we argue that the indicator above always equals 0 when there exists an index $i$ such that $ y_i<\phi$ for a small constant $\phi>0$. Let $(y_1,\ldots, y_h)$ be any vector such that $y_i<\phi$ for some $i$. By definition of $\CY_{\rho, q}$ in~\eqref{eq:y-rho-set}, we have to show that 
    \begin{equation}\label{eq:proof-limiting-constant-before-mid}
    \E\bigg[\prod_{x\in T_q, j\in[h]}\big(1-q(W_x^\sigma y_j\wedge 1)\big)\bigg] > 1-\rho.
    \end{equation}
    Bounding all factors with $j\neq i$ from below by $(1-q)$, and $y_i<\phi$, we obtain 
    \begin{equation}\label{eq:proof-limiting-constant-2}
    \E\bigg[\prod_{x\in T_q, j\in[h]}\big(1-q(W_x^\sigma y_j\wedge 1)\big)\bigg] \ge \E\bigg[(1-q)^{|T_q|(h-1)}\prod_{x\in T_q}\big(1-q(W_x^\sigma \phi\wedge 1)\big)\bigg].
    \end{equation}
    The expectation is non-increasing in $\phi$, so we evaluate the limit as $\phi\downarrow 0$.
    The argument of the expectation is a continuous and monotone function of $\phi$. Therefore, by the Monotone Convergence Theorem, 
    \[
    \lim_{\phi\downarrow 0}
    \E\bigg[(1-q)^{|T_q|(h-1)}\prod_{x\in T_q}\big(1-q(W_x^\sigma \phi\wedge 1)\big)\bigg]
    =\E\Big[(1-q)^{|T_q|(h-1)}\Big].
    \]
    By definition of $h=\lceil\mathrm{hubs}(\rho, q)\rceil$ in~\eqref{eq:hubs}, $\E\big[(1-q)^{|T_q|\mathrm{hubs}(\rho, q)}\big]=1-\rho$ when $q<1$. If $q=1$, then $\mathrm{hubs}(\rho, q)=1$. Therefore, regardless of $q$ (using for the case $q=1$ that $\rho\in(0,1)$),
    \[
    \lim_{\phi\downarrow 0}
    \E\bigg[(1-q)^{|T_q|(h-1)}\prod_{x\in T_q}\big(1-q(W_x^\sigma \phi\wedge 1)\big)\bigg]
    > 
    \E\bigg[(1-q)^{|T_q|\mathrm{hubs}(\rho, q)}\bigg]
    =1-\rho.
    \]
    Since the expectation on the right-hand side in~\eqref{eq:proof-limiting-constant-2} is a non-increasing function in $\phi$, there exists $\phi>0$ such that~\eqref{eq:proof-limiting-constant-before-mid} holds if there exists $y_i<\phi$.
    This proves that $(y_1,\ldots, y_h)\notin\CY_{\rho, q}$ if there exists $y_i<\phi$. Using this in \eqref{eq:proof-limiting-constant-middle}, we obtain 
\[
\begin{aligned}
&\frac{\Prob\big(|\CV_n[\phi n, \infty)|=h, \CV_n[\phi n,\infty)\in n\cdot\CY_{\rho, q}(h)\big)}{\big(n\Prob(W_1> n)\big)^{h}}\\&\hspace{40pt} \sim\frac{1}{h!}\int_{y_1=0}^\infty\cdots\int_{y_h=0}^\infty\frac{\alpha^h}{(y_1\cdot\ldots\cdot y_h)^{\alpha+1}}\mathds{1}\{( y_1,\ldots, y_h)\in\CY_{\rho, q}\}\rd y_1\cdot\ldots\cdot\rd y_h.\end{aligned}
\]
The integral on the right-hand side equals $C_{\rho, q}$ by its definition in~\eqref{eq:constant-value}, proving~\eqref{eq:lemma-leading-term-1}.
\end{proof}
As the bounds in Propositions~\ref{prop:lower-bound}--\ref{prop:upper-bound} are not stated for $\rho\in(\theta_q, 1)$, but for $r$ close to $\rho$, we analyze the scaling of $C_{\rho, q}(h)$. 
\begin{lemma}[Limit of the constant]\label{lemma:constant-limit}
Consider an inhomogeneous scale-free random graph as in \cref{def:irg}. Let $\rho\in(\theta_q,1)$, and set $h=\lceil\mathrm{hubs}{(\rho, q)}\rceil$. Then $\rho\mapsto C_{\rho, q}(h)$ is strictly decreasing. Moreover,
    \begin{align}
        \label{eq:lemma-leading-term-2}
        C_{\rho, q}(h)-\lim_{r\downarrow \rho}C_{r, q}(h)
        &= 
        \begin{dcases}
            0,&\text{if }\mathrm{hubs}(\rho, q)\notin\N \text{ or }q=1,\\
            \underline w^{\alpha\sigma}/h!,&\text{if }\mathrm{hubs}(\rho, q)\in\N \text{ and }q<1,
        \end{dcases}
        \intertext{and }
        \label{eq:lemma-leading-term-3}
        \lim_{r\uparrow \rho}C_{r, q}(h)&=C_{\rho, q}(h).
    \end{align}
\end{lemma}
\begin{proof}
We first compute the limit of $(h!/\alpha^h)\big(C_{\rho, q}(h)-\lim_{r\downarrow \rho}C_{r, q}(h)\big)$.
    By definition of $C_{\rho, q}$ in~\eqref{eq:constant-value}, this corresponds to evaluating 
\begin{equation}\label{eq:limiting-constant-pr-1}
\lim_{r\downarrow\rho}\int_{y_1=0}^\infty\cdots\int_{y_h=0}^\infty\mathds{1}\{(y_1,\ldots, y_h)\in\CY_{\rho, q}\setminus\CY_{r, q}\}\cdot (y_1\cdot\ldots\cdot y_h)^{-(\alpha+1)}\rd y_1\cdot\ldots\cdot\rd y_h,
\end{equation}
Assume $r\in(\rho, 1)$. If ${\bf{y}}\in\CY_{\rho, q}\setminus\CY_{r, q}$, then by definition of $\CY_{\rho, q}$ in~\eqref{eq:y-rho-set}, 
\[
1-r<\E\bigg[\prod_{x\in T_q, i\in[h]}\big(1-q(W_x^\sigma y_i\wedge 1)\big)\bigg]\le 1-\rho.
\]
Thus the indicator function in~\eqref{eq:limiting-constant-pr-1} is monotone in $r$.
By the Monotone Convergence Theorem, the integral in~\eqref{eq:limiting-constant-pr-1} corresponds to 
\begin{equation}
\int_{y_1=0}^\infty\cdots\int_{y_h=0}^\infty\mathds{1}\bigg\{{\bf{y}}: \E\bigg[\prod_{x\in T_q, i\in[h]}\big(1-q(W_x^\sigma y_i\wedge 1)\big)\bigg]=1-\rho\bigg\}\cdot (y_1\cdot\ldots\cdot y_h)^{-(\alpha+1)}\rd y_1\cdot\ldots\cdot\rd y_h.\label{eq:limiting-constant-pr-2}
\end{equation}
Given $y_1,\ldots, y_{h-1}$, we compute the integral over $y_h$. To do so, we compare the function 
\[
f_{h-1}(y_h):=\E\bigg[\prod_{x\in T_q, i\in[h]}\big(1-q(W_x^\sigma y_i\wedge 1)\big)\bigg]
\]
to $1-\rho$.
The function $f_{h-1}(y_h)$ is non-increasing for $y_h\ge 0$. We first rule out that $f_h(0)=1-\rho$:
regardless of $y_1,\ldots, y_{h-1}$ 
\begin{equation}\label{eq:limit-constant-0}
f_{h-1}(0)\ge \E\big[(1-q)^{|T_q|(\lceil\mathrm{hubs}(\rho, q)\rceil -1)}\big]>\E\big[(1-q)^{|T_q|\mathrm{hubs}(\rho, q) }\big]=1-\rho,
\end{equation}
where the equality follows from the definition of $\mathrm{hubs}$ in~\eqref{eq:hubs}. 
Next, we show that $f_{h-1}(y_h)$ is \emph{strictly} decreasing for $y_h\in(0, \underline w^{-\sigma})$. Indeed, fix $y_h, y_{h}'\in(0,\underline w^{-\sigma})$ with $y_h<y_h'$, and set $y_i':=y_i$ for $i<h$. The event that the root of the branching process has weight $W_\varnothing\in[\underline w, {y_h'}^{-1/\sigma})$ has strictly positive probability (note that the interval is indeed non-empty). Therefore, one can show that  
\[
f_{h-1}(y_h)-f_{h-1}(y_h')=\E\bigg[\prod_{x\in T_q, i\in[h]}\big(1-q(W_x^\sigma y_i\wedge 1)\big)-\prod_{x\in T_q, i\in[h]}\big(1-q(W_x^\sigma y'_i\wedge 1)\big)\bigg]>0.
\]
By strict monotonicity, the indicator in the integral in~\eqref{eq:limiting-constant-pr-2} over $y_h$ evaluates to one in at most a single point $y_h^\ast$, which is strictly positive by~\eqref{eq:limit-constant-0}. Hence, the integral vanishes for $y<\underline w^{-\sigma}$. Thus, the integral in~\eqref{eq:limiting-constant-pr-2} is equal to
\begin{equation}\label{eq:limiting-constant-pr-3}
\begin{aligned}
\int_{y_1=0}^\infty\cdots\int_{y_h=\underline w^{-\sigma}}^\infty\mathds{1}\big\{{\bf{y}}: f_{h-1}(y_h)=1-\rho\big\}\cdot (y_1\cdot\ldots\cdot y_h)^{-(\alpha+1)}\rd y_1\cdot\ldots\cdot\rd y_h.
\end{aligned}
\end{equation}
For $y_h\ge \underline w^{-\sigma}$, the function $f_{h-1}(y_h)$ is constant: each vertex $x\in T_q$ has weight at least $\underline w$ by \cref{def:branching}, so for all $y_h\ge \underline w^{-\sigma}$, 
\[
f_{h-1}(y_h)=\E\bigg[(1-q)^{|T_q|}\prod_{i\in[h-1], x\in T_q}\big(1-q(W_x^\sigma y_i\wedge 1)\big)\bigg],
\]
which only depends on $y_1,\ldots, y_{h-1}$.
Thus, also the indicator in~\eqref{eq:limiting-constant-pr-3} is a constant. Integrating over $y_h\ge \underline w^{-\sigma}$ yields that we have to evaluate 
\begin{equation}\label{eq:limiting-constant-pr-4}
\begin{aligned}
\frac{\underline w^{\sigma\alpha}}{\alpha}\int_{y_1=0}^\infty\cdots\int_{y_{h-1}=0}^\infty\mathds{1}\bigg\{{\bf{y}}: \E\bigg[(1-q)^{|T_q|}&\prod_{i\in[h-1], x\in T_q}\big(1-q(W_x^\sigma y_i\wedge 1)\big)\bigg]=1-\rho\bigg\}\\&\hspace{20pt}\cdot (y_1\cdot\ldots\cdot y_{h-1})^{-(\alpha+1)}\rd y_1\cdot\ldots\cdot\rd y_{h-1}.
\end{aligned}
\end{equation}
Now we iterate the reasoning from~\eqref{eq:limiting-constant-pr-2} until~\eqref{eq:limiting-constant-pr-4}
to compute the integral over the variables $y_{h-1},\ldots, y_1$. Thus, the limit in~\eqref{eq:limiting-constant-pr-1}, corresponding to $(h!/\alpha^h)\big(C_{\rho, q}(h)-\lim_{r\downarrow \rho}C_{r, q}(h)\big)$,  is equal to 
\[
\Big(\frac{\underline w^{\sigma\alpha}}{\alpha}\Big)^h\mathds{1}\big\{\E\big[(1-q)^{|T_q|h}\big]=1-\rho\big\}.
\]
By definition of $\mathrm{hubs}(\rho,q)$ in~\eqref{eq:hubs}, the indicator is one precisely when $\mathrm{hubs}(\rho,q)\in\N$ and $q<1$. This proves the limit in~\eqref{eq:lemma-leading-term-2}. We next prove that $\lim_{r\uparrow \rho}C_{r, q}(h)-C_{\rho, q}(h)=0$. Similar to~\eqref{eq:limiting-constant-pr-1}, this corresponds to evaluating
\begin{equation}\label{eq:limiting-constant-pr-5}
\lim_{r\uparrow\rho}\frac{\alpha^h}{h!}\int_{y_1=0}^\infty\cdots\int_{y_h=0}^\infty\mathds{1}\{(y_1,\ldots, y_h)\in\CY_{r, q}\setminus\CY_{\rho, q}\}\cdot (y_1\cdot\ldots\cdot y_h)^{-(\alpha+1)}\rd y_1\cdot\ldots\cdot\rd y_h,
\end{equation}
Assume $r\in(\theta_q, \rho)$. If ${\bf{y}}\in\CY_{r, q}\setminus\CY_{\rho, q}$, then by definition of $\CY_{\rho, q}$ in~\eqref{eq:y-rho-set}, 
\[
1-\rho<\E\bigg[\prod_{x\in T_q, i\in[h]}\big(1-q(W_x^\sigma y_i\wedge 1)\big)\bigg]\le 1-r.
\]
Thus, the indicator function in~\eqref{eq:limiting-constant-pr-5} is monotone in $r$, and converges to $0$ for all $y$. This proves the  limit in~\eqref{eq:lemma-leading-term-3}.

The statement that $\rho\mapsto C_{\rho, q}$ is strictly decreasing follows from the fact that the integral in~\eqref{eq:limiting-constant-pr-1} is strictly positive.
\end{proof}

\subsection{Proofs of the main results}\label{sec:main-proofs}
Given the above lemmas, we can prove the main results of the paper.
\begin{proof}[Proof of \cref{thm:upper-tail}]
We start with the upper bound. Let $\phi, \delta>0$ be two small constants such that $\lceil\mathrm{hubs}(\rho-\delta, q)\rceil=\lceil\mathrm{hubs}(\rho, q)\rceil=:h$, which is possible by the continuity of $\mathrm{hubs}$ in~\eqref{eq:hubs}, see also \cref{lemma:properties-hubs}(ii). We distinguish two cases for the set of the weights of the vertices of at least $\phi n$ and apply \cref{prop:lower-bound} and \cref{lemma:leading-term}:
\[
\begin{aligned}
\Prob\big(|\CC_n^\sss{(1)}|>\rho n\big) &\le \Prob\big(\{|\CC_n^\sss{(1)}|>\rho n\}\cap\neg\{|\CV_n[\phi n, \infty)|=h, \CV_n[\phi n,\infty)\!\in\! n\cdot\CY_{\rho-\delta, q}(h)\}\big)\\
&\hspace{15pt}+\Prob\big(|\CV_n[\phi n, \infty)|=h, \CV_n[\phi n,\infty)\!\in\! n\cdot\CY_{\rho-\delta, q}(h)\big) 
 \\
&= O\big((n\Prob(W_1>n))^{h+1}\big) + (1+ o(1)) C_{\rho-\delta,q}(h)\cdot \big(n\Prob(W_1> n)\big)^{h} \\
&\sim C_{\rho-\delta,q}(h)\cdot \big(n\Prob(W_1> n)\big)^{h}.
\end{aligned}
\]
Since $\delta>0$ was arbitrary, it follows by~\eqref{eq:lemma-leading-term-3} in \cref{lemma:constant-limit} that
\[
\Prob\big(|\CC_n^\sss{(1)}|>\rho n\big)\lesssim C_{\rho,q}(h)\cdot\big(n\Prob(W_1> n)\big)^{\lceil\mathrm{hubs}(\rho, q)\rceil},
\]
proving the upper bound of both~\eqref{eq:upper-tail-1} and~\eqref{eq:upper-tail-2}.
For the lower bound we invoke \cref{prop:upper-bound} and \cref{lemma:leading-term} to obtain that 
\[
\begin{aligned}
\Prob&\big(|\CC_n^\sss{(1)}|>\rho n\big)\\\ &\ge \Prob\big(\{|\CC_n^\sss{(1)}|>\rho n\}\cap\{|\CV_n[\phi n, \infty)|=\lceil\mathrm{hubs}(\rho+\delta, q)\rceil, \CV_n[\phi n,\infty)\!\in\! n\cdot\CY_{\rho+\delta, q}(\lceil\mathrm{hubs}(\rho+\delta, q)\rceil)\}\big)\\
&\sim C_{\rho+\delta,q}(\lceil\mathrm{hubs}(\rho+\delta, q)\rceil)\cdot \big(n\Prob(W_1> n)\big)^{\lceil\mathrm{hubs}(\rho+\delta, q)\rceil}.
\end{aligned}
\]
If $\mathrm{hubs}(\rho, q)\notin\N$ or $q=1$, it follows that $\lceil\mathrm{hubs}(\rho, q)\rceil=\lceil\mathrm{hubs}(\rho-\delta, q)\rceil=\lceil\mathrm{hubs}(\rho+\delta, q)\rceil$ for any $\delta>0$ sufficiently small. Hence,~\eqref{eq:upper-tail-1} follows by sending $\delta\to0$ and invoking \eqref{eq:lemma-leading-term-2} in \cref{lemma:constant-limit}. If $\mathrm{hubs}(\rho, q)\in\N$, then $\lceil\mathrm{hubs}(\rho+\delta, q)\rceil=\mathrm{hubs}(\rho, q)+1$ for any $\delta$ sufficiently small and~\eqref{eq:upper-tail-2} follows.
\end{proof}

Next, we formally prove also the upper bound for the lower tail as stated in \cref{thm:lower-tail}.
\begin{proof}[Proof of \cref{thm:lower-tail}]
    The bound immediately follows from \cref{lemma:giant-trunc-lower}, since the largest component in $\CG_n$ has at least the size of the largest component in an induced subgraph $\CG_n[\underline w,\overline w)$.
\end{proof}
Next, we derive the large deviation principle in \cref{cor:ldp-giant}.
\begin{proof}[Proof of \cref{cor:ldp-giant}]
We start with the upper bound. 
If $\theta_q\in \bar B$, then $\mathrm{hubs}(\rho,q)=0=\inf_{\rho\in \bar B}I_q(\rho)$, and the upper bound is trivial. Assume $\theta_q\notin\bar B$, and assume $B$ is such that $b_-:=\max_{x<\theta_q}\{x\in \bar B\}$ and $b_+:=\min_{x>\theta_q}\{x\in \bar B\}$ exist. Then, for any $\eps>0$ 
\begin{equation}\label{eq:cor-ldp-pr1}
\Prob\big(|\CC_n|/n\in B\big)\le \Prob\big(|\CC_n|/n\le b_-\big)+
\Prob\big(|\CC_n|/n> b_+-\eps\big).
\end{equation}
Since $b_-$ and $b_+$ are strictly smaller (resp. larger) than $\theta_q$, the first term decays exponentially in $n$ by \cref{thm:lower-tail}. If $b_+<1$, the second term is regularly varying with index $I_q(b_+)$ by \cref{thm:upper-tail} if $\eps$ is sufficiently small so that $\lceil\mathrm{hubs}(b_+,q)\rceil$ and $\lceil\mathrm{hubs}(b_+-\eps,q)\rceil$ agree, which is possible since $\rho\mapsto\mathrm{hubs}(\rho, q)$ is continuous and increasing. Thus, the second term dominates the right-hand side in~\eqref{eq:cor-ldp-pr1} as $n\to\infty$, proving the upper bound if $b_-$ and $b_+$ both exist and $b_+<1$. If $b_+\ge 1$, we use that $\lim_{\rho\uparrow1}\mathrm{hubs}(\rho,q)=\infty$ by \cref{lemma:properties-hubs}, and the result follows by \cref{thm:upper-tail} and taking $\eps$ arbitrarily small. If $b_-$ does not exist, only the second term on the right-hand side remains and the upper bound follows by the same reasoning. If $b_+$ does not exist, only the first term remains, which decays exponentially in $n$, so its logarithm tends to $-\infty$ much faster than $\log n$.

We turn to the lower bound. If $\theta_q\in B^\circ$, the proof is again trivial. Assume $\theta_q\notin B^\circ$. Let  $b_+:=\inf_{x>\theta_q}\{x\in B^\circ\}$ as before. If $b_+$ does not exist or is at least $1$, in which case $I_q(\rho)=\infty$ for all $\rho\in B^\circ$, the lower bound is trivial. Assume that $b_+<1$. Since $B^\circ$ is an open set, there exists $\eps>0$ such that $(b_++\eps, b_++2\eps]\in B^\circ\cap(\rho,1)$. Then, 
\[
\Prob\big(|\CC_n|/n\in B\big)\ge 
\Prob\big(|\CC_n|/n> b_++\eps\big) - \Prob\big(|\CC_n|/n> b_++2\eps\big).
\]
By continuity and monotonicity of $\mathrm{hubs}$, we may assume that $\eps$ is so small that  \[\lceil\mathrm{hubs}(b_++\eps, q)\rceil=\lceil\mathrm{hubs}(b_++2\eps, q)\rceil=\inf_{\rho\in B^{\circ}}\lceil\mathrm{hubs}(\rho, q)\rceil,\] and that for $b\in\{b_++\eps,b_++2\eps\}$, $\mathrm{hubs}(b, q)\notin\N$. By \cref{thm:upper-tail},
\[
\Prob\big(|\CC_n|/n\in B\big)\ge 
(1+o(1))\big(C_{b_++\eps, q}- C_{b_++2\eps, q}\big)(n\Prob(W>n))^{\inf_{\rho\in B^{\circ}}\lceil\mathrm{hubs}(\rho, q)\rceil}.
\]
The constant factor is positive since $C_{b_++\eps, q}> C_{b_++2\eps, q}$ by~\cref{lemma:constant-limit}.  The probability on the right-hand side is regularly varying with index $-\alpha$. Thus, the lower bound follows.
\end{proof}

\section{The graph conditional on a large giant}\label{sec:remainder}
We prove the remaining corollaries from Section~\ref{sec:ldp-giant}. We first start with a more general lemma. 
\begin{lemma}\label{lemma:limit-hubs}
Consider an inhomogeneous scale-free random graph as in \cref{def:irg}. Let $\rho\in(\theta_q,1)$, set $h=\lceil\mathrm{hubs}{(\rho, q)}\rceil$, and assume $\mathrm{hubs}{(\rho, q)}\notin\N$ or $q=1$. Then there exists a constant $\phi>0$ such that, as $n\to\infty$,
\begin{equation}\label{eq:upper-tail-asymp}
\Prob\big(|\CC_n^\sss{(1)}|>\rho n\big)\sim \Prob\big(|\CV_n[\phi n,\infty)|=h, \CV[\phi n, \infty)\in n\cdot \CY_{\rho, q}(h)\big),
\end{equation}
Moreover, for a sequence $(\CA_n)_{n\ge 1}$ of events, 
    \[
    \Prob\big(\CA_n\mid |\CC_n^\sss{(1)}|>\rho n\big) = \Prob\big(\CA_n\mid |\CV[\phi n, \infty)|=h, \CV[\phi n, \infty)\in\CY_{\rho, q}\big) + o(1).
    \]
\begin{proof}
    The first statement follows by reasoning analogous to the beginning of the proof of \cref{thm:upper-tail} above: we combine Propositions~\ref{prop:lower-bound} and~\ref{prop:upper-bound}, Lemmas~\ref{lemma:leading-term} and~\ref{lemma:constant-limit},  and take the limit $r\!\to\!\rho$. We prove the second statement.         We start with a lower bound. Let $\delta>0$ be an arbitrarily small constant. We write out the conditional expectation and apply \cref{lemma:limit-hubs} to obtain for some small constant $\phi>0$
        \[
        \begin{aligned}
        \Prob\big(\CA_n\mid |\CC_n^\sss{(1)}|>\rho n\big) &\gtrsim \frac{\Prob\big(\CA_n\cap\{ |\CC_n^\sss{(1)}|>\rho n\}\cap \{|\CV[\phi n, \infty)|=h, \CV[\phi n, \infty)\in\CY_{\rho+\delta, q}\}\big)}{\Prob\big(|\CV[\phi n, \infty)|=h, \CV[\phi n, \infty)\in\CY_{\rho, q}\big)}\\
        &\ge \frac{\Prob\big(\CA_n\cap \{|\CV[\phi n, \infty)|=h, \CV[\phi n, \infty)\in\CY_{\rho+\delta, q}\}\big)}{\Prob\big(|\CV[\phi n, \infty)|=h, \CV[\phi n, \infty)\in\CY_{\rho, q}\big)} \\
        &\hspace{15pt}-\frac{\Prob\big( |\CC_n^\sss{(1)}|\le\rho n\}\cap \{|\CV[\phi n, \infty)|=h, \CV[\phi n, \infty)\in\CY_{\rho+\delta, q}\}\big)}{\Prob\big(|\CV[\phi n, \infty)|=h, \CV[\phi n, \infty)\in\CY_{\rho, q}\big)}.
        \end{aligned}
        \]
        The second term vanishes as $n\to\infty$ by \cref{prop:upper-bound} and \cref{lemma:leading-term}. For the first term we use that $\CY_{\rho+\delta, q}\subseteq\CY_{\rho, q}$ by its definition in~\eqref{eq:y-rho-set}. Thus, 
        \[\begin{aligned}
        \Prob\big(\CA_n\mid |\CC_n^\sss{(1)}|>\rho n\big) &\ge 
         \frac{\Prob\big(\CA_n\cap \{|\CV[\phi n, \infty)|=h, \CV[\phi n, \infty)\in\CY_{\rho, q}\}\big)}{\Prob\big(|\CV[\phi n, \infty)|=h, \CV[\phi n, \infty)\in\CY_{\rho, q}\big)} \\
        &\hspace{15pt}-\frac{\Prob\big(|\CV[\phi n, \infty)|=h, \CV[\phi n, \infty)\in\CY_{\rho, q}\setminus\CY_{\rho+\delta, q}\big)}{\Prob\big(|\CV[\phi n, \infty)|=h, \CV[\phi n, \infty)\in\CY_{\rho, q}\big)}-o(1).
        \end{aligned}
        \]
        The lower bound follows by rewriting the term on the first line as a conditional probability, and applying \cref{lemma:leading-term} and \cref{lemma:constant-limit} to the first term on the second line which vanishes as $\delta\to0$ under the assumption that $\mathrm{hubs}(\rho, q)\notin\N$.

        For the upper bound we argue similarly. Writing out the conditional expectation and distinguishing whether the hubs have weight in the set $\CY_{\rho-\delta, q}$ yields 
        \[
        \begin{aligned}
        \Prob\big(\CA_n\mid |\CC_n^\sss{(1)}|>\rho n\big)
        &\le \frac{\Prob\big(\CA_n\cap \{|\CV[\phi n, \infty)|=h, \CV[\phi n, \infty)\in\CY_{\rho-\delta, q}\}\big)}{\Prob\big(|\CV[\phi n, \infty)|=h, \CV[\phi n, \infty)\in\CY_{\rho, q}\big)} \\
        &\hspace{15pt}+\frac{\Prob\big(\{|\CC_n^\sss{(1)}|>\rho n\}\cap \neg\{|\CV[\phi n, \infty)|=h, \CV[\phi n, \infty)\in\CY_{\rho-\delta, q}\}\big)}{\Prob\big(|\CV[\phi n, \infty)|=h, \CV[\phi n, \infty)\in\CY_{\rho, q}\big)}.
        \end{aligned}
        \]
        If $\delta$ is sufficiently small, then $\lceil\mathrm{hubs}(\rho-\delta, q)\rceil=\lceil\mathrm{hubs}(\rho, q)\rceil$. By \cref{lemma:leading-term} and \cref{prop:lower-bound} the second term is of order $O(n\Prob(W_1>n))=o(1)$. For the numerator in the first term we also distinguish whether the hubs are in $\CY_{\rho,q}\subseteq\CY_{\rho-\delta, q}$. This yields
        \[
        \begin{aligned}
        \Prob\big(\CA_n\mid |\CC_n^\sss{(1)}|>\rho n\big)
        &\le \frac{\Prob\big(\CA_n\cap \{|\CV[\phi n, \infty)|=h, \CV[\phi n, \infty)\in\CY_{\rho, q}\}\big)}{\Prob\big(|\CV[\phi n, \infty)|=h, \CV[\phi n, \infty)\in\CY_{\rho, q}\big)} \\
        &\hspace{15pt}+
         \frac{\Prob\big(|\CV[\phi n, \infty)|=h, \CV[\phi n, \infty)\in\CY_{\rho-\delta, q}\setminus\CY_{\rho,q}\big)}{\Prob\big(|\CV[\phi n, \infty)|=h, \CV[\phi n, \infty)\in\CY_{\rho, q}\big)}
         +
        o(1).
        \end{aligned}
        \]
        We apply \cref{lemma:leading-term} and \cref{lemma:constant-limit} to the first term on the second line which vanishes as $\delta\to0$ under the assumption that $\mathrm{hubs}(\rho, q)\notin\N$. The term on the first line can be rewritten as a conditional probability. This finishes the proof.
\end{proof}
\end{lemma}

\begin{proof}[Proof of \cref{cor:cond-dist-weight}]
Let $h:=\lceil\mathrm{hubs}(\rho, q)\rceil$. Let 
\begin{equation}\label{eq:cond-dist-probs}
\Prob_{\CV_n}(\,\cdot\,):=\Prob\big(\,\cdot\mid |\CV[\phi n, \infty)|=h, \CV[\phi n, \infty)\in\CY_{\rho, q}\big), \qquad
\Prob_{Y}(\,\cdot\,):=\Prob\big(\,\cdot\mid \{\phi Y_i\}\in\CY_{\rho, q}\big).
\end{equation}
For the conditional distributional convergence, it suffices by \cref{lemma:limit-hubs} to show that for any Borel set $\CB\subseteq(0,\infty)^h$
\begin{equation}\label{eq:cond-dist-weight-pr1}
\begin{aligned}
\Prob_{\CV_n}\big(\CV_n[n^{1-\eps},\infty)\in n\CB, |\CV_n[n^{1-\eps},\infty)|= h\big) + \Prob_{\CV_n}&\big(|\CV_n[n^{1-\eps},\infty)|> h\big)\\
&\overset{n\to\infty}\longrightarrow\,  
\Prob_{Y}\big((\phi Y_i)_{i\le h}\in \CB\big).
\end{aligned}
\end{equation}
Below we show that the second term on the left-hand side vanishes as $n\to\infty$. We start with the first term. On the conditional measure, and the event $\{|\CV_n[n^{1-\eps},\infty)|= h\}$, it follows that $\CV_n[\phi n,\infty)=\CV_n[n^{1-\eps},\infty)$. Therefore, 
\begin{equation}\label{eq:cond-weight-pr1}
\Prob_{\CV_n}\big(\CV_n[n^{1-\eps},\infty)\in n\CB, |\CV_n[n^{1-\eps},\infty)|= h\big)=
\Prob_{\CV_n}\big(\CV_n[\phi n,\infty)\in n\CB\big).
\end{equation}
Now we use that all weights are independent with distribution $F_W$ from~\eqref{eq:weight-dist}. So for $y\ge 1$, with $Y_i$ from the statement of \cref{cor:cond-dist-weight},
    \begin{equation}\label{eq:large-weight-limit}
    \Prob\big(W_u \ge y\cdot \phi n\mid W_u\ge \phi n\big) = \frac{L(x\phi n)}{L(\phi n)}(y/\phi)^{-\alpha}\sim (y/\phi)^{-\alpha}=\Prob\big(\phi Y_1\ge y\big).
    \end{equation}
    By the definition of the conditional probability measure in~\eqref{eq:cond-dist-probs}, we obtain by~\eqref{eq:cond-weight-pr1}
    \[
    \Prob_{\CV_n}\big(\CV_n[n^{1-\eps},\infty)\in n\CB, |\CV_n[n^{1-\eps},\infty)|=h\big)
    \overset{n\to\infty}\longrightarrow
    \Prob_{Y}\big((\phi Y_i)_{i\le h}\in \CB\big).
    \]
We substitute this limit into~\eqref{eq:cond-dist-weight-pr1}. In the remainder of the proof we show that the second term in~\eqref{eq:cond-dist-weight-pr1}, $\Prob_{\CV_n}\big(|\CV_n[n^{1-\eps},\infty)|> h\big)$, converges to zero.
On the conditional probability measure, there are exactly $h$ vertices with weight at least $\phi n$. So the weights of the remaining $n-h$ vertices have distribution $\Prob(W_1>w\mid W_1<\phi n)$. Since all weights are independent, using the weight distribution in~\eqref{eq:weight-dist},
\[
\begin{aligned}
\Prob_{\CV_n}\big(|\CV_n[n^{1-\eps},\infty)|> h\big) &= \Prob\big(\mathrm{Bin}\big(n-h, \Prob(W_1>n^{1-\eps}\mid W_1<\phi n)\big)\ge 1\big)\\
&\le n\Prob(W_1>n^{1-\eps}\mid W_1<\phi n)
=O\big(L(n^{1-\eps})n^{1-(1-\eps)\alpha}\big)
=o(1)
\end{aligned}
\]
 for any $\eps>0$ sufficiently small as $\alpha>1$ by assumption in \cref{def:irg}.
\end{proof}
\begin{proof}[Proof of \cref{cor:cond-dist-comp}]
Recall that $R^\infty$ is the space of sequences of real numbers, metrized by $d_\infty({\bf{x}}, {\bf{y}})=\sum_i(|x_i-y_i|\wedge1 )2^{-i}$. By~\cite[Example 1.2]{billingsley2013convergence} each probability measure on $R^\infty$ is tight, and hence it is also relatively compact~\cite[Theorem 5.1]{billingsley2013convergence}. Thus, it satisfies by~\cite[Theorem 2.6]{billingsley2013convergence} to show convergence of the finite-dimensional distributions. Fix $\ell_\ast\in\N$.
    Let $(x_\ell)_{\ell\le\ell_\ast}\in[0,\infty)^{\ell_\ast}$ be a continuity point of the distribution of $\big(g_{\ell}\big((\phi Y_i)_{i\le \lceil\mathrm{hubs}(\rho, q)\rceil}\big), \ell\in[\ell_\ast]\big)$ conditionally on $ \big(\phi Y_i\big)_{i\le \lceil\mathrm{hubs}(\rho, q)\rceil}\in \CY_{\rho, q}$.  Abbreviate $h=\lceil\mathrm{hubs}(\rho, q)\rceil$. By \cref{lemma:limit-hubs}, 
    \begin{equation}\label{eq:comp-dist-pr1}
    \begin{aligned}
    \Prob\big(\forall \ell\le\ell_\ast&: N_{n,\ell}/n\le x_\ell\, \big|\, |\CC_n^\sss{(1)}|>\rho n \big) \\&= 
    \Prob\big(\forall \ell\le\ell_\ast: N_{n,\ell}/n\le x_\ell\,\big|\, |\CV_n[\phi n,\infty)|=h, \CV_n[\phi n, \infty)\in \CY_{\rho , q} \big)+ o(1).\end{aligned}
    \end{equation}
    We condition on $\CV[\phi n,\infty)=n{\bf{y}}^{(h)}$ for some ${\bf{y}}^{(h)}\in \CY_{\rho , q}$ and establish a lower and upper bound on the probability on the right-hand side. We start with a lower bound. Let $\delta>0$ be an arbitrary small constant. Then
    \[
    \begin{aligned}
        \Prob&\big(\forall \ell\le\ell_\ast: N_{n,\ell}/n\le x_\ell\, \big|\, \CV[\phi n,\infty)=n{\bf{y}}^{(h)} \big)\\ 
        &\ge 
        \Prob\big(\forall \ell\le\ell_\ast: N_{n,\ell}/n\le x_\ell, 
        \big|N_{n,\ell}/n -\tfrac1\ell\E\big[\ind{|T_q|=\ell}\bar P\big((W_x)_{x\in T_q}, {\bf{y}}^{(h)}\big)\big]\big|\le\delta
        \,\Big|\, \CV[\phi n,\infty)=n{\bf{y}}^{(h)} \big) \\
        &\ge 
        \Prob\big(\forall \ell\le\ell_\ast: \tfrac1\ell\E\big[\ind{|T_q|=\ell}\bar P\big((W_x)_{x\in T_q}, {\bf{y}}^{(h)}\big)\big]\le x_\ell-\delta
        \,\Big|\, \CV[\phi n,\infty)=n{\bf{y}}^{(h)} \big)\\
        &\hspace{15pt}-
        \Prob\big(\exists \ell\le\ell_\ast:  
        \big|N_{n,\ell}/n -\tfrac1\ell\E\big[\ind{|T_q|=\ell}\bar P\big((W_x)_{x\in T_q}, {\bf{y}}^{(h)}\big)\big]\big|>\delta
        \,\Big|\, \CV[\phi n,\infty)=n{\bf{y}}^{(h)} \big)
        .        
    \end{aligned}
    \]
    The negative term tends to $0$ by \cref{lemma:size-ell-hubs}. So, 
    \[
    \begin{aligned}
    \Prob&\big(\forall \ell\le\ell_\ast: N_{n,\ell}/n\le x_\ell\, \big|\, |\CC_n^\sss{(1)}|>\rho n \big) \\&\ge 
  \Prob\Big(\forall \ell\le\ell_\ast: \tfrac1\ell\E\big[\ind{|T_q|=\ell}\bar P\big((W_x)_{x\in T_q}, {\bf{y}}^{(h)}\big)\big]\le x_\ell-\delta
        \,\big|\, |\CV_n[\phi n,\infty)|=h, \CV_n[\phi n, \infty)\in \CY_{\rho , q} \Big)\\ &\hspace{15pt}- o(1).\end{aligned}
    \]
    Recall that by definition of $\bar P$ in~\eqref{eq:p-bar}, the expectation inside the probability corresponds to $g_\ell({\bf{y}}^\sss{(\ell)})$ defined in~\eqref{eq:g-comp-size}. By the weak convergence of the weights of the vertices in $\CV[\phi n,\infty)$, see~\eqref{eq:large-weight-limit}, we obtain 
    \[
    \begin{aligned}
    \Prob&\big(\forall \ell\le\ell_\ast: N_{n,\ell}/n\le x_\ell\, \big|\, |\CC_n^\sss{(1)}|>\rho n \big) \\&\ge 
  \Prob\Big(\forall \ell\le\ell_\ast: \tfrac1\ell\E\big[\ind{|T_q|=\ell}\bar P\big((W_x)_{x\in T_q}, (\phi Y_i)_{i\le \lceil\mathrm{hubs}(\rho,q)\rceil}\big)\big]\le x_\ell-\delta
        \,\big|\, (\phi Y_i)_{i\le \lceil\mathrm{hubs}(\rho,q)\rceil}\in \CY_{\rho ,q}) \Big) \\
        &\hspace{15pt}- o(1),\end{aligned}
    \]
    where $(Y_i)_{i\ge 1}$ are independent copies of $Y$ following distribution $\Prob(Y\ge y)=y^{-\alpha}$. Recall that $(x_\ell)_{\ell\le\ell_\ast}$ is a continuity point of the (finite-dimensional version of) the limiting distribution in~\cref{cor:cond-dist-comp}. Since $\delta>0$ was arbitrary, for each $\eps>0$ there exists $\delta>0$ such that 
    \[
    \begin{aligned}
      \Prob&\Big(\forall \ell\le\ell_\ast: \tfrac1\ell\E\big[\ind{|T_q|=\ell}\bar P\big((W_x)_{x\in T_q}, (\phi Y_i)_{i\le \lceil\mathrm{hubs}(\rho,q)\rceil}\big)\big]\le x_\ell-\delta
        \,\big|\, (\phi Y_i)_{i\le \lceil\mathrm{hubs}(\rho,q)\rceil}\in \CY_{\rho ,q}) \Big)\\
        &\ge
          \Prob\Big(\forall \ell\le\ell_\ast: \tfrac1\ell\E\big[\ind{|T_q|=\ell}\bar P\big((W_x)_{x\in T_q}, (\phi Y_i)_{i\le \lceil\mathrm{hubs}(\rho,q)\rceil}\big)\big]\le x_\ell
        \,\big|\, (\phi Y_i)_{i\le \lceil\mathrm{hubs}(\rho,q)\rceil}\in \CY_{\rho ,q}) \Big)\\
        &
        \hspace{15pt}-\eps/2.
        \end{aligned}
    \]
    Thus, for any $\eps>0$, when $n$ is sufficiently large 
        \[
    \begin{aligned}
    \Prob&\big(\forall \ell\le\ell_\ast: N_{n,\ell}/n\le x_\ell\, \big|\, |\CC_n^\sss{(1)}|>\rho n \big) \\&\hspace{-3pt}\ge 
            \Prob\Big(\forall \ell\le\ell_\ast: \tfrac1\ell\E\big[\ind{|T_q|=\ell}\bar P\big((W_x)_{x\in T_q}, (\phi Y_i)_{i\le \lceil\mathrm{hubs}(\rho,q)\rceil}\big)\big]\le x_\ell
        \,\big|\, (\phi Y_i)_{i\le \lceil\mathrm{hubs}(\rho,q)\rceil}\in \CY_{\rho ,q}) \Big)-\eps.\end{aligned}
    \]
    We leave it to the reader to prove an upper bound (almost analogously), so that weak convergence of the finite-dimensional distribution follows. By the reasoning above~\eqref{eq:comp-dist-pr1}, this suffices for the proof of the conditional component-size distribution in \cref{cor:cond-dist-comp}. The proof of the conditional distribution of $|\CC_n^\sss{(1)}|/n$ in~\eqref{eq:cond-size-giant} follows immediately from \cref{thm:upper-tail}.
\end{proof}

\appendix
\section{Postponed proofs}
We present the proofs of Lemmas~\ref{lemma:properties-hubs},~\ref{lemma:approx-non-connections}, and~\ref{lemma:size-ell-hubs}.
\begin{proof}[Proof of \cref{lemma:properties-hubs}]\label{proof:properties-hubs}
    The first statement follows from a rearrangement of~\eqref{eq:hubs}, and the fact that the generating function $H_{T_q}(z)$ is increasing in $z$. 

    We proceed to (ii). Continuity follows from continuity of the generating function $H_{T_q}(z)$. Assume $q$ is fixed, then by (i) it follows that $\rho\mapsto \mathrm{hubs}(\rho, q)$ is non-decreasing. 
    Assume $\rho$ is fixed. 
    The branching processes with different percolation parameters $q'<q$ can be coupled such that $\Prob\big(|T_q|\ge |T_{q'}| \big)=1$: $T_{q'}$ is obtained from $T_q$ by removing independently each edge and the entire subtree with probability $q'/q$. Hence, 
    \[
    \E\big[(1-q)^{|T_q|h}\big]\le \E\big[(1-q)^{|T_{q'}|h}\big]<\E\big[(1-q')^{|T_{q'}|h}\big].
    \]    
    As a result, 
    \[
    \inf\Big\{h'>0: \E\big[(1-q)^{|T_{q}|h'}\big]\le 1-\rho\Big\}< \inf\Big\{h'>0: \E\big[(1-q')^{|T_{q'}|h'}\big]\le 1-\rho\Big\},
    \]
    and (ii) follows by (i).
    
    We proceed to (iii) when $q\downarrow0$.  We analyze the generating function $\E[z^{|T_q|}\ind{|T_q|<\infty}]$ appearing in the definition of $\mathrm{hubs}$ in~\eqref{eq:hubs}. As $q\downarrow0$, the probability that the root of the branching process has degree 0, tends to 1. Thus, $\Prob\big(|T_q|=1\big)\to 1$ as well. Therefore, for any $z$ as $q\downarrow 0$,
    \[
    H_{T_q}(z)=\E[z^{|T_q|}\ind{|T_q|<\infty}] = \sum_{k=0}^{\infty}\Prob\big(|T_q|=k\big)z^k= (1+o(1))z.
    \]
    Inverting $H_{T_q}$ yields that also $H_{T_q}^{(-1)}(z)=(1+o(1))z$ as $q\downarrow 0$. Substituting this limit into~\eqref{eq:hubs} yields (iii) when $q\downarrow 0$. When $q$ is fixed and $\rho\uparrow1$, the other statement in (iii) follows immediately from part (iv), which we prove now. We analyze the generating function $\E[z^{|T_q|}\ind{|T_q|<\infty}]$ as $z\downarrow 0$. Since the total progeny of a branching process is at least $1$, 
    \[
    H_{T_q}(z)=\E[z^{|T_q|}\ind{|T_q|<\infty}] = z\Prob\big(|T_q|=1\big)(1+o(1)) = z\Prob\big(D_\varnothing=0\big)(1+o(1)),\qquad \text{as }z\downarrow0,
    \]
    where $D_\varnothing$ is the degree of the root of the branching process, see \cref{def:branching}. By definition, $D_\varnothing$ is a compound Poisson random variable with mean $qW_\varnothing\E[\kappa_\sigma(W_\varnothing, W)\mid W_\varnothing]$, where $W_\varnothing$ has distribution $F_W$. Therefore, $\Prob\big(D_\varnothing=0\big)=\E\big[\exp\big(-qW_\varnothing\E[\kappa_\sigma(W_\varnothing, W)\mid W_\varnothing]\big)\big]$. Thus, as $z\downarrow 0$
    \[
    H_{T_q}(z)=z\E\big[\exp\big(-qW_\varnothing\E[\kappa_\sigma(W_\varnothing, W)\mid W_\varnothing]\big)\big](1+o(1)).
    \]
    Inverting the formula, and substituting the limit into~\eqref{eq:hubs} yields the first limit in (iv). For the second limit, we use that $\kappa_1(W_\varnothing, W)=W_\varnothing W$ and that $W_\varnothing$ and $W_1$ are independent.
\end{proof}
\begin{proof}[Proof of \cref{lemma:approx-non-connections}]\phantomsection\label{proof:truncation-effect} We first give a probabilistic proof of the statement~\eqref{eq:non-conn-exp-lower}.
By \cref{def:type-w-comp}, \eqref{eq:non-conn-exp-lower} is equivalent to showing that there exist constants $\bar R, \ell_\ast$ such that
\begin{equation}\label{eq:non-conn-pr0}
\Prob\big(|T_q|\le \ell_\ast, \max_{x\in T_q}W_x\le \bar R\big)\ge \Prob\big(|T_q|<\infty\big)-\psi,
\end{equation}
which is equivalent to showing for some $\bar R, \ell_\ast$  that
\begin{equation}\label{eq:non-conn-pr-1}
\Prob\big(\ell_\ast<|T_q|<\infty\big)+\Prob\big(|T_q|\le \ell_\ast, \max_{x\in T_q}W_x> \bar R\big)\le \psi.
\end{equation}
The first term on the left-hand side tends to 0 as $\ell_\ast$ tends to 0. Let $\ell_0$ be such that for any $\ell_\ast\ge \ell_0$ the left-hand side is at most $\psi/2$. To bound the second term, we use that 
\begin{equation}\label{eq:non-conn-pr-2}
\Prob\big(|T_q|\le \ell_\ast, \max_{x\in T_q}W_x> \bar R\big)
\le 
\Prob\big(|T_q|\le \ell_\ast \mid \max_{x\in T_q}W_x> \bar R\big).
\end{equation}
We argue now that the right-hand side tends to zero as $\bar R\to\infty$. Let $x_{\ge \bar R}\in T_q$ be a vertex that has weight at least $\bar R$. Then the number of offspring of $x_{\ge \bar R}$ with weight at most $\bar R$, stochastically dominates a Poisson random variable with mean 
\[
q\bar R\int_{\underline w}^{\bar R} w^{\sigma}\rd F_W(w)=:C\bar R.
\]
So, the probability that $x_{\ge \bar R}$ has at least $(C/2)\bar R$ offspring tends to 1 as $\bar R\to\infty$. If $\bar R>2\ell_\ast/C$, the size of the total progeny $T_q$ is also at least $\ell_\ast$ on this event. Thus, there exists $R_0=R_0(\ell_0)$ such that for $\ell_\ast=\ell_0$ and $\bar R=R_0$ also the right-hand side in~\eqref{eq:non-conn-pr-2} is at most $\psi/2$. Thus, both terms in~\eqref{eq:non-conn-pr-1} are at most $\psi/2$ for these values $\ell_\ast=\ell_0$ and $R=R_0$. This proves~\eqref{eq:non-conn-exp-lower}. Since the left-hand side in~\eqref{eq:non-conn-pr0}, corresponding to the left-hand side in \eqref{eq:non-conn-exp-lower}, is increasing in $\ell_\ast$, and $\bar R$,~\eqref{eq:non-conn-exp-lower} also holds for any $\bar R\ge R_0, \ell_\ast\ge \ell_0$.

    To prove~\eqref{eq:non-conn-upper}, we need to show that we can choose $\eps, \ell_\ast\ge \ell_0, \bar R\ge R_0$ such that  
    \begin{equation}\label{eq:truncation-s}
    s(\eps, \ell_\ast, \bar R):=\sup_{{\bf{y}}>\phi{\bf{1}}}\bigg|\E\left[\bar P\big((W_x)_{x\in T_q}, {\bf{y}}\big)\right] - \sum_{\ell=1}^{\ell^\ast}\sum_{{\bf{w}}\in\mathrm{CT}_\ell(\eps, \bar R)} \bar P({\bf w}, {\bf y})\theta({\bf w},\eps)\bigg|
    \end{equation}
    can be made arbitrarily small. We truncate the expectation on the right-hand side in~\eqref{eq:truncation-s} by considering four events for the branching process. First define the constant 
    \[
    c:=\frac q2\int_{\underline w}^{2\underline w}w^\sigma \rd F_W(w)>0.\] Define the events
    \begin{align*}
        \CA_1(\bar R)&:=\Big\{\max_{x\in T_q}W_x >\bar R, \big|\{x\in T_q: W_x\in[\underline w, 2\underline w]\}\big|<c\bar R \Big\}, \\
        \CA_2(\bar R)&:=\Big\{\max_{x\in T_q}W_x >\bar R, \big|\{x\in T_q: W_x\in[\underline w, 2\underline w]\}\big|\ge c\bar R \Big\}, \\
        \CA_3(\ell_\ast, \bar R)&:=\Big\{\max_{x\in T_q}W_x \le \bar R, |T_q|> \ell_\ast\Big\}, \\
        \CA_4(\ell_\ast, \bar R)&:=\Big\{\max_{x\in T_q}W_x \le \bar R, |T_q|\le \ell_\ast\Big\}.
    \end{align*}
    Using that $\bar P\le 1$ by definition in~\eqref{eq:p-bar},  
    \begin{equation}\label{eq:truncation-4-events}
    \begin{aligned}
    \E\left[\bar P\big((W_x)_{x\in T_q}, {\bf{y}}\big)\right] &\le 
    \Prob\big(\CA_1(\bar R)\big) 
    +\E\left[\bar P\big((W_x)_{x\in T_q}, {\bf{y}}\big)\ind{\CA_2(\bar R)}\right] \\
    &\hspace{15pt}+ \Prob\big(\CA_3(\ell_\ast, \bar R)\big)
    + \E\left[\bar P\big((W_x)_{x\in T_q}, {\bf{y}}\big)\ind{\CA_4(\ell_\ast, \bar R)}\right].
    \end{aligned}
    \end{equation}
    The fourth term considers only the types of branching processes that are considered in the sums in~\eqref{eq:truncation-s}. We will show that it forms the main contribution to the right-hand side in~\eqref{eq:truncation-4-events}. We first analyze the three other terms. For the first term, we consider the offspring of the branching process of a vertex of weight at least $\bar R$, assuming that $\bar R\ge 2\underline w$. By \cref{def:branching}, the number of offspring in the interval $[\underline w, 2\underline w]$ dominates a Poisson random variable with mean $2c\bar R$. By concentration of Poisson random variables (\cref{lemma:chernoff}), we obtain that as $\bar R\to\infty$, 
    \[
    \Prob\big(\CA_1(\bar R)\big) \le \Prob\big(|\{x\in T_q: W_x\in[\underline w, 2\underline w]\}<c\bar R\,\big|\, \max_{x\in T_q}W_x >\bar R\big) \le \Prob\big(\mathrm{Poi}(2c\bar R)<c\bar R\big)=o(1).
    \]
    For the second term in~\eqref{eq:truncation-4-events}, we use that there are at least $c\bar R$ vertices with weight in $[\underline w, 2\underline w]$ on the event $\CA_2(\bar R)$. Let $c_1=\min\big(\underline w^\sigma, (2\underline w)^\sigma\big)$. By definition of $\bar P$ in~\eqref{eq:p-bar}, and that $y_i\ge \phi$ for all $i$, we obtain 
    \[
    \E\left[\bar P\big((W_x)_{x\in T_q}, {\bf{y}}\big)\ind{\CA_2(\bar R)}\right] \le \big(1-q(c_1\phi\wedge 1)\big)^{c\bar R} \longrightarrow 0,\qquad\text{as }\bar R\to\infty.
    \]
    We turn to the third term in~\eqref{eq:truncation-4-events} in which we assume that $\bar R$ is fixed so that the first two terms in~\eqref{eq:truncation-4-events} are at most $\psi/4$. We consider again the branching process. The number of offspring with weight at least $\bar R$ generated from a particle with weight at most $\bar R$, dominates a Poisson random variable with mean 
    \[
    q\min(\underline w^\sigma, \bar R^\sigma)\int_{\bar R}^\infty w\rd F_W(w)=:c_{\bar R}.
    \]
    The offspring of the first (at least) $\ell_\ast$ particles of the branching process are formed by independent Poisson point processes. Hence, the probability that none of them generates a particle of at least $\bar R$, decays exponentially in $\ell_\ast$, i.e., 
    \begin{equation}\nonumber
    \Prob\big(\CA_3(\ell_\ast, \bar R)\big)\le \exp\big(-\ell_\ast c_{\bar R}\big) \longrightarrow 0,\qquad\text{as }\ell_\ast\to\infty.
    \end{equation}
    This bounds the first three terms on the right-hand side in~\eqref{eq:truncation-4-events}. 
    Substituting the bounds on~\eqref{eq:truncation-4-events} into~\eqref{eq:truncation-s} yields that for all $\psi$ there exists a constant $R_1\ge R_0$, such that for all $R_2\ge R_1$ there exists a constant $\ell_{1}=\ell_1(R_2)\ge \ell_0$ such that for all $\ell_2\ge \ell_1$, 
    \begin{align}
    \inf_{\ell_\ast\ge 0, \bar R\ge 0} &s(\eps, \ell_\ast, \bar R)\label{eq:truncation-sup-three-lines}\\
        &
     \le 3\psi/4+\sup_{{\bf{y}}>\phi{\bf{1}}}  \Big|\E\big[\bar P\big((W_x)_{x\in T_q}, {\bf{y}}\big)\ind{|T_q|\le \ell_2}\big]-\sum_{\ell=1}^{\ell_2}\sum_{{\bf{w}}\in\mathrm{CT}_\ell(\eps, R_2)}  \theta({\bf w},\eps)\bar P({\bf{w}}, {\bf{y}}) \Big| 
        \nonumber\\
    &\le 3\psi/4+ \sup_{{\bf{y}}>\phi{\bf{1}}}\sum_{\ell=1}^{\ell_2}\sum_{{\bf{w}}\in\mathrm{CT}_\ell(\eps, R_2)}\hspace{-12pt}\theta({\bf w},\eps)\cdot\big|\E\big[\bar P({\bf w}, {\bf y})- \bar P\big((W_x)_{x\in T_q}, {\bf{y}}\big) \,\big|\, T_q\text{ has type }{(\bf{w}, \eps)}\big]\big| .\nonumber
    \end{align}
    Next, we will find a suitable upper bound on the expectations in the third line.
    Define 
    \begin{equation}\label{eq:truncation-single-type-1}
    \delta_\eps:= \sup_{{\bf{w}},{\bf{y}}}\bigg|\E\bigg[\frac{\bar P\big((W_x)_{x\in T_q}, {\bf{y}}\big)}{\bar P({\bf{w}}, {\bf{y}} )}\,\Big|\, T_q\text{ has type }{(\bf{w}, \eps)}\bigg]-1\bigg|.
    \end{equation}
    Here the supremum runs over all component types in $\mathrm{CT}_\ell(\eps, R_2)$ with $\ell\le\ell_2$, and ${\bf{y}}\ge\phi{\bf{1}}^{(h)}$.
    We show that $\delta_\eps$ tends to 0 as $\eps$ tends to 0. First assume that $\sigma\ge 0$.
    By construction of the types in \cref{def:type-w-comp}, the numerator is always smaller than the denominator, and for each $i\in [\ell]$ there exists $x\in T_q$ such that $w_i\le W_x< w_i+\eps$. Using the definition of $\bar P$ in~\eqref{eq:p-bar}, we obtain that for all ${\bf{y}}, {\bf{w}}$, 
    \begin{equation}\label{eq:truncation-single-type-2}
    \begin{aligned}
    1\ge \E\left[\frac{\bar P\big((W_x)_{x\in T_q}, {\bf{y}}\big)}{\bar P({\bf{w}}, {\bf{y}} )}\right]&\ge \prod_{i\in[\ell], j\in[h]}\frac{1-q\big((w_i+\eps)^\sigma y_j\wedge 1\big)}{1-q\big(w_i^\sigma y_j\wedge 1\big)} \\
    &\ge \left(\inf_{w\in[\underline w, R_2], y\ge \phi}\frac{1-q\big((w+\eps)^\sigma y\wedge 1\big)}{1-q\big(w y\wedge 1\big)}\right)^{\ell_2}.
    \end{aligned}
    \end{equation}
    When $\sigma<0$, these bounds hold in the opposite direction, replacing the infimum by a supremum.
    Since $\ell_2$ is a constant, it is elementary to verify that the right-hand side tends to 1 as $\eps\to 0$ for any $\sigma\in\R$. So $\delta_\eps$ tends to 0 as $\eps\to0$.  Thus, each expectation in the third line in~\eqref{eq:truncation-sup-three-lines} is at most $\delta_\eps\bar P({\bf{w}}, {\bf{y}})\le\delta_\eps$. This proves that 
    \begin{equation}\nonumber
    \begin{aligned}
    \inf_{\ell_\ast\ge 0, \bar R\ge 0} s(\eps, \ell_\ast, \bar R) &\le 3\psi/4 + \delta_\eps\sum_{\ell=1}^{\ell_2}\sum_{{\bf{w}}\in\mathrm{CT}_\ell(\eps, R_2)}\theta({\bf w},\eps)\\
    &= 3\psi/4+\delta_\eps\Prob\big(|T_q|\le \ell_2, \max_{x\in T_q}W_x\le R_2\big)\le 3\psi/4+\delta_\eps.
    \end{aligned}
    \end{equation}
    We conclude that for each $\psi>0$, there exist choices of $\ell_\ast, \bar R, \eps$ such that the right-hand side is at most $\psi$, so that~\eqref{eq:non-conn-upper} holds. The bound~\eqref{eq:non-conn-component} follows from the same reasoning as from~\eqref{eq:truncation-single-type-1} until the lines below ~\eqref{eq:truncation-single-type-2}. \end{proof}
\begin{proof}[Proof of \cref{lemma:size-ell-hubs}]\phantomsection\label{proof:concentration-components}
Let $\tilde\psi$ be a sufficiently small constant depending on $\psi$.
    We first claim that for any $\ell\in\N$, if $R$ is sufficiently large and $\eps$ sufficiently small that 
    \begin{equation}\label{eq:size-ell-hubs-pr1}
    \sup_{{\bf{y}}^{(h)}>\phi{\bf{1}}^{(h)}}\bigg|\E\Big[\ind{|T_q|=\ell}\bar P\big((W_x)_{x\in T_q}, {\bf{y}}^{(h)}\big)\Big]\,-\, \sum_{{\bf{w}}^\sss{(\ell)}\in\mathrm{CT}_\ell(\eps, R)}\bar P({\bf{w}}^\sss{(\ell)}, {\bf{y}}^{(h)})\theta({\bf{w}}^\sss{(\ell)}, \eps) \bigg|<\tilde\psi.
    \end{equation}
    We leave it to the reader to verify that this follows from the same reasoning as the proof of~\eqref{eq:non-conn-upper}, starting from~\eqref{eq:truncation-s}.  
    
    Moreover, let $\ell_\ast\ge \ell$ and $R$ be at least so large, and $\eps$ so small that we may apply \cref{lemma:approx-non-connections} with $\psi_{\ref{lemma:approx-non-connections}}=\tilde\psi$, and that $\Prob(|T_q|\le \ell_\ast)\ge \Prob(|T_q|<\infty)-\psi'=1-\theta_q-\tilde\psi$.
    By definition of $\theta({\bf{w}}^\sss{(\ell)}, \eps)$ in~\eqref{eq:branching-type-prob}, $\Prob(|T_q|=\ell)=\sum_{{\bf{w}}^\sss{(\ell)}\in\mathrm{CT}_\ell(\eps)}\theta({\bf{w}}^\sss{(\ell)}, \eps)$ for any $\eps>0$. Thus,
    \[
    \sum_{k\in[\ell_\ast]} \sum_{{\bf{w}}^\sss{(k)}\in\mathrm{CT}_k(\eps, R)}\hspace{-15pt}\theta({\bf{w}}^\sss{(k)},\eps)
    =
    \Prob\big(|T_q|\le\ell_\ast\big) - \sum_{k\in[\ell_\ast]} \sum_{{\bf{w}}^\sss{(k)}\in\mathrm{CT}_k(\eps)\setminus{CT}_k(\eps, R)}\hspace{-15pt}\theta({\bf{w}}^\sss{(k)},\eps).
    \]
    Since the sum over $\theta({\bf{w}}^\sss{(k)},\eps)$ for all $k<\infty$ and ${\bf{w}}^\sss{(k)}\in\mathrm{CT}_k(\eps)$ is at most one, if $R$ is sufficiently large, the double sum on the right-hand side is at most $\tilde\psi$. Thus, 
    \begin{equation}\label{eq:vertices-in-small}
        \sum_{k\in[\ell_\ast]} \sum_{{\bf{w}}^\sss{(k)}\in\mathrm{CT}_k(\eps, R)}\hspace{-15pt}\theta({\bf{w}}^\sss{(k)},\eps) \ge 1-\theta_q-2\tilde\psi. 
    \end{equation}
    
    We now adjust the proof of \cref{prop:upper-bound}. Let $\phi=\phi(\psi, C)$ be a sufficiently small constant. We define
    \begin{align}\CA_\mathrm{comp}&:=\Bigg\{\sum_{\substack{\ell\le \ell_\ast,\\ {\bf{w}}^\sss{(\ell)}\in \mathrm{CT}_\ell(\eps, R)}} \hspace{-15pt}\Big|\frac{\ell N_n({\bf{w}}^\sss{(\ell)}, \varepsilon, \phi n)}n - \theta({\bf{w}}^\sss{(\ell)}, \varepsilon)\Big|\le\tilde\psi\Bigg\}\cap\bigg\{\frac{|\CC_n^\sss{(1)}[\underline w, \phi n)}n\ge \theta_q-\tilde\psi\bigg\}\nonumber\\
    &\hspace{20pt}\cap 
    \Big\{\forall\ell\le \ell_\ast, {\bf{w}}^\sss{(\ell)}\in \mathrm{CT}_\ell(\eps, R): \big(|\ell N_n({\bf w}^\sss{(\ell)}, \eps, \phi n)/n\big) \big/ \theta({\bf w}^\sss{(\ell)}, \eps) - 1|\le \tilde \psi\Big\}
    ,\label{eq:a-comp-size-ell} \\
    \CA_\mathrm{hubs}&:=\big\{\forall v\in\CV_n[\phi n,\infty): v\sim \CC_n^\sss{(1)}[\underline w,\phi n)\big\}\label{eq:a-hubs-size-ell}
    \end{align}
    For each $\ell\le \ell_\ast$,
    we write $M_{n}({\bf{w}}^\sss{(\ell)},\eps)$ for the number of components of type ${\bf{w}}^\sss{(\ell)}$ in the induced subgraph $\CG_n[\underline w,\phi n)$ that are \emph{not} connected by an edge to the hubs in $\CG_n$. Let $h:=|\CV_n[\phi n,\infty)|$ denote the number of hubs, which is equal to $\lceil\mathrm{hubs}(r, q)\rceil$ by assumption. By the conditioning in~\eqref{eq:prop-lower} all rescaled weights in ${\bf{y}}^{(h)}=\{y_1,\ldots, y_h\}$ are at least $\phi$.
    Define 
    \begin{equation}
    \begin{aligned}
         \CA_\mathrm{conn}:=&\Bigg\{\begin{aligned}\forall \ell\!\in\![\ell_\ast], {\bf{w}}^\sss{(\ell)}\!\in\!\mathrm{CT}_\ell(\varepsilon,\bar R): M_{n}({\bf{w}}^\sss{(\ell)},\eps)&\ge (1-\tilde\psi)^3 \cdot(n/\ell)\cdot  \theta({\bf{w}}^\sss{(\ell)}, \eps)\\
    &\hspace{15pt}\cdot \bar P({\bf{w}}^\sss{(\ell)}, {\bf{y}}^{(h)}) \ind{\bar P({\bf{w}}^\sss{(\ell)}, {\bf{y}}^{(h)})\ge \tilde\psi}
    \end{aligned}\Bigg\}\\
    &\cap \Bigg\{\begin{aligned}\forall \ell\!\in\![\ell_\ast], {\bf{w}}^\sss{(\ell)}\!\in\!\mathrm{CT}_\ell(\varepsilon,\bar R):  M_{n}({\bf{w}}^\sss{(\ell)},\eps)&\le (1+\tilde\psi)^2\cdot\big(\tilde\psi+ \bar P({\bf w}^\sss{(\ell)}, {\bf{y}}^{(h)})\big)\\&\hspace{15pt}\cdot N_{n}({\bf{w}}^\sss{(\ell)}, \eps, \phi n)\end{aligned}\Bigg\}.
    \end{aligned}\label{eq:conn-size-ell}
    \end{equation}
    We next analyze the number of size-$\ell$ components on the intersection $\CA_\mathrm{comp}\cap\CA_\mathrm{hubs}\cap \CA_\mathrm{conn}$.  We first establish a lower bound. The number of remaining components of size $\ell$ is at least the number of components of size $\ell$ with type in $\mathrm{CT}_\ell(\eps, R)$ that does not connect by an edge to one of the hubs. Thus, we find for each $\ell\in[\ell_\ast]$,
    \[
    \begin{aligned}
    N_{n,\ell}\ge \hspace{-10pt}\sum_{{\bf{w}}^\sss{(\ell)}\in\mathrm{CT}_\ell(\eps, R)}\hspace{-10pt}M_{n}({\bf{w}}^\sss{(\ell)},\eps)&\ge \hspace{-10pt}\sum_{{\bf{w}}^\sss{(\ell)}\in\mathrm{CT}_\ell(\eps, R)}\hspace{-10pt}(1-\tilde\psi)^3 \cdot(n/\ell)\cdot  \theta({\bf{w}}^\sss{(\ell)}, \eps)\cdot \bar P({\bf{w}}^\sss{(\ell)}, {\bf{y}}^{(h)}) \ind{\bar P({\bf{w}}^\sss{(\ell)}, {\bf{y}}^{(h)})\ge\tilde\psi}\\
    &\ge (1-\tilde\psi)^3 \cdot(n/\ell)\hspace{-10pt}\sum_{{\bf{w}}^\sss{(\ell)}\in\mathrm{CT}_\ell(\eps, R)}\hspace{-10pt} \theta({\bf{w}}^\sss{(\ell)}, \eps)\cdot \bar P({\bf{w}}^\sss{(\ell)}, {\bf{y}}^{(h)}) \\
    &\hspace{15pt}- \tilde\psi(1-\tilde\psi)^3 \cdot(n/\ell)\hspace{-10pt}\sum_{{\bf{w}}^\sss{(\ell)}\in\mathrm{CT}_\ell(\eps, R)}\hspace{-10pt} \theta({\bf{w}}^\sss{(\ell)}, \eps)\\
    &\ge (1-\tilde\psi)^3 \cdot(n/\ell)\hspace{-10pt}\sum_{{\bf{w}}^\sss{(\ell)}\in\mathrm{CT}_\ell(\eps, R)}\hspace{-10pt} \theta({\bf{w}}^\sss{(\ell)}, \eps)\cdot \bar P({\bf{w}}^\sss{(\ell)}, {\bf{y}}^{(h)}) - \tilde\psi n.
    \end{aligned}
    \]
    In the last step we used that the sum over the probabilities on the third line is at most one. Applying the bound from~\eqref{eq:size-ell-hubs-pr1}, we obtain on $\CA_\mathrm{comp}\cap\CA_\mathrm{hubs}\cap \CA_\mathrm{conn}$ that
    \begin{equation}
        N_{n,\ell}/n\ge (1-\tilde\psi)^3 \cdot(1/\ell)\Big(\E\Big[\ind{|T_q|=\ell}\bar P\big((W_x)_{x\in T_q}, {\bf{y}}^{(h)}\big)\Big] - \tilde\psi\Big) - \tilde\psi.\nonumber
    \end{equation}
    We now subtract $(1/\ell)\Big(\E\Big[\ind{|T_q|=\ell}\bar P\big((W_x)_{x\in T_q}, {\bf{y}}^{(h)}\big)\Big]$ from both sides, and use that the expectation on the right-hand side is at most $1$ by definition of $\bar P$ in~\eqref{eq:p-bar}. As a result, for $\tilde\psi$ sufficiently small depending on $\psi$, we obtain that 
    \begin{equation}\label{eq:lower-size-ell}
    \begin{aligned}
    N_{n,\ell}&/n - (1/\ell)\E\Big[\ind{|T_q|=\ell}\bar P\big((W_x)_{x\in T_q}, {\bf{y}}^{(h)}\big)\Big]\\
    &\ge (-3\tilde\psi +3\tilde\psi^2-\tilde\psi^3 ) \cdot(1/\ell)\E\Big[\ind{|T_q|=\ell}\bar P\big((W_x)_{x\in T_q}, {\bf{y}}^{(h)}\big)\Big]  - (1-\tilde \psi)^3\tilde\psi/\ell-\tilde\psi \ge -\psi,
    \end{aligned}
    \end{equation}

    We next establish an upper bound on the number of size-$\ell$ components on $\CA_\mathrm{comp}\cap\CA_\mathrm{hubs}\cap \CA^\ge_\mathrm{conn}$. Since all hubs connect to a component of size at least $(\theta_q-\tilde\psi)n$, no component of size less than $\ell$ in $\CG_n[1,\phi n)$ is contained in a component of size exactly $\ell$ in $\CG_n$. Thus, the number of components of size-$\ell$ is at most the number of components with type in $\mathrm{CT}_\ell(\eps, R)$ that do not connect by an edge to one of the hubs, plus the number of components of size $\ell$ with type \emph{not} in $\mathrm{CT}_\ell(\eps, R)$, i.e., 
    \begin{equation}\label{eq:n-n-ell-upper-1}
    N_{n,\ell}\le \sum_{{\bf{w}}^\sss{(\ell)}\in\mathrm{CT}_\ell(\eps, R)}\hspace{-10pt} M_{n}({\bf{w}}^\sss{(\ell)},\eps) + \sum_{{\bf{w}}^\sss{(\ell)}\in\mathrm{CT}_\ell(\eps)\setminus\mathrm{CT}_\ell(\eps, R)}\hspace{-10pt}N_{n}({\bf{w}}^\sss{(\ell)},\eps, \phi n).
    \end{equation}
    We first bound the second term from above, multiplied by a factor $\ell$ to count vertices. The number of vertices in components with type in $\mathrm{CT}_\ell(\eps)\setminus\mathrm{CT}_\ell$ is at most the number of components that is not in the largest component, and not in a component with type in $\mathrm{CT}_k(\eps, R)$ for some $k\in\N$. We obtain
    \[
    \sum_{{\bf{w}}^\sss{(\ell)}\in\mathrm{CT}_\ell(\eps)\setminus\mathrm{CT}_\ell(\eps, R)}\hspace{-15pt}\ell N_{n}({\bf{w}}^\sss{(\ell)},\eps, \phi n) \le n-|\CC_n^\sss{(1)}|-\sum_{k\in[\ell_\ast]} \sum_{{\bf{w}}^\sss{(K)}\in\mathrm{CT}_k(\eps, R)}\hspace{-15pt}k N_{n}({\bf{w}}^\sss{(k)},\eps, \phi n).
    \]
    We use the bounds from the definition of $\CA_\mathrm{comp}$ in~\eqref{eq:a-comp-size-ell}, yielding 
    \[
    \sum_{{\bf{w}}^\sss{(\ell)}\in\mathrm{CT}_\ell(\eps)\setminus\mathrm{CT}_\ell(\eps, R)}\hspace{-15pt}\ell N_{n}({\bf{w}}^\sss{(\ell)},\eps, \phi n) \le  (1-\theta_q+\tilde\psi)n + \tilde\psi n - \sum_{k\in[\ell_\ast]}\sum_{{\bf{w}}^\sss{(K)}\in\mathrm{CT}_k(\eps, R)}\hspace{-15pt}k \theta({\bf{w}}^\sss{(k)},\eps)n. 
    \]
    By~\eqref{eq:vertices-in-small}, the double sum on the right-hand side is at least $(1-\theta_q-2\tilde\psi)n$, so that the right-hand side is in total at most $4\tilde\psi n$. We substitute this bound in~\eqref{eq:n-n-ell-upper-1}, and use the upper bound on $M_n({\bf w}^\sss{(\ell)}, \eps)$ from~\eqref{eq:conn-size-ell} to obtain 
    \[
    \begin{aligned}
    N_{n,\ell}&\le 4\tilde\psi n + \sum_{{\bf{w}}^\sss{(\ell)}\in\mathrm{CT}_\ell(\eps, R)}\hspace{-10pt}(1+\tilde\psi)^2\cdot\big(\tilde\psi+ \bar P({\bf w}^\sss{(\ell)}, {\bf{y}}^{(h)})\big)\cdot N_{n}({\bf{w}}^\sss{(\ell)}, \eps, \phi n)\\
    &\le 4\tilde\psi n + (1+\tilde\psi)^2\tilde\psi n + 
    (1+\tilde\psi)^2\sum_{{\bf{w}}^\sss{(\ell)}\in\mathrm{CT}_\ell(\eps, R)}\hspace{-10pt}\bar P({\bf w}^\sss{(\ell)}, {\bf{y}}^{(h)})\cdot N_{n}({\bf{w}}^\sss{(\ell)}, \eps, \phi n)
    \end{aligned}
    \]
    as the total number of components is at most $n$. Without loss of generality, we may assume that $\tilde\psi$ is at most one. Next, we use the upper bound on $N_n({\bf w}^\sss{(\ell)}, \eps, \phi n)$ by $\CA_\mathrm{comp}$ in~\eqref{eq:a-comp-size-ell}, which yields, 
    \[
    \begin{aligned}
    N_{n,\ell}&\le 8\tilde\psi n + (1+\tilde\psi)^2 (n/\ell)\sum_{{\bf{w}}^\sss{(\ell)}\in\mathrm{CT}_\ell(\eps, R)}\hspace{-10pt}\bar P({\bf w}^\sss{(\ell)}, {\bf{y}}^{(h)})\cdot\big(\tilde\psi \theta({\bf w}^\sss{(\ell)}, \eps) + \theta({\bf w}^\sss{(\ell)}, \eps)\big).
    \end{aligned}
    \]
    Since $\bar P$ is at most $1$ by~\eqref{eq:p-bar}, and the sum over all $\theta({\bf w}^\sss{(\ell)}, \eps)$ is at most one, we obtain by~\eqref{eq:size-ell-hubs-pr1}
    \[
    \begin{aligned}
    N_{n,\ell}&\le 8\tilde\psi n + (1+\tilde\psi)^2\tilde\psi + (1+\tilde\psi)^2 (n/\ell)\sum_{{\bf{w}}^\sss{(\ell)}\in\mathrm{CT}_\ell(\eps, R)}\hspace{-10pt}\bar P({\bf w}^\sss{(\ell)}, {\bf{y}}^{(h)})\cdot\theta({\bf w}^\sss{(\ell)}, \eps)\\
    &\le 12\tilde \psi n + (n/\ell)(1+\tilde\psi)^2\Big(\E\Big[\ind{|T_q|=\ell}\bar P\big((W_x)_{x\in T_q}, {\bf{y}}^{(h)}\big)\Big] + \tilde\psi\Big).
    \end{aligned}
    \]
    Thus, there exists a constant $c>0$ such that 
    \[
    \le N_{n,\ell}/n - (1/\ell)\E\Big[\ind{|T_q|=\ell}\bar P\big((W_x)_{x\in T_q}, {\bf{y}}^{(h)}\big)\Big]\le c\psi
    \]
    We combine this upper bound with the lower bound in~\eqref{eq:lower-size-ell}, and obtain on the event $\CA_\mathrm{comp}\cap\CA_\mathrm{hubs}\cap \CA_\mathrm{conn}$
    \[
    \Big|N_{n,\ell}/n - (1/\ell)\E\Big[\ind{|T_q|=\ell}\bar P\big((W_x)_{x\in T_q}, {\bf{y}}^{(h)}\big)\Big]\Big|\le \max(c, 1)\tilde \psi.
    \]
    Thus, for proving \cref{lemma:size-ell-hubs}, it suffices to show that 
    \begin{equation}
        \Prob\big(\CA_\mathrm{comp}\cap\CA_\mathrm{hubs}\cap \CA_\mathrm{conn}\mid \CV_n[\phi n,\infty)=n{\bf y}^\sss{(h)}\big)=1-o(n^{-C}).
    \end{equation}
    We argue similar as around~\eqref{eq:cond-meas-lb}. We define the conditional probability measure
    We condition on the graph $\CG_n[\underline w,\bar R)$ satisfying $\CA^\ge_\mathrm{comp}$ and the realization of $\CV_n[\phi n,\infty)=n{\bf{y}}^{(h)}$ satisfying $\CA_\mathrm{comp}$. We abbreviate 
    \begin{equation}\label{eq:cond-meas-comp}
        \Prob_{{\bf{y}}, \CG}\big(\,\cdot\,\big):=\Prob\big(\, \cdot \mid \CG_n[\underline w,\bar R), \CV_n[\phi n,\infty)=n{\bf{y}}^{(h)}, \CA_\mathrm{comp}\big).
    \end{equation}    
    Then, 
    \begin{equation}\nonumber
    \begin{aligned}
    \Prob\big(\CA_\mathrm{comp}\cap\CA_\mathrm{hubs}&\cap \CA_\mathrm{conn}\mid \CV_n[\phi n,\infty)=n{\bf y}^\sss{(h)}\big)\\&=
    \E\Big[\ind{\CA^\ge_\mathrm{comp}}\Prob_{{\bf{y}}, \CG}\big(\CA_\mathrm{hubs}\cap\CA_\mathrm{conn}\big)\,\big|\,\CV_n[\phi n, \infty)=n{\bf{y}}^\sss{(\ell)}\Big].
    \end{aligned}
    \end{equation}
    To bound the conditional probability, the same reasoning as below~\eqref{eq:cond-meas-lb} applies, using for bounding the probability of the event $\CA_\mathrm{conn}$ also the reasoning from~\eqref{eq:upper-upper-pr2}. As a result, the conditional probability is at least $1-\exp(-c'n)$ for some $c'>0.$ Recalling $\CA_\mathrm{comp}$ from~\eqref{eq:a-comp-size-ell}, this leaves to show that 
    \[
    \begin{aligned}
        \Prob&\Bigg(\sum_{\substack{\ell\le \ell_\ast,\\ {\bf{w}}^\sss{(\ell)}\in \mathrm{CT}_\ell(\eps, R)}} \hspace{-15pt}\Big|\frac{\ell N_n({\bf{w}}^\sss{(\ell)}, \varepsilon, \phi n)}n - \theta({\bf{w}}^\sss{(\ell)}, \varepsilon)\Big|>\tilde\psi\Bigg) \\&\hspace{20pt}+ \Prob\bigg(\frac{|\CC_n^\sss{(1)}[\underline w, \phi n)}n\ge \theta_q-\tilde\psi\bigg)\nonumber\\
    &\hspace{20pt}+ 
    \Prob\Big(\forall\ell\le \ell_\ast, {\bf{w}}^\sss{(\ell)}\in \mathrm{CT}_\ell(\eps, R): \big(|\ell N_n({\bf w}^\sss{(\ell)}, \eps, \phi n)/n\big) \big/ \theta({\bf w}^\sss{(\ell)}, \eps) - 1|> \tilde \psi\Big)=o(n^{-C}).
    \end{aligned}
    \]
    The first two terms are $o(n^{-C})$ by Lemmas~\ref{lemma:ell-trunc-n}--\ref{lemma:giant-trunc-lower}. The third term is also of order $o(n^{-C})$ by \cref{lemma:ell-trunc-n} by choosing $\psi$ dependent on the finitely many $\theta({\bf w}^\sss{(\ell)}, \eps)$. This finishes the proof.
\end{proof}

\section{Preliminaries}
The following lemma is a straightforward application of the Chernoff bound. 
\begin{lemma}[Concentration bounds]\label{lemma:chernoff}
Let $X$ be Poisson or Binomial with mean $\mu$. Then, for every $\delta>0$ there exists a constant $c_\delta>0$ such that $\Prob(|X-\mu|> \delta \mu) \leq \re^{-c_\delta \mu}$.
\end{lemma}

The next result provides a useful estimate for sums of truncated heavy-tailed random variables; it is a version of Lemma 3 of \cite{resnick1999} adapted to our setting. 

\begin{lemma}
\label{lem-truncated}
Let $(W_i)_{i\in[n]}$ be iid random variables with regularly varying distribution as in Definition~\ref{def:irg}. For every $\psi, C, R>0$ there exists $\phi_0>0$ such that for all $\phi\in(0,\phi_0)$,
\begin{equation}
    \Prob\bigg(\sum_{i=1}^nW_i\ind{W\in[R,\phi n)}> 2\E[W\ind{W\ge R}] n \bigg) = o(n^{-C}). 
\end{equation}
\end{lemma}
\begin{proof}
    Let $(W'_i)_{i\ge 1}$ be iid copies of $W$ conditionally on $W\!\ge\! R$. Let $n_\delta\!:=\! (1\!+\!\delta)n\Prob(W\!\ge\! R)$ for  $\delta\in(0,1)$. Conditionally on $|\CV_n[R, \infty)|\le n_\delta$, $\sum_{i\in[n]}W_i\ind{W_i\ge R}$ is stochastically dominated by $\sum_{i\in[n_\delta]}W'_i\ind{W_i'<\phi n}$. Since $|\CV_n[R, \infty)|\sim\mathrm{Bin}(n, \Prob(W\ge R))$ and $R$ is a constant, we obtain by a Chernoff bound, 
\[
\begin{aligned}
\Prob\bigg(\sum_{i=1}^nW_i&\ind{W_i\in[R,\phi n)}> 2\E[W\ind{W\ge R}] n \bigg) \\
&\le 
\Prob\bigg(\sum_{i=1}^{n_\delta}W'_i\ind{W_i'<\phi n}> 2\E[W\ind{W\ge R}] n \bigg) + \Prob\Big(|\CV_n[R, \phi n)|> n_\delta\Big) \\
&=\Prob\bigg(\sum_{i=1}^{n_\delta}W'_i\ind{W_i'<\phi n}> \frac{2}{1+\delta}\E[W'] n_\delta \bigg) + \exp(-\Omega(n)).
\end{aligned}
\]
Let $(W^\sss{(\phi)}_i)_{i\ge 1}$ be iid copies of $W'$ conditional upon $W'< \phi n$. Then $\sum_{i\in[n_\delta]}W'_i\ind{W'_i<\phi n}$ is stochastically dominated by $\sum_{i\in[n_\delta]}W^\sss{(\phi)}_i$. As a result, 
    by Lemma 3 in~\cite{resnick1999}, also the first term is of order $o(n_\delta^{-C})=o(n^{-C})$ for any $\phi$ sufficiently small.
\end{proof}

Our final auxiliary result is a bound for sums of independent Bernoulli random variables.

\begin{lemma}[{\cite[Theorem A.1.4]{AlonSpencer}}]
\label{lem-bernoulli}
Let $B_i,i\geq 1,$ be a sequence of independent Bernoulli random variables with $p_i=\Prob(B_i=1)=1-\Prob(B_i=0)$.
Set $\mu_n = \sum_{i=1}^n p_i$. For every $b>0$ we have 
\begin{equation}
    \Prob\bigg(\sum_{i=1}^n B_i > (1+b) \mu_n\bigg) \leq \re^{-\mu_n I_B(b)}, \hspace{1cm} \Prob\bigg(\sum_{i=1}^n B_i < (1-b) \mu_n\bigg)\leq \re^{-\mu_n I_B(-b)},
\end{equation}
with $I_B(b) = (1+b) \log (1+b)-b$.
\end{lemma}

\end{document}